\documentclass[11pt,a4paper]{article}
\usepackage{amsmath,amssymb,amsthm,amscd,mathrsfs}
\usepackage{indentfirst}
\usepackage[top=25mm, bottom=30mm, left=30mm, right=30mm]{geometry}
\usepackage{color}
\usepackage{comment}

\usepackage{hyperref}

\usepackage{enumerate}

\newtheorem{df}{\bf Definition}[section]
\newtheorem{thm}[df]{\bf Theorem}
\newtheorem{lem}[df]{\bf Lemma}
\newtheorem{cor}[df]{\bf Corollary}
\newtheorem{prop}[df]{\bf Proposition}
\newtheorem{rem}[df]{\bf Remark}

\newtheorem*{claim}{\bf Claim}

\newtheorem{thmA}{\bf Theorem}
\newtheorem{propA}[thmA]{\bf Proposition}

\newtheorem{corA}[thmA]{\bf Corollary}

\newcommand{\R}{\mathbb{R}}
\newcommand{\C}{\mathbb{C}}
\newcommand{\Z}{\mathbb{Z}}

\newcommand{\N}{\mathbb{N}}

\newcommand{\B}{\mathbb{B}}

\newcommand{\M}{\mathbb{M}}

\newcommand{\W}{\mathrm{W}}

\newcommand{\ri}{\mathrm{i}}
\newcommand{\Ad}{\operatorname{Ad}}
\newcommand{\id}{\text{\rm id}}
\newcommand{\Int}{\operatorname{Int}}
\newcommand{\Aut}{\operatorname{Aut}}
\newcommand{\Out}{\operatorname{Out}}
\newcommand{\Tr}{\mathord{\text{\rm Tr}}}

\newcommand{\ovt}{\mathbin{\overline{\otimes}}}

\newcommand{\otm}{\otimes_{\rm min}}
\newcommand{\ota}{\otimes_{\rm alg}}
\newcommand{\rtm}{\rtimes_{\rm red}}

\renewcommand{\Re}{\operatorname{Re}}

\newcommand{\alg}{\mathrm{alg}}
\newcommand{\op}{\mathrm{op}}
\title{\bf A boundary approach to rigidity of $q$-Araki--Woods factors}
\author{Yusuke Isono\thanks{Research Institute for Mathematical Sciences, Kyoto University, 606-8502, Kyoto, Japan}}
\date{}

\begin{document}

\maketitle

\renewcommand{\thefootnote}{}
\footnote[0]{E-mail: \texttt{isono@kurims.kyoto-u.ac.jp}}
\footnote[0]{2020 \textit{Mathematics Subject Classification.} 46L10, 46L36, 46L54.}
\footnote[0]{\textit{Key words}: $q$-Araki--Woods factors; Type $\rm III$ factors; strong solidity; Connes' $\tau$-invariant.}
\footnote[0]{YI is supported by JSPS KAKENHI Grant Number 24K06759.}
\renewcommand{\thefootnote}{\arabic{footnote}}

\begin{abstract}
	In this paper, we study structural and rigidity properties of $q$-Araki--Woods factors arising from strongly continuous orthogonal representations of $\R$.
	We first introduce an AO-type boundary property that is formulated in terms of length projections on the $q$-Fock spaces, and prove that all $q$-Araki--Woods factors satisfy this condition. 
	We next characterize uniform convergence of the associated malleable
	deformations. We then use these two ingredients as a substitute for weak containment
	in the coarse bimodule and deduce several rigidity theorems in the
	framework of Popa's deformation/rigidity theory.
	As a main consequence, we prove that every $q$-Araki--Woods factor is stably strongly solid. 
	We also prove a normalizer rigidity theorem for inclusions arising from direct sum representations. Finally, we compute Connes' $\tau$-invariant, and construct a continuum of pairwise non-isomorphic type $\mathrm{III}_1$ $q$-Araki--Woods factors with the same $\tau$-invariant.
\end{abstract}

\section{Introduction}

The $q$-Araki--Woods von Neumann algebras were introduced by Hiai
\cite{Hi03} as $q$-deformations of Shlyakhtenko's free Araki--Woods
factors \cite{Sh96}. They are associated with a strongly continuous
orthogonal representation $U:\mathbb R\to\mathcal O(H_{\mathbb R})$
and a parameter $-1<q<1$. 
When $q=0$, one recovers the free Araki--Woods factors, while for the
trivial representation $U$ one recovers the $q$-Gaussian von Neumann
algebras \cite{BS90,BKS96}.

In this paper, we study $q$-Araki--Woods factors from the viewpoint of
Popa's \emph{deformation/rigidity theory}. In this framework, a basic
difference between the free and the $q$-deformed settings appears for
the bimodules associated with the \emph{malleable deformations}.
Let us briefly recall the relevant point. 
A typical malleable deformation involves an inclusion
$M\subset \widetilde M$ together with a one-parameter family
$(\alpha_t)_{t\in\mathbb R}$ of automorphisms of $\widetilde M$.
The orthogonal complement
$L^2(\widetilde M)\ominus L^2(M)$ is naturally an $M$-$M$ bimodule, and
one particularly useful condition is
\[
L^2(\widetilde M)\ominus L^2(M)
\prec
L^2(M)\otimes L^2(M)
\qquad\text{as $M$-$M$ bimodules},
\]
where the right-hand side is the coarse bimodule.
Such weak containment allows nonamenability to be exploited through
spectral-gap arguments, which can then be combined with malleability to
obtain uniform control of the deformation. This mechanism plays a basic
role in Popa's deformation/rigidity theory; see
\cite{Po03,IPP05,Po06,OP07}.

This coarse behavior is particularly transparent for free products.
Indeed, for an inclusion $M\subset M*N$,
the reduced-word decomposition of $L^2(M*N)$ shows that
\[
L^2(M*N)\ominus L^2(M)
\prec
L^2(M)\otimes L^2(M)
\]
as $M$-$M$ bimodules. For free Araki--Woods factors, the required coarse
weak containment follows directly from freeness.

For $q\neq0$, this direct free-product argument is no longer available.
Nevertheless, when the underlying Hilbert space is finite-dimensional,
one can recover sufficient coarse behavior by isolating suitable
subbimodules \cite{Av11,Wi20}. The resulting coarse-bimodule estimates
are dimension-dependent and do not directly extend to the
infinite-dimensional setting. This distinction is not merely technical:
in the infinite-dimensional $q$-Gaussian case, the factors for
$q\neq0$ are not isomorphic to the corresponding free case $q=0$
\cite{Ca22}; see also \cite{DY26}. Related non-isomorphism results
covering a broad class of $q$-Araki--Woods factors were also obtained in
\cite{DT25}.

Significant progress in the infinite-dimensional setting was made by
Ding and Peterson \cite{DP23}, who proved strong solidity for arbitrary
$q$-Gaussian algebras. Their proof crucially uses finite-dimensional
approximation together with the finite-dimensional coarse-bimodule
estimates described above.
Using related finite-dimensional input, Ding and Tan \cite{DT25}
established strong solidity for $q$-Araki--Woods factors
associated with almost periodic representations. This approximation
scheme, however, does not apply when the underlying orthogonal
representation has a nonzero weakly mixing part, since the weakly
mixing component cannot be approximated by finite-dimensional invariant
subrepresentations.

One of the main purposes of this paper is to overcome this obstruction
and develop a deformation/rigidity method that applies in the presence
of a nonzero weakly mixing part.
Our boundary approach replaces the missing coarse-bimodule input by
two complementary ingredients. The first is an \emph{AO-type boundary 
property} adapted to the length decomposition of the $q$-Fock space.
After passing to an appropriate quotient, this property retains a
suitable form of coarse-bimodule control. The second is a characterization
of \emph{uniform convergence of the associated malleable deformation}.
In the non-uniform case, this characterization yields completely positive
maps and states that vanish on the ideal defining the quotient and hence
factor through it. The two ingredients can then be combined in
deformation/rigidity arguments, allowing us to treat arbitrary
$q$-Araki--Woods factors.

\medskip

We introduce the notation needed to state our main results.
Fix $-1<q<1$. Let $\Gamma_q(H_{\mathbb R},U)''$ be the $q$-Araki--Woods factor associated with a strongly continuous orthogonal representation
$U:\mathbb R\to\mathcal O(H_{\mathbb R})$.
It contains a canonical $C^*$-subalgebra and we denote it by $A_q$; see Subsection \ref{$q$-Araki--Woods factors}.
Let $G$ be a second countable locally compact group and let
$\pi:G\to\mathcal O(H_{\mathbb R})$ be a strongly continuous orthogonal representation commuting with $U$.
We denote by
\[
\sigma^\pi:G\curvearrowright \Gamma_q(H_{\mathbb R},U)''
\]
the associated Bogoliubov action. We consider a $\sigma$-weakly dense inclusion 
\[
A_q(G):=A_q\rtimes_{\sigma^\pi,\mathrm{red}}G
\subset 
\Gamma_q(H_{\mathbb R},U)''\rtimes_{\sigma^\pi} G
=:M_q(G)\subset \mathbb B(\mathcal F_q(H)\otimes L^2(G)),
\]
where $H=H_\R\otimes_\R\C$ and $\mathcal F_q(H)$ is the $q$-Fock space.
We denote by
\[
\mathcal K_\ell(G)\subset \mathbb B(\mathcal F_q(H)\otimes L^2(G))
\]
the canonical hereditary $C^*$-subalgebra generated by the length
projections associated with the $q$-Fock space. We regard
$\mathcal K_\ell(G)$ as the analogue of the compact operators adapted
to the length decomposition. When $G$ is trivial and $H$ is
finite-dimensional, $\mathcal K_\ell(G)$ is precisely the algebra of
compact operators.

Our first main theorem is an AO-type property for the inclusion
$A_q(G)\subset M_q(G)$. 
Below, $\mathcal M(\mathcal K_\ell(G))\subset \B(\mathcal F_q(H)\otimes L^2(G))$ is the multiplier algebra, and $J_G$ is the modular conjugation of $M_q(G)$ on $\mathcal F_q(H)\otimes L^2(G)$.

\begin{thmA}\label{thmA}
	Assume that $G$ is amenable. Then the multiplication map
	\[
	\nu:
	A_q(G)\otimes_{\mathrm{alg}}J_GA_q(G)J_G
	\longrightarrow
	\mathcal M(\mathcal K_\ell(G))/\mathcal K_\ell(G),
	\]
	given by
	\[
	\nu(a\otimes J_GbJ_G)
	=
	aJ_GbJ_G+\mathcal K_\ell(G),
	\]
	is a well-defined $*$-homomorphism and is continuous with respect to
	the minimal tensor norm.
\end{thmA}

We in fact prove a stronger statement involving approximate u.c.p.\
lifts;
see Theorem~\ref{thm-length-AO} and Corollary~\ref{cor-stable-AO}.

When $G$ is trivial and $H_{\mathbb R}$ is finite-dimensional,
$\mathcal K_\ell(G)$ is the algebra of compact operators, and
the conclusion of Theorem~\ref{thmA} is precisely the multiplication-map
condition appearing in Ozawa's property (AO) \cite{Oz03}.
Related results for finite-dimensional $q$-Gaussian algebras were
obtained in \cite{CIW19}, whereas in the infinite-dimensional setting
this multiplication-map condition with respect to $A_q$ fails for
$q\neq0$ \cite{BCKW22}.
Thus, the use of the length-compact ideal $\mathcal K_\ell(G)$ is
essential in the infinite-dimensional setting.
For related von Neumann algebraic variants of property (AO), see
\cite{Ca22,DP23,DT25}. A related formulation also appears in the very
recent preprint \cite{DY26}.

Let us briefly comment on the proof of Theorem~\ref{thmA}.
When $U$ is almost periodic, we combine estimates from
\cite{Wi20,DT25} with an argument inspired by the deformation approach
to biexactness in \cite{DP23}.
For an arbitrary orthogonal representation, we approximate it by almost
periodic representations and use second quantization
\cite{HR10,Wa17} to pass from the almost periodic case to the general one,
in a way compatible with the $q$-Fock space construction and the
length-compact ideals.
This approximation argument is the main new ingredient in our
proof of Theorem~\ref{thmA}.
Since it may be useful independently, we collect the relevant results in
Appendix~\ref{Almost periodic approximation}.

The crossed-product formulation of Theorem~\ref{thmA} is needed for
our later applications to continuous cores. Moreover, the approximation
used in its proof naturally leads to intermediate crossed products
arising from commuting Bogoliubov actions, which need not themselves
be continuous cores. This is why we formulate the theorem for general
amenable crossed products rather than only for continuous cores.

\medskip

We next turn to the second ingredient in our approach. Let
$\widetilde M_q(G)$ be the doubled $q$-Araki--Woods crossed product and
let $(\alpha_\theta)_{\theta\in\mathbb R}$ be the associated malleable
deformation; see Subsection~\ref{Malleable deformations}. Our next main
result characterizes when this deformation converges uniformly on a
finite von Neumann subalgebra.

We use Popa's intertwining notation $\preceq_M$; see
Subsection~\ref{Further preliminaries on von Neumann algebras}. 
We say that an inclusion $A\subset N$ of von Neumann algebras is \emph{with expectation} if there exists a faithful normal conditional expectation from $N$ onto $A$.

\begin{thmA}\label{thmB}
	Let $p\in M_q(G)$ be a nonzero projection and let
	$A\subset pM_q(G)p$ be a finite von Neumann subalgebra with expectation.
	Assume that either
	\begin{enumerate}
		\item[$\rm (a)$] $G=\{e\}$; or
		\item[$\rm (b)$]
		$A'\cap p\widetilde M_q(G)p\subset pM_q(G)p$ or
		$(A'\cap pM_q(G)p)'\cap p\widetilde M_q(G)p
		\subset pM_q(G)p$.
	\end{enumerate}
	Then the following conditions are equivalent:
	\begin{enumerate}
	\item $\alpha_\theta$ converges to $\id$ uniformly on the unit ball of $A$ as $\theta\to0$;	
	\item for every nonzero projection
		$r\in A'\cap pM_q(G)p$, we have
		$Ar\preceq_{M_q(G)}LG$.
	\end{enumerate}
	Here $LG\subset M_q(G)$ is the group von Neumann algebra of $G$.
\end{thmA}

When $G$ is trivial, condition~(2) is equivalent to $A$ being completely
atomic. Thus, in this case, Theorem~\ref{thmB} extends the corresponding
result of Ding and Peterson \cite{DP23}. The proof of case~(a) follows a
strategy similar to \cite{DP23}, using Haagerup-type completely positive maps arising
from second quantization \cite{HR10,Wa17}. In case~(b), the additional relative commutant
assumptions allow us to apply a standard deformation/intertwining
argument familiar from free-product settings.
The technical result needed for this argument is established in Appendix~\ref{Intertwining between orthogonal summands in $q$-Araki--Woods factors}.

\subsubsection*{Application 1: Stable strong solidity}

Recall that a diffuse von Neumann algebra $M$ is \emph{strongly solid}
if the normalizer of every diffuse amenable von Neumann subalgebra
$A\subset M$ with expectation generates an amenable von Neumann algebra.
Strong solidity is a strong indecomposability property and is one of the
central rigidity phenomena in Popa's deformation/rigidity theory.

Strong solidity was introduced by Ozawa and Popa, who proved it for free
group factors \cite{OP07}. For free Araki--Woods factors, partial results
were obtained in \cite{HR10}, and strong solidity was established in full
generality in \cite{BHV15}. For $q$-Gaussian factors, the finite-dimensional
case was treated in \cite{Av11,Wi20} and the general case in \cite{DP23}.
More recently, Ding and Tan \cite{DT25} proved strong solidity for
$q$-Araki--Woods factors associated with almost periodic representations.

We prove the stronger stable version introduced in \cite{BHV15}. For an
inclusion $A\subset M$, set
\[
s\mathcal N_M(A)
:=
\{x\in M\mid xAx^*\subset A,\ x^*Ax\subset A\}.
\]
We say that $M$ is \emph{stably strongly solid} if
$s\mathcal N_M(A)''$ is amenable for every diffuse amenable von Neumann
subalgebra $A\subset M$ with expectation. The unitary elements of
$s\mathcal N_M(A)$ are precisely the usual normalizer
$\mathcal N_M(A)$, so stable strong solidity implies strong solidity.
This property also behaves well under corners and amplifications.

\begin{thmA}\label{thmC}
	Every $q$-Araki--Woods factor is stably strongly solid. 
	
	More precisely, $p\Gamma_q(H_{\mathbb R},U_t)''p$ is stably strongly solid 
	for every $(H_{\mathbb R},U_t)$ and for every nonzero projection
	$p\in \Gamma_q(H_{\mathbb R},U_t)''$.
\end{thmA}

Thus the theorem covers arbitrary orthogonal representations, including those
with a weakly mixing part. 
The proof combines the weak compactness argument
of \cite{OP07,BHV15} with Theorem~\ref{thmA} and case (a) of Theorem \ref{thmB}, following ideas in \cite{PV12,Is12}.

\subsubsection*{Application 2: Connes' $\tau$-invariant and non-isomorphism}

Our second application concerns the classification of type $\mathrm{III}_1$
$q$-Araki--Woods factors.
Recall that, for a full type $\mathrm{III}_1$ factor $M$, Connes'
\emph{$\tau$-invariant} $\tau(M)$ \cite{Co74} is the weakest topology on
$\mathbb R$ for which
\[
\R\ni t\longmapsto [\sigma_t^\varphi]\in
\Out(M):=\Aut(M)/\Int(M)
\]
is continuous, where $\varphi$ is a faithful normal state.
This topology does not depend on the choice of $\varphi$.
For free Araki--Woods factors, the $\tau$-invariant was computed in
\cite{Sh97,Va05}.

We first record the corresponding computation for
$q$-Araki--Woods factors. Fullness in the type $\mathrm{III}_1$ case follows
from \cite{HI20,KSW23,KW24}.

\begin{propA}\label{propD}
	Let $(H_{\mathbb R},U_t)$ be such that
	$\Gamma_q(H_{\mathbb R},U_t)''$ is a type $\mathrm{III}_1$ factor.
	Then $\Gamma_q(H_{\mathbb R},U_t)''$ is full, and its
	$\tau$-invariant is the weakest topology on $\mathbb R$ for which
	\[
	\R\ni t\longmapsto U_t\in \mathcal O(H_\R)
	\]
	is continuous. Here $\mathcal O(H_\R)$ is equipped with the strong operator topology.
\end{propA}

Our main classification result in this direction is the following.
It shows, in particular, that even within the class of $q$-Araki--Woods
factors, the $\tau$-invariant is far from being a complete isomorphism invariant.

\begin{thmA}\label{thmE}
	There exists a continuum of strongly continuous orthogonal representations
	$(H_{\mathbb R}^{(i)},U_t^{(i)})$, $i\in I$, on separable Hilbert spaces such that the factors
	\[
	\bigl\{\Gamma_q(H_{\mathbb R}^{(i)},U_t^{(i)})''\bigr\}_{i\in I}
	\]
	are all of type $\mathrm{III}_1$, have the same $\tau$-invariant, and are
	pairwise non-isomorphic.
\end{thmA}

This is a $q$-deformed analogue of the corresponding result for free
Araki--Woods factors in \cite{HSV16}, and the proof follows the strategy
developed there. The key input from \cite{HSV16} is a criterion relating
Popa's intertwining between the canonical copies of $L(\mathbb R)$ in
continuous cores to the corresponding faithful normal states. The crossed-product
forms of Theorems~\ref{thmA} and~\ref{thmB} provide the ingredients needed to
apply this strategy in the present setting.

\subsubsection*{Application 3: Normalizer rigidity for orthogonal sums}

Our third application concerns inclusions arising from orthogonal direct sum
decompositions. The basic orthogonality argument behind this phenomenon goes
back to Ricard \cite{Ri04} in the tracial $q$-Gaussian setting. Closely related
finite-dimensional arguments appear in \cite{SW16,BM16}.

We obtain a normalizer rigidity result without any finite-dimensionality
assumption. Let
\[
(H_{\mathbb R},U_t)
=
(H_{\mathbb R}^{(1)},U_t^{(1)})
\oplus
(H_{\mathbb R}^{(2)},U_t^{(2)})
\]
be an orthogonal direct sum of invariant subspaces, and put
\[
M=\Gamma_q(H_{\mathbb R},U_t)'',
\qquad
M_i=\Gamma_q(H_{\mathbb R}^{(i)},U_t^{(i)})'',
\quad i=1,2.
\]

\begin{thmA}\label{thmF}
	Let $p\in M_1$ be a nonzero projection and let
	$A\subset pM_1p$ be a diffuse von Neumann subalgebra with expectation.
	Then we have
	$s\mathcal N_{pMp}(A)''
	\subset pM_1p$.
\end{thmA}

In the proof, we introduce an analogue of weak mixing for such inclusions
and combine it with Theorem~\ref{thmB} and the orthogonality arising from
the direct sum decomposition.

As a consequence, we also obtain an outerness result for Bogoliubov
automorphisms.

\begin{corA}\label{corG}
	Let $U:\mathbb R\to\mathcal O(H_{\mathbb R})$ be a strongly continuous
	orthogonal representation with
	$\dim H_{\mathbb R}\geq 2$, and let
	$V\in\mathcal O(H_{\mathbb R})$ commute with $(U_t)_{t\in\mathbb R}$.
	If $V\neq 1$, then the Bogoliubov automorphism $\sigma^V$ of
	$\Gamma_q(H_{\mathbb R},U_t)''$ is outer.
	
	In particular, if $\pi:G\to\mathcal O(H_{\mathbb R})$ is a faithful
	orthogonal representation of a discrete group $G$ commuting with
	$(U_t)_{t\in\mathbb R}$, then the crossed product
	$\Gamma_q(H_{\mathbb R},U_t)''
	\rtimes_{\sigma^\pi}G$ is a factor.
\end{corA}

For $q=0$, related outerness and factoriality results are known in the free
Araki--Woods setting \cite{HS09,HT19}.

\bigskip
\noindent{\bf Declaration on the use of AI.}
In an earlier stage of this project, the author used ChatGPT Plus (GPT-5.5) to discuss possible approaches to several proofs.
The main ideas behind the AO-type boundary property,
the characterization of uniform convergence of the malleable
deformation, and the strong solidity argument originated with
the author.
In a later stage of the project, after the author switched to ChatGPT Pro (GPT-5.6 Sol and GPT-6 Astra), an initial proof of Theorem~\ref{thmF} was provided by ChatGPT.
The author subsequently developed the underlying ideas through
further discussion with ChatGPT and reformulated the argument into the present form.
ChatGPT also provided substantial assistance in carrying out and checking
the lengthy calculations in Appendix~\ref{Intertwining between orthogonal summands in $q$-Araki--Woods factors}.
It was also used for language editing and to assist with the writing of the manuscript,
including a few refinements and simplifications of the exposition and proofs.
The author takes full responsibility for all mathematical content.

\tableofcontents

\section{Preliminaries}
\label{Preliminaries}

Throughout this paper, we fix a parameter $q\in(-1,1)$ and do not repeat this assumption in the statements of individual results. The real Hilbert spaces considered throughout the paper are not assumed to be separable.
All von Neumann algebras are assumed to be
$\sigma$-finite, and all locally compact groups are assumed to be second
countable. All inner products on complex Hilbert spaces are linear in the second
variable.

For a faithful normal semifinite weight $\varphi$ on a von Neumann algebra
$M$, we write
\[
\|x\|_\varphi:=\varphi(x^*x)^{1/2},
\qquad \text{for }x\in \mathfrak n_\varphi
:=
\{a\in M:\varphi(a^*a)<\infty\},
\]
and denote by $L^2(M,\varphi)$ the corresponding GNS Hilbert space.
We freely use standard facts and notation from \emph{Tomita--Takesaki modular
theory}, for which we refer to \cite{Ta03}. 
We write $(M)_r:=\{x\in M:\|x\|\leq r\}$ for $r>0$.

We write $\mathbb N_0:=\mathbb N\cup\{0\}$; indices associated with the degree decomposition of the Fock space will always range over $\mathbb N_0$.

\subsection{$q$-Fock spaces}
\label{$q$-Fock spaces}

We recall the construction of $q$-Fock spaces \cite{BS90,BS93,BKS96}.
Let $H$ be a complex Hilbert space.
For $n\geq1$, the $q$-symmetrizer on the usual Hilbert tensor product
$H^{\otimes n}$ is given by
\[
P_q^{(n)}(\xi_1\otimes\cdots\otimes\xi_n)
=
\sum_{\sigma\in S_n}
q^{\operatorname{inv}(\sigma)}
\xi_{\sigma(1)}\otimes\cdots\otimes\xi_{\sigma(n)},
\]
where
\[
\operatorname{inv}(\sigma)
=
\#\{(i,j):1\leq i<j\leq n,\ \sigma(i)>\sigma(j)\}.
\]
The operator $P_q^{(n)}$ is positive and boundedly invertible.
We denote by $H_q^{\otimes n}$ the space $H^{\otimes n}$ equipped with
the equivalent inner product
\[
\langle\xi,\eta\rangle_q
=
\langle\xi,P_q^{(n)}\eta\rangle.
\]
The \emph{$q$-Fock space} over $H$ is
\[
\mathcal F_q(H)
=
\bigoplus_{n=0}^\infty H_q^{\otimes n},
\qquad
H_q^{\otimes0}=\mathbb C\Omega,
\]
where $\Omega$ is the vacuum vector. 
For $n,N\in \N_0$, let
\[
P_n:\mathcal F_q(H)\to H_q^{\otimes n},
\qquad
P_{\leq N}:=\sum_{n=0}^N P_n,
\qquad
P_{>N}:=1-P_{\leq N}
\]
be the corresponding orthogonal projections. 
We refer to these as \emph{length projections}.

For each $\xi\in H$, the left and right creation operators
$\ell_q(\xi),r_q(\xi)\in\mathbb B(\mathcal F_q(H))$ are determined by
$\ell_q(\xi)\Omega=r_q(\xi)\Omega=\xi$ and
\[
\begin{aligned}
	\ell_q(\xi)(\xi_1\otimes\cdots\otimes\xi_n)
	&=\xi\otimes\xi_1\otimes\cdots\otimes\xi_n,\\
	r_q(\xi)(\xi_1\otimes\cdots\otimes\xi_n)
	&=\xi_1\otimes\cdots\otimes\xi_n\otimes\xi
\end{aligned}
\]
for $n\geq1$ and $\xi_1,\ldots,\xi_n\in H$.
They satisfy
\[
\|\ell_q(\xi)\|=\|r_q(\xi)\|
\leq
\begin{cases}
	\|\xi\|, & -1<q\leq0,\\[2mm]
	\dfrac{\|\xi\|}{\sqrt{1-q}}, & 0\leq q<1.
\end{cases}
\]
Their adjoints annihilate $\Omega$ and satisfy
\[
\begin{aligned}
	\ell_q(\xi)^*(\xi_1\otimes\cdots\otimes\xi_n)
	&=
	\sum_{k=1}^n q^{k-1}\langle\xi,\xi_k\rangle\,
	\xi_1\otimes\cdots\widehat{\xi_k}\cdots\otimes\xi_n,\\
	r_q(\xi)^*(\xi_1\otimes\cdots\otimes\xi_n)
	&=
	\sum_{k=1}^n q^{n-k}\langle\xi,\xi_k\rangle\,
	\xi_1\otimes\cdots\widehat{\xi_k}\cdots\otimes\xi_n.
\end{aligned}
\]
In particular, for all $\xi,\eta\in H$, the following $q$-commutation relation holds:
\[
\ell_q(\xi)^*\ell_q(\eta)
-q\ell_q(\eta)\ell_q(\xi)^*
=
\langle\xi,\eta\rangle1.
\]

For a contraction $T:H\to K$ between complex Hilbert spaces,
its Fock-space second quantization is the contraction
$\mathcal F_q(T):\mathcal F_q(H)\to\mathcal F_q(K)$ 
determined by $\mathcal F_q(T)\Omega=\Omega$ and 
\[
\mathcal F_q(T)(\xi_1\otimes\cdots\otimes\xi_n)
=
T\xi_1\otimes\cdots\otimes T\xi_n
\]
for $n\geq1$ and $\xi_1,\ldots,\xi_n\in H$.

\subsection{Length-compact operators}
\label{Length-compact operators}

Let $G$ be a second countable locally compact group.
We define the $C^*$-algebra of \emph{length-compact operators} by
\[
\mathcal K_\ell(G)
:=
\overline{
	\bigcup_{N\geq0}
	(P_{\leq N}\otimes 1)\mathbb B(\mathcal F_q(H)\otimes L^2(G))(P_{\leq N}\otimes 1)
}^{\|\cdot\|}.
\]
This is the hereditary $C^*$-subalgebra generated by $(P_{\leq N}\otimes 1)_{N\geq0}$,
and $(P_{\leq N}\otimes 1)_{N\geq0}$ is an approximate identity.
Equivalently, for
$T\in\mathbb B(\mathcal F_q(H)\otimes L^2(G))$,
\[
T\in\mathcal K_\ell(G)
\quad\Longleftrightarrow\quad
\lim_{N\to\infty} \|T-T(P_{\leq N}\otimes 1)\|=\lim_{N\to\infty} \|T-(P_{\leq N}\otimes 1)T\|=0.
\]
The terminology reflects the fact that $P_{\leq N}$ cuts off tensors of length greater than $N$.
Thus, operators in $\mathcal K_\ell(G)$ vanish in norm as the tensor length tends to infinity, both on the left and on the right.
Note that if $G$ is trivial  and $H$ is finite dimensional, $\mathcal K_\ell(G)$ coincides with the set of all compact operators.

We identify the multiplier algebra $\mathcal M(\mathcal K_\ell(G))$ with its idealizer in
$\mathbb B(\mathcal F_q(H)\otimes L^2(G))$.
When $G$ is trivial, we simply write $\mathcal K_\ell:=\mathcal K_\ell(G)$. 

We say that $T\in\mathbb B(\mathcal F_q(H)\otimes L^2(G))$ is \emph{finite-band}  with respect to the length decomposition if there exists $d>0$ such that
\[
(P_n\otimes1)T(P_m\otimes1)=0
\qquad\text{whenever }n,m\in \N_0,\ |n-m|>d.
\]
Every finite-band operator belongs to $\mathcal M(\mathcal K_\ell (G))$,
since multiplication on either side by $P_{\le N}\otimes1$ gives
an operator in $\mathcal K_\ell (G)$. 
For example, all left and right creation and annihilation operators, tensored with the identity on $L^2(G)$, are finite-band.

\subsection{$q$-Araki--Woods factors}
\label{$q$-Araki--Woods factors}

We recall the construction of $q$-Araki--Woods factors
\cite{Sh96,Hi03}.
Let $H_{\mathbb R}$ be a real Hilbert space and let
$U:\mathbb R\to\mathcal O(H_{\mathbb R})$ be a strongly continuous
orthogonal representation.
Write $H=H_{\mathbb R}\otimes_{\mathbb R}\mathbb C $
for its complexification and let $I$ be the canonical conjugation on $H$.
We use the same notation for the complexification of $U$.
Let $A$ be the positive nonsingular self-adjoint operator on $H$ such that
$U_t=A^{it}$ for all $t\in\mathbb R$.
Since $IU_tI=U_t$ and $I$ is antilinear, we have $IAI=A^{-1}$.
Put
\[
j:=\frac{\sqrt2}{\sqrt{1+A^{-1}}} \in \B(H).
\]
We have $jH=\operatorname{dom}(A^{-1/2})$, and the densely defined closed
antilinear involution $S:=IA^{-1/2}$ with $\operatorname{dom}(S)=jH$ 
satisfies
\[
S(j(\xi)+\ri j(\eta))=j(\xi)-\ri j(\eta)
\qquad(\xi,\eta\in H_{\mathbb R}).
\]
For each $\zeta\in jH$, put
\[
W(\zeta):=\ell_q(\zeta)+\ell_q(S\zeta)^*.
\]
For $\zeta\in jH$ and $\xi,\eta\in H_\R$, we have
\[
W(\zeta)^* = W(S\zeta),\qquad 
W(j(\xi)+\ri j(\eta))
=
W(j(\xi))+\ri W(j(\eta)).
\]
In particular, $W(j(\xi))$ is self-adjoint for every
$\xi\in H_{\mathbb R}$.
The \emph{$q$-Araki--Woods von Neumann algebra} associated with $(H_{\mathbb R},U_t)$ is defined by
\[
\Gamma_q(H_{\mathbb R},U_t)''
:=
\{W(\zeta):\zeta\in jH\}''.
\]
We often use the following notation:
\[
A_q
:=
C^*\{1,W(\zeta):\zeta\in jH\}
\subset\mathbb B(\mathcal F_q(H)),
\qquad
M_q:=A_q''=\Gamma_q(H_{\mathbb R},U_t)''.
\]

\medskip
\noindent
\textbf{Modular objects.}
The vacuum vector $\Omega$ is cyclic and separating for $M_q$.
The corresponding faithful normal state
\[
\varphi_q(x):=\langle\Omega,x\Omega\rangle
\qquad(x\in M_q)
\]
is called the \emph{$q$-quasi-free state}. We always identify $L^2(M_q,\varphi_q)=\mathcal F_q(H)$.
Then the modular conjugation $J_q$ and modular operator $\Delta_q$ associated with $(M_q,\varphi_q)$ satisfy
\[
J_q\Omega=\Omega,
\qquad
J_q(\xi_1\otimes\cdots\otimes\xi_n)
=
I\xi_n\otimes\cdots\otimes I\xi_1
\]
for $n\geq1$ and $\xi_1,\ldots,\xi_n\in H$, and
\[
\Delta_q^{it}=\mathcal F_q(U_{-t}),
\qquad
\sigma_t^{\varphi_q}(W(j(\xi)))
=
W(j(U_{-t}\xi))
\quad(t\in\mathbb R,\ \xi\in H).
\]
In particular, $J_q$ depends only on $I$, not on the representation $U$.

For $\eta\in IjH$, we define the right Wick element by
\[
W_r(\eta)
:=
r_q(\eta)+r_q(S^*\eta)^*
=
J_qW(I\eta)J_q
\in M_q',
\]
where $S^*=A^{-1/2}I=IA^{1/2}$ and $\operatorname{dom}(S^*)=IjH$. 

We say that $\zeta\in H$ is \emph{analytic} if it is entire analytic
for $U$, equivalently, if
$\zeta\in\bigcap_{n\in\mathbb Z}\operatorname{dom}(A^n)$.
In this case, we have 
\[
J_qW(A^{1/2}\zeta)^*J_q
=
J_qW(I\zeta)J_q
=
W_r(\zeta).
\]

\medskip
\noindent
\textbf{Wick polynomials.}
For every $n\geq 1$ and $\zeta=\zeta_1\otimes \cdots \otimes \zeta_n \in (jH)^{\otimes_\alg n}$, there is a unique $W(\zeta)\in A_q$ such that
\[
W(\zeta)\Omega=\zeta.
\]
We also put $W(\Omega)=1$. These operators are called \emph{Wick words}.
We have 
\[
W(\zeta_1\otimes\cdots\otimes\zeta_n)^*
=
W(S\zeta_n\otimes\cdots\otimes S\zeta_1).
\]
The set of \emph{Wick polynomials} $A_{q,\mathrm{alg}}$ is defined as the linear span of all Wick words.
Then $A_{q,\mathrm{alg}}$ is a norm-dense unital $*$-subalgebra of $A_q$. 
The right Wick words $W_r(\zeta)$ are defined similarly,
with $jH$ replaced by $IjH$.

\medskip
\noindent
\textbf{Bogoliubov actions and crossed products.}
Let $G$ be a second countable locally compact group and let
$\pi:G\to\mathcal O(H_{\mathbb R})$ be a strongly continuous
orthogonal representation such that
\[
\pi_gU_t=U_t\pi_g
\qquad(g\in G,\ t\in\mathbb R).
\]
We use the same notation $\pi$ for its complexification.
Put $U_g^\pi:=\mathcal F_q(\pi_g)$ for $g\in G$.
These unitaries implement the $\varphi_q$-preserving action
\[
\sigma^\pi:G\curvearrowright M_q,
\qquad
\sigma^\pi_g:=\Ad(U_g^\pi),
\qquad
\sigma^\pi_g(W(j(\xi)))=W(j(\pi_g\xi))
\]
for $g\in G$ and $\xi\in H$. We call $\sigma^\pi$ the \emph{Bogoliubov action} associated with $\pi$.
Its restriction to $A_q$ is point-norm continuous.

Fix a left Haar measure on $G$.
Let $\delta_G$ be its modular function, and let $\lambda$ and $\rho$
be the left and right regular representations on $L^2(G)$, with
\[
(\lambda_gf)(s)=f(g^{-1}s),
\qquad
(\rho_gf)(s)=\delta_G(g)^{1/2}f(sg)
\quad(g\in G,\ f\in L^2(G)).
\]
For $x\in M_q$, define
$\pi_\sigma(x)\in\mathbb B(L^2(G,\mathcal F_q(H)))$ by
\[
(\pi_\sigma(x)\zeta)(s)
=
\sigma^\pi_{s^{-1}}(x)\zeta(s)
\qquad(\zeta\in L^2(G,\mathcal F_q(H)))
\]
for almost every $s\in G$. We realize the crossed product $M_q\rtimes_{\sigma^\pi} G$ on
$\mathcal F_q(H)\otimes L^2(G)
=
L^2(G,\mathcal F_q(H))$ as the von Neumann algebra generated by
$\pi_\sigma(M_q)$ and $1\otimes\lambda_g$ for $g\in G$. 
In this notation, we also have
\[
A_q\rtimes_{\sigma^\pi,\mathrm{red}}G
=
\overline{
	\left\{
	\int_G\pi_\sigma(a(g))(1\otimes\lambda_g)\,dg
	:
	a\in C_c(G,A_q)
	\right\}
}^{\|\cdot\|}.
\]
We often use the following notation:
\[
M_q(G):=M_q\rtimes_{\sigma^\pi}G,
\qquad
A_q(G):=A_q\rtimes_{\sigma^\pi,\mathrm{red}}G.
\]
Then, $A_q(G)$ is a $\sigma$-weakly dense $C^*$-subalgebra
of $M_q(G)$. It need not be unital.

Let $\widehat\varphi_q$ be the \emph{dual weight} on $M_q(G)$; see e.g.\ \cite[Chapter X]{Ta03}.
Under the above realization, we identify
$L^2(M_q(G),\widehat\varphi_q)$ with $\mathcal F_q(H)\otimes L^2(G)$.
Its modular conjugation $J_{G}$ is given by
\[
(J_{G}\zeta)(s)
=
\delta_G(s)^{-1/2}
U^\pi_{s^{-1}}J_q\zeta(s^{-1})
\qquad(\zeta\in L^2(G,\mathcal F_q(H)))
\]
for almost every $s\in G$.
Consequently, for $x\in M_q$ and $g\in G$,
\[
J_{G}\pi_\sigma(x)J_{G}
=
J_qxJ_q\otimes1,
\qquad
J_{G}(1\otimes\lambda_g)J_{G}
=
U^\pi_g\otimes\rho_g.
\]
In particular, $M_q(G)'=J_GM_q(G)J_G$ is generated by these elements.

Finally, if $G=\mathbb R$ and $\pi_s=U_{-s}$, then
$\sigma_s^\pi=\sigma_s^{\varphi_q}$, and hence
\[
M_q(\mathbb R)
=
M_q\rtimes_{\sigma^{\varphi_q}}\mathbb R
\]
is the continuous core of $M_q$.

\medskip
\noindent
\textbf{Second quantization.}
Let $(K_{\mathbb R},V_t)$ be another orthogonal representation of
$\mathbb R$, with complexification $K$ and embedding $j_K$ defined
as above.
If $T:H\to K$ is a contraction satisfying
\[
T(j(H_{\mathbb R}))\subset j_K(K_{\mathbb R}),
\]
then there is a unique normal unital completely positive map
\[
\Gamma_q(T):
\Gamma_q(H_{\mathbb R},U_t)''
\longrightarrow
\Gamma_q(K_{\mathbb R},V_t)''
\]
such that
\[
\Gamma_q(T)(W(\xi_1\otimes\cdots\otimes\xi_n))
=
W(T\xi_1\otimes\cdots\otimes T\xi_n)
\]
for $n\geq1$ and $\xi_1,\ldots,\xi_n\in jH$.
This map preserves the vacuum states, and its GNS implementation is
$\mathcal F_q(T)$.
See \cite[Corollary 3.16]{HR10} and \cite[Theorem~3.4]{Wa17}.

\subsection{Malleable deformations}
\label{Malleable deformations}

Let $M\subset \widetilde M$ be an inclusion of $\sigma$-finite von Neumann
algebras with expectation $E_M$, and let
$\alpha\colon \mathbb R\curvearrowright \widetilde M$ be a continuous action.
We call $(\alpha_t)_{t\in\mathbb R}$ a \emph{malleable deformation} of the
inclusion $M\subset\widetilde M$ \cite{Po03}.
If there exists $\beta\in\Aut(\widetilde M)$ such that
\[
\beta|_M=\id,\qquad
\beta^2=\id,\qquad
\beta\circ\alpha_t=\alpha_{-t}\circ\beta
\quad(t\in\mathbb R),
\]
we call $(\alpha,\beta)$ an \emph{s-malleable deformation}.

Following \cite[Section~4]{Ma16a}, we say that
\emph{$\alpha$ converges uniformly} on a bounded subset
$S\subset\widetilde M$ if, for every strong$^*$ neighborhood
$\mathcal V$ of $0$ in $\widetilde M$, there exists $t_0>0$ such that
\[
|t|<t_0
\quad\Longrightarrow\quad
\alpha_t(x)-x\in\mathcal V
\quad\text{for all }x\in S.
\]
Equivalently, for any fixed faithful normal state $\psi$ on $\widetilde M$,
\[
\sup_{x\in S}
\|\alpha_t(x)-x\|_\psi^\#
\longrightarrow0
\qquad(t\to0),
\]
where $\|x\|_\psi^\#
:=
\bigl(\psi(x^*x)+\psi(xx^*)\bigr)^{1/2}$.
In particular, this notion is independent of the choice of $\psi$.
If $S=S^*$, this is also equivalent to uniform convergence
with respect to $\|\cdot\|_\psi$.

\medskip
For $q$-Araki--Woods factors, we use the $q$-deformed version of the
s-malleable deformation introduced in \cite{HR10} and subsequently used
in the $q$-deformed setting; see, e.g., \cite{Av11,DP23,DT25}.
Let $G$ be a second countable locally compact group and let
$\pi\colon G\curvearrowright H_{\mathbb R}$ be a strongly continuous
orthogonal representation commuting with $(U_t)_{t\in\mathbb R}$.
Set $\widetilde{H}_\R:= H_\R\oplus H_\R$ and $\widetilde{H}=H\oplus H$, 
and extend $U_t$ and $\pi$ diagonally to these spaces.
We put
\[
M_q(G)
:=
\Gamma_q(H_{\mathbb R},U_t)''\rtimes_{\sigma^\pi}G
\subset 
\Gamma_q(\widetilde{H}_{\mathbb R},U_t)''
\rtimes_{\sigma^{\pi}}G
=:
\widetilde M_q(G),
\]
where we identify $M_q(G)$ with the subalgebra with expectation 
corresponding to the first copy of $H_{\mathbb R}$.
When $G$ is trivial, we simply write $M_q\subset\widetilde M_q$.
For $\theta\in\mathbb R$, let
\[
R_\theta
=
\begin{pmatrix}
	\cos\theta & -\sin\theta\\
	\sin\theta & \cos\theta
\end{pmatrix},
\qquad
\beta_0
=
\begin{pmatrix}
	1&0\\
	0&-1
\end{pmatrix}
\]
on $H_{\mathbb R}\oplus H_{\mathbb R}$.
Since $R_\theta$ and $\beta_0$ commute with $U_t\oplus U_t$ and
$\pi_g\oplus\pi_g$, their second quantizations induce automorphisms
\[
\alpha_\theta
:=
\Ad(\mathcal F_q(R_\theta)\otimes 1),
\qquad
\beta
:=
\Ad(\mathcal F_q(\beta_0)\otimes 1)
\]
of $\widetilde M_q(G)$.
Then $(\alpha,\beta)$ is an s-malleable deformation of the inclusion
$M_q(G)\subset\widetilde M_q(G)$.
The automorphisms $\alpha_\theta$ for $\theta\in \R$ and $\beta$ preserve the dual weight
associated with the vacuum state.

\subsection{Further preliminaries on von Neumann algebras}
\label{Further preliminaries on von Neumann algebras}

\medskip
\noindent
\textbf{Further facts from modular theory.}
Let $M$ be a von Neumann algebra and fix a faithful normal semifinite weight
$\varphi$ on $M$. We realize the \emph{standard form} of $M$ on the GNS Hilbert
space $L^2(M,\varphi)$, which we simply denote by $L^2(M)$, and denote its
natural positive cone by $L^2(M)^+$.
We denote the modular conjugation by $J$.
We write $C_\varphi(M):=M\rtimes_{\sigma^\varphi}\R$ for the continuous
core and denote by $L_\varphi\R$ the canonical copy of $L(\R)$
inside $C_\varphi(M)$.

For each positive functional $\psi\in M_\ast$, there exists a unique
$\xi_\psi\in L^2(M)_+$ such that
\[
\psi(x)=\langle\xi_\psi,x\xi_\psi\rangle
\qquad(x\in M).
\]
If $\psi$ is faithful, we denote its centralizer by
\[
M_\psi
:=
\{x\in M\mid \sigma_t^\psi(x)=x\text{ for all }t\in\mathbb R\}.
\]
By the KMS condition,
$\psi(xy)=\psi(yx)$ for $x\in M_\psi$ and $y\in M$,
and hence
\[
x\xi_\psi=\xi_\psi x \ (:=Jx^*J\xi_\psi)
\qquad(x\in M_\psi).
\]
In particular, $\psi|_{M_\psi}$ is a trace.

Let $p\in M$ be a projection and suppose that
$A\subset pMp$ is a finite von Neumann subalgebra with expectation $E$.
Then one can choose a faithful normal state $\psi$ on $M$ such that
$p\in M_\psi$ and $A\subset pM_\psi p$.
Indeed, let $\tau$ be a faithful normal tracial state on $A$.
If $p=1$, take $\psi=\tau\circ E$.
Otherwise, let $\omega$ be a faithful normal state on
$(1-p)M(1-p)$ and put
\[
\psi(x)
:=
\frac12\bigl(\tau(E(pxp))+\omega((1-p)x(1-p))\bigr),
\qquad x\in M.
\]
It is a faithful normal state with the desired properties.

Finally, every $\alpha\in\operatorname{Aut}(M)$ admits a canonical
unitary implementation $U_\alpha\in\mathcal U(L^2(M))$ satisfying
\[
U_\alpha xU_\alpha^*=\alpha(x)
\qquad(x\in M),
\]
\[
U_\alpha J=JU_\alpha,
\qquad
U_\alpha L^2(M)_+=L^2(M)_+.
\]
For these facts, we refer the reader to \cite{Ta03}.

\medskip
\noindent
\textbf{Popa's intertwining technique.}
We briefly recall \emph{Popa's intertwining technique} \cite{Po01,Po03}. 
Let $M$ be a $\sigma$-finite von Neumann algebra and let
$A\subset 1_AM1_A$ and $B\subset 1_BM1_B$ be von Neumann
subalgebras with expectation, where $1_A,1_B\in M$ are nonzero projections.
Fix a faithful normal conditional expectation
$E_B:1_BM1_B\to B$, and assume that $A$ is finite.

\begin{thm}
	The following conditions are equivalent.
	\begin{enumerate}
		\item there exist $n\geq1$, a nonzero projection $f\in B\otimes \M_n(\C)$, a normal unital $*$-homomorphism
		$\theta:A\to f(B\otimes \M_n(\C))f$, and a nonzero partial isometry $v\in \M_{1,n}(1_AM1_B)f$
		such that 
		\[
		av=v\theta(a)
		\qquad(a\in A),
		\]
		and the inclusion
		$\theta(A)\subset f\bigl(B \otimes\M_n(\C)\bigr)f$
		is with expectation.
		
		\item there exists no net $(u_i)_i\subset\mathcal U(A)$ such that
		\[
		E_B(y^*u_i x)\longrightarrow0
		\qquad(i\to\infty)
		\]
		in the $\sigma$-strong topology for all $x,y\in1_AM1_B$.
	\end{enumerate}
We write $A\preceq_M B$ if one of these conditions holds.
\end{thm}

Earlier versions of this result, including the tracial case, were established by Popa in
\cite{Po01,Po03}. For the general version stated above, see
\cite[Theorem~4.3]{HI15b}.
We note that if $B=\C1_M$, then $A\not\preceq_M \C$ holds if and only if $A$ is diffuse. Thus, $Ar \preceq_M \C$ holds for any nonzero projection $r\in A' \cap 1_AM1_A$ if and only if $A$ is completely atomic.

\medskip
\noindent
\textbf{Ultraproducts.}
Let $M$ be a $\sigma$-finite von Neumann algebra and let
$\omega\in\beta\mathbb N\setminus\mathbb N$ be a free ultrafilter.
Put
\[
\mathcal I_\omega(M)
:=
\left\{
(x_n)_n\in\ell^\infty(\mathbb N,M)
:
x_n\to0 \text{ strongly* as }n\to\omega
\right\};
\]
\[
\mathcal M^\omega(M)
:=
\left\{
(x_n)_n\in\ell^\infty(\mathbb N,M)
:
(x_n)_n\mathcal I_\omega(M)\subset\mathcal I_\omega(M),\
\mathcal I_\omega(M)(x_n)_n\subset\mathcal I_\omega(M)
\right\}.
\]
The \emph{Ocneanu ultraproduct}~\cite{Oc85} of $M$ is defined by
\[
M^\omega:=\mathcal M^\omega(M)/\mathcal I_\omega(M).
\]
It is a von Neumann algebra. We write $(x_n)_\omega$ for the class of
$(x_n)_n\in\mathcal M^\omega(M)$.
For any faithful normal state $\varphi$ on $M$, the formula
\[
\varphi^\omega((x_n)_\omega)
:=
\lim_{n\to\omega}\varphi(x_n)
\]
defines a faithful normal state on $M^\omega$.

If $M$ is finite, then every bounded sequence in $M$ defines an element of $M^\omega$.
More generally, if $N\subset M$ is a von Neumann subalgebra with
expectation, then the inclusion $N\subset M$ induces a canonical inclusion
$N^\omega\subset M^\omega$ with expectation.
Consequently, if $N\subset M$ is finite and with expectation, every bounded
sequence $(x_n)_n\in\ell^\infty(\mathbb N,N)$ defines an element
$(x_n)_\omega\in M^\omega$.

We identify $M$ with its diagonal copy in $M^\omega$. There is a canonical
faithful normal conditional expectation 
$E_\omega:M^\omega\to M$ given by 
$E_\omega((x_n)_\omega)
=
\lim_{n\to\omega}x_n$ in the $\sigma$-weak topology.

\section{An AO-type property for length projections}
\label{An AO-type property for length projections}

We first introduce some notation and terminology used in this section.
Let $K\subset H$ be a closed subspace of a complex Hilbert space and let $p_K\in\mathbb B(H)$ be the orthogonal projection onto $K$.
We identify $\mathcal F_q(K)$ with its natural subspace of
$\mathcal F_q(H)$, and denote by
$e_K=\mathcal F_q(p_K)$ the corresponding orthogonal projection.

For $r\in\mathbb N_0$, let $Q_r$ be the orthogonal
projection onto the closed span of elementary tensors whose legs belong
to $K$ or $K^\perp$, with exactly $r$ legs in $K^\perp$.
Thus
\[
Q_0=e_K,
\qquad
\sum_{r=0}^\infty Q_r=1,
\qquad
Q_rP_n=P_nQ_r,
\qquad
Q_rP_n=0\quad(r>n),
\]
where the sum converges strongly.
For $m\in\mathbb N_0$, we also put $Q_{\leq m}:=\sum_{s=0}^m Q_s  $ and $Q_{> m}:=\sum_{s>m} Q_s = 1-Q_{\leq m}$. We refer to $Q_r$, $Q_{\leq m}$, and $Q_{>m}$ as \emph{external projections}.
Each $Q_r$ commutes with $\ell_q(\xi)$ and $r_q(\xi)$ for every $\xi\in K$.

\subsection{Weak containment and length estimates}

In this subsection, we prove lemmas used in the proof of Theorem \ref{thmA}.
Let $H$ be a complex Hilbert space and identify $H$ with the first
summand of $\widetilde H:=H\oplus H$.
Let $e_H=\mathcal F_q(p_H)$ be the corresponding orthogonal projection, as above.
Let $R_\theta$ be a rotation matrix as in Subsection \ref{Malleable deformations}, and put $V_\theta:= \mathcal F_q(R_\theta)$ for $\theta\in \R$.

\begin{lem}\label{lem-ao-low-external-legs}
	Let $\pi\colon G\to \mathcal U(H)$ be a strongly continuous representation of a second countable locally compact group. Let $K\subset H$ be a finite-dimensional complex subspace, and put
	\[
	K(s)=K+\pi_{s^{-1}}K,
	\qquad s\in G.
	\]
	For $r\in \N_0$ and $s\in G$, let $Q_r^{K(s)} \in \B(\mathcal F_q(\widetilde{H}))$ be the external projection for the inclusion $K(s)\subset \widetilde{H}$.
	
	Then the following hold.
	\begin{enumerate}
		\item
		For every $r\in \N_0$, the formula
		\[
		(Q_r^G\xi)(s)=Q_r^{K(s)}\xi(s)
		\qquad\text{for almost every }s\in G
		\]
		for $\xi\in L^2(G,\mathcal F_q(\widetilde H))$,
		defines a projection $Q_r^G\in\mathbb B(\mathcal F_q(\widetilde H)\otimes L^2(G))$.
		
		\item 
		The projections $\{Q_r^G\}_{r\in \N_0}$ are mutually orthogonal, commute with the
		length projections, and their sum converges strongly to $1$.
		For each $m\in \N_0$, we set
		\[
		Q_{\leq m}^G=\sum_{r=0}^{m}Q_r^G,
		\qquad
		Q_{> m}^G=\sum_{r>m}Q_r^G=1-Q_{\leq m}^G.
		\]
		
		\item
		For every $m\in \N_0$ and $0<|\theta|<\pi/2$, we have
		\[
		\lim_{N\to\infty}
		\bigl\|Q_{\leq m}^GV_\theta e_H(1-P_{\le N})\bigr\|=0,
		\]
		where we suppress tensoring $V_\theta$, $e_H$ and $P_{\le N}$ with $1_{L^2(G)}$ for notational simplicity.
		
		Consequently, for every finite-band operator
		$T\in\mathbb B(\mathcal F_q(\widetilde H)\otimes L^2(G))$,
		\[
		e_HV_\theta^*TV_\theta e_H
		-
		e_HV_\theta^*Q_{> m}^GTQ_{> m}^GV_\theta e_H
		\in\mathcal K_\ell(G),
		\]
		where the compressed operators are regarded as operators on $\mathcal F_q(H)\otimes L^2(G)$.
		Both compressed operators belong to
		$\mathcal M(\mathcal K_\ell(G))$.
		
	\end{enumerate}
\end{lem}

\begin{proof}
	(1)
	Recall that an operator field $G\ni s\mapsto T(s)\in\mathbb B(\widetilde H)$
	is said to be strongly Borel measurable if
	$s\mapsto T(s)\xi$ is Borel measurable for every
	$\xi\in\widetilde H$.
	We use the same terminology for operator fields on the Fock space.
	
	Let $p_K$ be the projection onto $K$.
	The map
	\[
	s\longmapsto p_K+\pi_{s^{-1}}p_K\pi_s
	\]
	is norm continuous, since $p_K$ has finite rank.
	Its support projection is the projection $p(s)$ onto
	$K(s)=K+\pi_{s^{-1}}K$.
	Indeed, if we put
	$b(s):=p_K+\pi_{s^{-1}}p_K\pi_s$, then
	\[
	p(s)=\lim_{n\to\infty} b(s)^{1/n}\qquad \text{strongly}.
	\]
	This implies that $s\mapsto p(s)$ is strongly Borel measurable.
	
	We regard $p(s)$ as the projection on $\widetilde H$ with range
	$K(s)\subset H\subset\widetilde H$.
	For $z\in\mathbb T$, put
	\[
	D_z(s)=p(s)+z(1-p(s)).
	\]
	Then $s\mapsto\mathcal F_q(D_z(s))$ is strongly Borel measurable.
	On a simple tensor containing exactly $k$ legs in
	$\widetilde H\ominus K(s)$, the operator
	$\mathcal F_q(D_z(s))$ acts as multiplication by $z^k$.
	Hence
	\[
	\mathcal F_q(D_z(s))
	=
	\sum_{k\ge0}z^kQ_k^{K(s)}
	\]
	strongly, and therefore
	\[
	Q_r^{K(s)}
	=
	\int_{\mathbb T}
	z^{-r}\mathcal F_q(D_z(s))\,dz,
	\]
	where $dz$ denotes normalized Haar measure on $\mathbb T$ and the integral is taken in the strong operator topology.
	It follows that $s\mapsto Q_r^{K(s)}$ is strongly Borel measurable.
	Since $\|Q_r^{K(s)}\|\le1$, this defines the decomposable projection
	\[
	Q_r^G=\int_G^\oplus Q_r^{K(s)}\,ds
	\in\mathbb B(\mathcal F_q(\widetilde H)\otimes L^2(G)),
	\]
	which satisfies the desired condition. The assertions in (2) follow fibrewise.

	(3) Fix $m\in\N_0$ and $\theta\in\mathbb R$ with
	$0<|\theta|<\pi/2$.
	For each $r\in\N_0$, let $F_r\subset\mathcal F_q(\widetilde H)$
	be the range of the external projection associated with the inclusion
	$H\subset\widetilde H$, and denote by $P_{F_r}$ the corresponding
	orthogonal projection.
	Since
	\[
	Q_{\leq m}^{K(s)}
	\le
	\sum_{r=0}^{m}P_{F_r}
	\qquad (s\in G),
	\]
	it suffices to estimate the right-hand side.
	
	The computation at the beginning of the proof of
	\cite[Proposition~6.10]{Wi20} gives, for all
	$n,r\in\mathbb N_0$ with $r\le n$ and all
	$\xi,\eta\in H_q^{\otimes n}$,
	\[
	\bigl\langle
	P_{F_r}V_\theta e_H\xi,
	P_{F_r}V_\theta e_H\eta
	\bigr\rangle
	=
	\binom nr
	|\cos\theta|^{2(n-r)}
	|\sin\theta|^{2r}
	\langle\xi,\eta\rangle.
	\]
	Indeed, this is the same computation as in \cite[Proposition~6.10]{Wi20},
	with $e^{-t}$ and $(1-e^{-2t})^{1/2}$ replaced by
	$|\cos\theta|$ and $|\sin\theta|$, respectively.
	Consequently,
	\[
	\|P_{F_r}V_\theta e_HP_n\|^2
	=
	\binom nr
	|\cos\theta|^{2(n-r)}
	|\sin\theta|^{2r}.
	\]
	Hence for every integer $n\ge0$ and every $s\in G$,
	\[
	\|Q_{\leq m}^{K(s)}V_\theta e_HP_n\|^2
	\le
	\sum_{r=0}^{\min\{m,n\}}
	\binom nr
	|\cos\theta|^{2(n-r)}
	|\sin\theta|^{2r}.
	\]
	The right-hand side is independent of $s$ and converges to $0$ as $n\to\infty$. Applying this estimate fibrewise and using the fact that $Q_{\leq m}^G$, $V_\theta$, and $e_H$ preserve tensor length, we obtain
	\[
	\|Q_{\leq m}^GV_\theta e_H(1-P_{\le N})\|
	=
	\sup_{n>N}\|Q_{\leq m}^GV_\theta e_HP_n\|
	\longrightarrow0
	\qquad\text{as }N\to\infty.
	\]
	For the last assertion, fix a finite-band operator $T\in\mathbb B(\mathcal F_q(\widetilde H)\otimes L^2(G))$ and put
	\[
	S=V_\theta e_H,\qquad
	Q=Q_{>m}^G,\qquad
	L=(1-Q)S=Q_{\leq m}^GV_\theta e_H.
	\]
	By the preceding estimate, we have $LP_{\le N}\to L$ in norm as $N\to\infty$.
	Observe that
	\[
	T-QTQ=(1-Q)T+QT(1-Q),
	\]
	and therefore the difference in question is equal to
	\[
	\begin{aligned}
		S^*TS-S^*QTQS
		&=S^*(T-QTQ)S\\
		&=S^*(1-Q)TS+S^*QT(1-Q)S\\
		&=L^*TS+S^*QTL.
	\end{aligned}
	\]
	By the above convergence, we have
	\[
	P_{\le N}L^*TS\to L^*TS,
	\qquad
	S^*QTL P_{\le N}\to S^*QTL
	\]
	in norm as $N\to\infty$. Since $L^*TS$ and $S^*QTL$ are finite-band, 
	each of the operators on the left is supported
	between finite length projections, and hence belongs to
	$\mathcal K_\ell(G)$. We obtain 
	$L^*TS,\ S^*QTL\in\mathcal K_\ell(G)$.
	
	Finally, both compressed operators are finite-band and therefore
	belong to $\mathcal M(\mathcal K_\ell(G))$.
\end{proof}

The following is the form of weak containment that we will use.
It follows from the estimates in
\cite[Section 6]{Wi20} and \cite[Section 3]{DT25}.

\begin{lem}\label{lem-ao-weak-containment}
	Let $(H_\R,U_t)$ be a strongly continuous orthogonal representation such that $U$ is almost periodic.
	Let $K_{\mathbb R}\subset H_{\mathbb R}$ be a finite-dimensional
	$U$-invariant real subspace, and put
	\[
	N=\Gamma_q(K_{\mathbb R},U_t)''\subset \Gamma_q(H_{\mathbb R},U_t)''\subset \Gamma_q(\widetilde{H}_{\mathbb R},U_t)''.
	\]
	Here $\widetilde H_{\mathbb R}=H_{\mathbb R}\oplus H_{\mathbb R}$,
	with $H_{\mathbb R}$ identified with the first summand.
	If $m\in \N_0$ satisfies
	$(\dim K_{\mathbb R})|q|^{m+1}<1$, then
	\[
	Q_{> m} L^2(\widetilde M_q)
	\prec L^2(N)\otimes L^2(N)
	\]
	as $N$-$N$ bimodules. 
	Here $Q_{>m}$ is the external projection for the inclusion
	$K\subset\widetilde H$, where
	$K:=K_{\mathbb R}\otimes_{\mathbb R}\mathbb C$ is viewed as a subspace
	of the first summand of $\widetilde H$.
\end{lem}

\begin{proof}
	The case $K_{\mathbb R}=\{0\}$ is immediate.
	For $q\neq0$, this is exactly the weak containment estimate used in
	the last part of the proof of \cite[Lemma~3.9]{DT25}; see also
	\cite[Lemmas~3.4 and~3.5]{DT25}.
	The condition appearing there is
	$(\dim K_{\mathbb R})|q|^{m+1}<1$.
	The case $q=0$ follows from the free product decomposition.
\end{proof}

\subsection{Proof of Theorem~\ref{thmA}}

Let $G$ be a second countable locally compact group, and let
$U:\mathbb R\to\mathcal O(H_{\mathbb R})$ and
$\pi:G\to\mathcal O(H_{\mathbb R})$ be strongly continuous orthogonal
representations such that $\pi_gU_t=U_t\pi_g$ for all
$g\in G$ and $t\in\mathbb R$.
As in Subsection~\ref{$q$-Araki--Woods factors}, consider the inclusion
\[
A_q(G):=A_q\rtimes_{\sigma^\pi,\mathrm{red}}G
\subset
M_q\rtimes_{\sigma^\pi}G=:M_q(G).
\]
We also set
\[
A^{\mathrm u}_q(G):=A_q(G)+\mathbb C1\subset M_q(G).
\]
Thus $A^{\mathrm u}_q(G)=A_q(G)$ if $G$ is discrete, whereas
$A^{\mathrm u}_q(G)$ is the unitization of $A_q(G)$ if $G$ is
nondiscrete.

We first verify that the map appearing in Theorem~\ref{thmA} are well defined.
Let $J:=J_q$ and $J_G$ denote the modular conjugations associated with
$M_q$ and $M_q(G)$, respectively.

\begin{lem}\label{lem-ao-map}
	The algebras
	\[
	\pi_\sigma(A_q),\qquad
	JA_qJ\otimes1,\qquad
	A^{\mathrm u}_q(G),\qquad
	J_GA^{\mathrm u}_q(G)J_G
	\]
	and the unitaries
	$1\otimes\lambda_g$, $U_g^\pi\otimes\rho_g$ for $g\in G$
	are contained in $\mathcal M(\mathcal K_\ell(G))$.
	In particular, the formulas
	\[
	\begin{aligned}
		\nu_0:A_q\ota JA_qJ
		&\longrightarrow\mathcal M(\mathcal K_\ell(G)),\\
		a\otimes JbJ
		&\longmapsto\pi_\sigma(a)(JbJ\otimes1),
	\end{aligned}
	\]
	for $a,b\in A_q$, and
	\[
	\begin{aligned}
		\nu_0^G:
		A^{\mathrm u}_q(G)\ota J_GA^{\mathrm u}_q(G)J_G
		&\longrightarrow\mathcal M(\mathcal K_\ell(G)),\\
		x\otimes J_GyJ_G
		&\longmapsto xJ_GyJ_G,
	\end{aligned}
	\]
	for $x,y\in A^{\mathrm u}_q(G)$, define $*$-homomorphisms.
\end{lem}

\begin{proof}
	For every Wick polynomial $a\in A_q$, the operators
	$\pi_\sigma(a)$ and $JaJ\otimes1$ are finite-band.
	The unitaries $1\otimes\lambda_g$ and
	$U_g^\pi\otimes\rho_g$, for $g\in G$, commute with the length
	projections and hence are finite-band.
	Moreover, $A_q(G)$ is the norm-closed linear span of the operators
	\[
	\pi_\sigma(a)(1\otimes\lambda(f)),
	\qquad
	a\in A_q\text{ a Wick polynomial},\quad f\in C_c(G),
	\]
	where $\lambda(f):=\int_G f(g)\lambda_g\,dg$.
	These operators and their $J_G$-conjugates are finite-band.
	Hence
	$A_q(G)$ and $J_GA_q(G)J_G$ are contained in $\mathcal M(\mathcal K_\ell(G))$.
	Since $1\in\mathcal M(\mathcal K_\ell(G))$, the same holds for
	$A^{\mathrm u}_q(G)$ and $J_GA^{\mathrm u}_q(G)J_G$.
	
	Finally, since $M_q(G)$ commutes with $J_GM_q(G)J_G$, $\nu_0$ and $\nu_0^G $ define unital $\ast$-homomorphisms.
\end{proof}

Let $\nu_0$ and $\nu_0^G$ be the maps from
Lemma~\ref{lem-ao-map}, and let
\[
q_G:
\mathcal M(\mathcal K_\ell(G))
\longrightarrow
\mathcal M(\mathcal K_\ell(G))/\mathcal K_\ell(G)
\] be the quotient map. Put 
$\nu:=q_G\circ\nu_0$ and 
$\nu^G:=q_G\circ\nu_0^G$.

For later applications, we need not only the min-continuity of $\nu$
but also an approximate u.c.p.\ lifting formulation.
We therefore prove the following strengthening of
Theorem~\ref{thmA}.
The proof below does not rely on nuclearity of the $q$-creation algebra.

\begin{thm}\label{thm-length-AO}
	The following hold.
	\begin{enumerate}
		\item
		The map $\nu$ is continuous with respect to the minimal tensor norm, and hence it
		extends uniquely to a unital $*$-homomorphism
		\[
		\nu:D:=A_q\otimes_{\min}JA_qJ
		\longrightarrow
		\mathcal M(\mathcal K_\ell(G))/\mathcal K_\ell(G).
		\]
		
		\item
		There exists a net of u.c.p.\ maps
		$\Theta_i:D\longrightarrow\mathcal M(\mathcal K_\ell(G))$ such that, for every $z\in D$,
		\[
		\|q_G(\Theta_i(z))-\nu(z)\|\longrightarrow0
		\qquad (i\to\infty).
		\]
		
		\item
		If $G$ is amenable, then $\nu^G$ is continuous with respect to the minimal tensor norm, and
		extends uniquely to a unital $*$-homomorphism
		\[
		\nu^G:
		A^{\mathrm u}_q(G)\otimes_{\min}J_GA^{\mathrm u}_q(G)J_G
		\longrightarrow
		\mathcal M(\mathcal K_\ell(G))/\mathcal K_\ell(G).
		\]
	\end{enumerate}
\end{thm}

\begin{proof}[{\bf Proof of Theorem \ref{thm-length-AO}(1) and (2) in the almost periodic case}]
We first prove (1) and (2) under the assumption that $U$ is almost periodic.

Fix a finite-dimensional $U$-invariant real subspace
$K_{\mathbb R}\subset H_{\mathbb R}$, and put $K:=K_{\mathbb R}\otimes_{\mathbb R}\mathbb C$,
\[
A_{q,K}
:=
C^*\{1,W(j(\xi))\mid \xi\in K_{\mathbb R}\},
\qquad
D_K:=A_{q,K}\otimes_{\min}JA_{q,K}J.
\]
We prove the following claim.

\begin{claim}
There exists a net of u.c.p.\ maps $\Psi^K_i\colon D_K \to \mathcal M(\mathcal K_\ell(G))$ such that, for every $z\in A_{q,K}\ota JA_{q,K}J$,
\[
\|q_G(\Psi^K_i(z))-q_G(\nu_0(z))\|\to 0,
\qquad
\text{as }i\to \infty.
\]
\end{claim}
\begin{proof}[{\bf Proof of Claim}]	
	Let $\alpha$ be the malleable deformation for the inclusion
	\[
	M_q(G)
	=
	\Gamma_q(H_{\mathbb R},U_t)''\rtimes_{\sigma^\pi}G
	\subset 
	\Gamma_q(\widetilde{H}_{\mathbb R},U_t)''
	\rtimes_{\sigma^{\pi}}G
	=:
	\widetilde M_q(G),
	\]
	as in Subsection \ref{Malleable deformations}.
	Let $D_{K,\mathrm{alg}}\subset D_K$ denote the algebraic tensor product of the Wick polynomial $\ast$-subalgebra of $A_{q,K}$ and its conjugate by $J$.
For $s\in G$, set
\[
K_{\mathbb R}(s)
:=
K_{\mathbb R}+\pi_{s^{-1}}K_{\mathbb R},
\qquad
K(s):=K_{\mathbb R}(s)\otimes_{\mathbb R}\mathbb C.
\]
Since $\pi$ commutes with $U$, the space $K_{\mathbb R}(s)$ is
$U$-invariant, and
$\dim K_{\mathbb R}(s)\le 2\dim K_{\mathbb R}$.
Put $N_s:=\Gamma_q(K_{\mathbb R}(s),U_t)''$. Choose $m\in \N_0$ such that
\[
2(\dim K_{\mathbb R})|q|^{m+1}<1,
\]
and let $Q_{> m}^G$ be the projection as in
Lemma~\ref{lem-ao-low-external-legs}(2).

Since $\mathcal F_q(\widetilde H)\otimes L^2(G)$ is a $\widetilde{M}_q(G)$-bimodule, the restriction to $A_{q,K}$ gives the following $*$-homomorphism
\[
\widetilde\nu_{0,K}:
A_{q,K}\ota JA_{q,K}J
\longrightarrow
\mathbb B(\mathcal F_q(\widetilde H)\otimes L^2(G)).
\]
We show that $D_{K,\mathrm{alg}}\ni z \mapsto Q_{> m}^G\widetilde\nu_{0,K}(z)Q_{> m}^G$ is a $*$-homomorphism and is bounded with respect to the minimal tensor norm.

For each $z=\sum_{i=1}^k a_i\otimes Jb_iJ\in D_{K,\mathrm{alg}}$ with Wick polynomials $a_i,b_i\in A_{q,K}$, we have
\[
\widetilde\nu_{0,K}(z)
=
\sum_{i=1}^k
\pi_\sigma(a_i)(\widetilde Jb_i\widetilde J\otimes1)
\in\mathbb B(\mathcal F_q(\widetilde{H})\otimes L^2(G)).
\]
Here $\widetilde J$ is the modular conjugation of $\widetilde{M}_q(G)$.
Then $\widetilde z:=\widetilde\nu_{0,K}(z)$ is a decomposable operator with fibres
\[
\widetilde z(s)
=
\sum_{i=1}^k
\sigma^\pi_{s^{-1}}(a_i)\widetilde Jb_i\widetilde J\qquad (s\in G).
\]
Since $\sigma^\pi_{s^{-1}}(a_i),b_i\in N_s$, these fibres
preserve the number of legs outside $K(s)$.
Thus $\widetilde z$ commutes with every $Q_r^G$ and the above map is a $\ast$-homomorphism.
By Lemma~\ref{lem-ao-weak-containment}, for every $s\in G$,
\[
Q_{> m}^{K(s)}L^2(\widetilde M_q)
\prec L^2(N_s)\otimes L^2(N_s)
\]
as $N_s$-$N_s$ bimodules.
Consequently, we have
\[
\begin{aligned}
	\bigl\|
	Q_{> m}^{K(s)}\widetilde z(s)Q_{> m}^{K(s)}
	\bigr\|
	&\le
	\left\|
	\sum_{i=1}^k
	\sigma^\pi_{s^{-1}}(a_i)\otimes Jb_iJ
	\right\|_{\min}\\
	&=\|(\sigma^\pi_{s^{-1}}\otimes \id )(z)\|_{\min}=\|z\|_{\min}.
\end{aligned}
\]
Applying this fibrewise gives
$\|Q_{> m}^G\widetilde zQ_{> m}^G\|
\le\|z\|_{\min}$. We get the desired boundedness.

Let $e_H,V_\theta$ be as in the statement of 
Lemma~\ref{lem-ao-low-external-legs}(3), and we regard them as operators on $\mathcal F_q(\widetilde{H})\otimes L^2(G)$, by suppressing $\otimes 1_{L^2(G)}$. Then $e_H$ is the Jones projection for the inclusion $M_q(G)\subset \widetilde{M}_q(G)$.
For $0<\theta<\pi/2$, define
\[
\Psi_{K,\theta}(z)
:=
e_HV_\theta^*
Q_{> m}^G\widetilde zQ_{> m}^G
V_\theta e_H,
\qquad z\in D_{K,\mathrm{alg}}.
\]
The preceding observation shows that $\Psi_{K,\theta}$ extends to a
c.c.p.\ map
\[
\Psi_{K,\theta}:
D_K\longrightarrow\mathcal M(\mathcal K_\ell(G)).
\]
Indeed, for every $z\in D_{K,\mathrm{alg}}$, both $\widetilde z$ and $\Psi_{K,\theta}(z)$ are finite-band, and the range assertion follows by norm closure.

Lemma~\ref{lem-ao-low-external-legs}(3) gives
\[
q_G(\Psi_{K,\theta}(z))
=
q_G\bigl(e_HV_\theta^*\widetilde zV_\theta e_H\bigr)
\qquad
(z\in D_{K,\mathrm{alg}})
\]
and $q_G(\Psi_{K,\theta}(1))=1$. Choose a state $\omega_K$ on $D_K$ and put
\[
\widetilde\Psi_{K,\theta}(z)
:=
\Psi_{K,\theta}(z)
+
\omega_K(z)\bigl(1-\Psi_{K,\theta}(1)\bigr)\qquad (z\in D_K).
\]
Then $\widetilde\Psi_{K,\theta}$ is u.c.p., and
$q_G\circ\widetilde\Psi_{K,\theta}
=
q_G\circ\Psi_{K,\theta}$.

We next compare this map with $\nu_0$.
For each $z\in D_{K,\mathrm{alg}}$, 
the commutation of $V_\theta$ with $\widetilde J$ and with the diagonal action of $G$ gives
\[
e_HV_\theta^*\widetilde zV_\theta e_H
\longrightarrow
\nu_0(z)
\qquad \text{as }\theta\to0
\]
in operator norm.
Indeed, for an elementary tensor $z=a\otimes JbJ\in D_{K,\mathrm{alg}}$ with Wick polynomials $a,b\in A_{q,K}$, the corresponding
fibrewise difference is bounded uniformly in $s\in G$ by
\[
\|\alpha_{-\theta}(a)-a\|\,\|b\|
+
\|a\|\,\|\alpha_{-\theta}(b)-b\|,\qquad \text{where }\alpha_\theta:=\Ad(V_\theta),
\]
which tends to zero as $\theta\to0$ for Wick polynomials.
Consequently, for every $z\in D_{K,\mathrm{alg}}$,
\[
q_G(\widetilde\Psi_{K,\theta}(z))
\longrightarrow
q_G(\nu_0(z)),
\qquad
\text{as }\theta\to 0.
\]
By norm approximation, the same holds for every
$z\in A_{q,K}\ota JA_{q,K}J$.
Hence the family $(\widetilde\Psi_{K,\theta})_{0<\theta<\pi/2}$,
with $\theta$ directed towards $0$, is the required net.
\end{proof}

We now deduce assertion~(1).
For every finite-dimensional $U$-invariant real subspace $K_{\mathbb R}\subset H_{\mathbb R}$ and every $z\in A_{q,K}\ota JA_{q,K}J$, the claim and contractivity give
\[
\|\nu(z)\|
=
\lim_{i\to\infty}\|q_G(\Psi_i^K(z))\|
\le
\|z\|_{\min}.
\]
By almost periodicity, $\bigcup_K A_{q,K}$ is norm dense in $A_q$, where the union is taken over all finite-dimensional $U$-invariant real subspaces $K_{\mathbb R}\subset H_{\mathbb R}$.
Approximating the tensor factors in norm, we obtain the same inequality for every $z\in A_q\ota JA_qJ$.
Thus $\nu=q_G\circ\nu_0$ extends uniquely to a unital $*$-homomorphism on $D$.
This proves~(1).

It remains to obtain the maps in~(2).
For each finite-dimensional $U$-invariant real subspace $K_{\mathbb R}$, let $E_K:M_q\longrightarrow\Gamma_q(K_{\mathbb R},U_t)''$ be the vacuum-state-preserving conditional expectation.
Its restriction to $A_q$ takes values in $A_{q,K}$.
Define $E_K^J:JA_qJ\to JA_{q,K}J$ by $E_K^J(JbJ):=JE_K(b)J$ for $b\in A_q$, and put
\[
\mathcal E_K:=E_K\otimes E_K^J:D\longrightarrow D_K.
\]
These maps are u.c.p.
Direct the finite-dimensional $U$-invariant real subspaces $K_{\mathbb R}\subset H_{\mathbb R}$ by inclusion.
Then, for every $z\in D$, we have
\[
\|\mathcal E_K(z)-z\|\longrightarrow0
\qquad\text{as }K_{\mathbb R}\nearrow H_{\mathbb R},
\]
since $\bigcup_KD_K$ is norm dense in $D$ and $\mathcal E_K$ fixes $D_K$ pointwise.

For each such $K_\R$, take a net $(\Psi_i^K)_i$ as in the claim. Since $\nu$ is now contractive on $D$, contractivity and density imply that $q_G(\Psi_i^K(w))\to\nu(w)$ in norm as $i\to\infty$ for every $w\in D_K$. For each index $i$ of this net, put
\[
\Phi_{K,i}
:=
\Psi_i^K\circ\mathcal E_K:
D\longrightarrow\mathcal M(\mathcal K_\ell(G)).
\]
These maps are u.c.p., and, for every $z\in D$,
\[
\begin{aligned}
	\|q_G(\Phi_{K,i}(z))-\nu(z)\|
	&\le
	\bigl\|
	q_G\bigl(\Psi_i^K(\mathcal E_K(z))\bigr)
	-\nu(\mathcal E_K(z))
	\bigr\|\\
	&\quad+\|\mathcal E_K(z)-z\|.
\end{aligned}
\]
Given a finite subset $F\subset D$ and $\varepsilon>0$, first choose $K_{\mathbb R}$ so that the second term is smaller than $\varepsilon/2$ for every $z\in F$, and then choose $i$ so that the first term is smaller than $\varepsilon/2$ for every $z\in F$.
Indexing the resulting maps by finite subsets $F\subset D$ ordered by inclusion and positive tolerances $\varepsilon$ decreasing to $0$ yields the net in~(2).
This completes the almost periodic case.
\end{proof}

\begin{proof}[{\bf Proof of Theorem~\ref{thm-length-AO}(1) and (2) in the general case}]
	We use the almost periodic approximation established in
	Appendix~\ref{Almost periodic approximation}.
	Let $U^{(n)}$ be the almost periodic representations and
	$\theta_n:A_q\to A_q^{(n)}$ the u.c.p.\ maps given by
	Proposition~\ref{prop-almost-periodic-approximation-qaw}.
	Each $U^{(n)}$ commutes with $\pi$.
	
	All the corresponding $q$-Araki--Woods algebras are realized on the
	same Fock space $\mathcal F_q(H)$, with the same modular conjugation
	$J$ and the same implementing unitaries
	$U_g^\pi=\mathcal F_q(\pi_g)$.
	Hence we may use the common representation
	\[
	\pi_\sigma:
	\mathbb B(\mathcal F_q(H))
	\longrightarrow
	\mathbb B(\mathcal F_q(H)\otimes L^2(G)),
	\]
	\[
	(\pi_\sigma(T)\xi)(s)
	=
	\Ad(U_{s^{-1}}^\pi)(T)\xi(s)\quad (s\in G),
	\qquad
	T\in\mathbb B(\mathcal F_q(H)).
	\]
	Its restrictions give the left representations of all the crossed
	products under consideration.
	The length projections, and hence $\mathcal K_\ell(G)$ and $q_G$,
	are also independent of $n$.
	
	Define the u.c.p.\ map
	\[
	\psi_n:=\theta_n\otimes\theta_n^J:
	D\longrightarrow
	D^{(n)}:=A_q^{(n)}\otimes_{\min}JA_q^{(n)}J,
	\]
	where
	\[
	\theta_n^J:JA_qJ\longrightarrow JA_q^{(n)}J,
	\qquad
	\theta_n^J(JbJ):=J\theta_n(b)J
	\quad (b\in A_q).
	\]
	By the almost periodic case proved above, the corresponding AO maps
	extend to unital $*$-homomorphisms
	\[
	\nu^{(n)}:D^{(n)}
	\longrightarrow
	\mathcal M(\mathcal K_\ell(G))/\mathcal K_\ell(G).
	\]
	For each
	$z=\sum_{k=1}^r a_k\otimes Jb_kJ
	\in A_q\ota JA_qJ$,
	the norm convergence $\theta_n(a)\to a$ for $a\in A_q$ in the common Fock
	representation gives
	\[
	\begin{aligned}
		\nu^{(n)}(\psi_n(z))
		&=
		q_G\left(
		\sum_{k=1}^r
		\pi_\sigma(\theta_n(a_k))
		(J\theta_n(b_k)J\otimes1)
		\right)\\
		&\longrightarrow
		q_G\left(
		\sum_{k=1}^r
		\pi_\sigma(a_k)(Jb_kJ\otimes1)
		\right)
		=
		\nu(z),\qquad \text{as }n\to \infty
	\end{aligned}
	\]
	in norm.
	Since $\nu^{(n)}\circ\psi_n$ is contractive, we obtain
	$\|\nu(z)\|\le\|z\|_{\min}$.
	Thus $\nu$ extends uniquely to a unital $*$-homomorphism on $D$,
	which proves~(1).
	By contractivity and density, the same convergence holds for all $z\in D$.
	
	For each $n$, the almost periodic case provides a net of
	u.c.p.\ maps
	$\Theta_i^{(n)}:
	D^{(n)}\to\mathcal M(\mathcal K_\ell(G))$
	such that
	$q_G\circ\Theta_i^{(n)}
	\to
	\nu^{(n)}$
	in the point-norm topology.
	The compositions
	\[
	\Theta_i^{(n)}\circ\psi_n:
	D\longrightarrow\mathcal M(\mathcal K_\ell(G))
	\]
	are again u.c.p., and for each $z\in D$,
	\[
	\begin{aligned}
		\bigl\|
		q_G(\Theta_i^{(n)}(\psi_n(z)))-\nu(z)
		\bigr\|
		&\le
		\bigl\|
		q_G(\Theta_i^{(n)}(\psi_n(z)))
		-\nu^{(n)}(\psi_n(z))
		\bigr\|\\
		&\quad+
		\|\nu^{(n)}(\psi_n(z))-\nu(z)\|.
	\end{aligned}
	\]
	Given a finite subset $F\subset D$ and $\varepsilon>0$, first choose
	$n$ so that the second term is smaller than $\varepsilon/2$ for all
	$z\in F$, and then choose $i$ so that the first term is smaller than
	$\varepsilon/2$ for all $z\in F$.
	Directing the resulting maps by pairs $(F,\varepsilon)$ gives the
	required net in~(2).
\end{proof}

\begin{proof}[{\bf Proof of Theorem~\ref{thm-length-AO}(3)}]
	Assume that $G$ is amenable, and put $A_q(G)=A_q\rtm G$.
	Consider the action
	\[
	\beta:G\times G\curvearrowright D=A_q\otimes_{\min}JA_qJ,
	\qquad
	\beta_{(g,h)}=\sigma^\pi_g\otimes\sigma^{\pi,\mathrm r}_h\quad \text{for }g,h\in G,
	\]
	where $\sigma^{\pi,\mathrm r}_h(JaJ)=J\sigma^\pi_h(a)J$ for $a\in A_q$.
	Set
	\[
	W_{g,h}:=(1\otimes\lambda_g)(U_h^\pi\otimes\rho_h)\qquad (g,h\in G).
	\]
	These unitaries implement $\beta$ on the image of $\nu_0$.
	
	Consider the canonical inclusion
	\[
	\mathscr C
	:=
	C_c(G,A_q)\ota C_c(G,JA_qJ)
	\subset C_c(G\times G,D),
	\]
	and note that 
	$\mathscr C$ is an $L^1$-dense $*$-subalgebra of
	$L^1(G\times G,D,\beta)$.
	The regular representations give the canonical identification
	\[
	D\rtimes_{\beta,\mathrm{red}}(G\times G)
	\simeq A_q(G)\otimes_{\min}J_GA_q(G)J_G.
	\]
	Under this identification, $\mathscr C$ is also a dense
	$*$-subalgebra of the right-hand side, and the restriction of
	$\nu^G$ to $\mathscr C$ is given by
	\[
	\nu^G(F)=q_G(X_F),
	\qquad
	X_F:=
	\int_{G\times G}
	\nu_0(F(g,h))W_{g,h}\,dg\,dh
	\in\mathcal M(\mathcal K_\ell(G)),
	\]
	for $F\in \mathscr C$. 
	Here the integral is taken in the strong operator topology.
	For an elementary tensor $F=x\otimes y\in \mathscr C$, it is the product of
	the integrated left representation of $x$ and the integrated
	right representation of $y$.

	We show that
	\[
	\|\nu^G(F)\|\le \|F\|_1
	\qquad(F\in\mathscr C).
	\]
	Put $e_N:=P_{\le N}\otimes1$.
	Since $(e_N)_N$ is an approximate identity for
	$\mathcal K_\ell(G)$, for every
	$F\in\mathscr C$ and $(g,h)\in G\times G$, we have
	\[
	\|\nu^G(F)\|
	=
	\lim_{N\to\infty}\|(1-e_N)X_F\|,
	\qquad 
	\|\nu(F(g,h))\|
	=
	\lim_{N\to\infty}
	\|(1-e_N)\nu_0(F(g,h))\|.
	\]
	Writing 
	$F=\sum_{j=1}^k x_j\otimes y_j$ for some 
	$x_j\in C_c(G,A_q)$ and $y_j\in C_c(G,JA_qJ)$, 
	for every $N$, we have
	\[
	\|(1-e_N)\nu_0(F(g,h))\|
	\le
	\sum_{j=1}^k\|x_j(g)\|\,\|y_j(h)\|\qquad (g,h\in G).
	\]
	The right-hand side is integrable and independent of $N$.
	Combining the preceding pointwise limit with dominated convergence,
	we obtain
	\[
	\begin{aligned}
		\|\nu^G(F)\|
		&=
		\lim_{N\to\infty}\|(1-e_N)X_F\|\\
		&\le
		\lim_{N\to\infty}
		\int_{G\times G}
		\|(1-e_N)\nu_0(F(g,h))\|\,dg\,dh\\
		&=
		\int_{G\times G}
		\|\nu(F(g,h))\|\,dg\,dh\\
		&\le
		\int_{G\times G}
		\|F(g,h)\|_{\min}\,dg\,dh
		=
		\|F\|_1,
	\end{aligned}
	\]
	where the last inequality follows from~(1).

	Thus $\nu^G|_{\mathscr C}$ extends to a contractive
	$*$-homomorphism on $L^1(G\times G,D,\beta)$, and hence to a
	$*$-homomorphism on the full crossed product
	$D\rtimes_{\beta,\mathrm{full}}(G\times G)$.
	Since $G\times G$ is amenable, the full and reduced crossed
	products coincide.
	Using the canonical identification above, we obtain a
	$*$-homomorphism on $A_q(G)\otimes_{\min}J_GA_q(G)J_G$.
	Approximating the two tensor factors separately by integrated
	elements of $C_c(G,A_q)$ and $C_c(G,JA_qJ)$ shows that this
	extension agrees with $\nu^G$ on $A_q(G)\ota J_GA_q(G)J_G$.
	This proves its min-continuity.

	It remains to consider the unitizations.
	There is nothing to prove if $A_q(G)$ is unital.
	Otherwise, let $A^{\mathrm u}_q(G)=A_q(G)+\mathbb C1$ and let
	$\varepsilon:A^{\mathrm u}_q(G)\to\mathbb C$ be the quotient character.
	For $z\in A^{\mathrm u}_q(G)\ota J_GA_q(G)J_G$, put
	\[
	w=(\varepsilon\otimes\id)(z)\in J_GA_q(G)J_G.
	\]
	Since
	\[
	z-1\otimes w\in A_q(G)\ota J_GA_q(G)J_G
	\quad\text{and}\quad
	\|w\|\le\|z\|_{\min},
	\]
	the preceding result gives
	\[
	\|\nu^G(z)\|
	\le
	\|z-1\otimes w\|_{\min}+\|w\|
	\le
	3\|z\|_{\min}.
	\]
	Hence the multiplication map remains min-continuous after
	adjoining a unit to the first factor.
	Applying the same argument to the second factor yields the
	required unital $*$-homomorphism on
	$A^{\mathrm u}_q(G)\otimes_{\min}J_GA^{\mathrm u}_q(G)J_G$.
\end{proof}

We will later need the following stabilized form in the proof of
strong solidity. This is also the reason for formulating
Theorem~\ref{thm-length-AO} in terms of approximate u.c.p.\ lifts,
rather than merely proving the min-continuity of the multiplication
map. Indeed, after tensoring with an arbitrary unital $C^*$-algebra
$B$, the min-continuity of $\nu$ alone does not directly give the
required multiplication map into the corresponding quotient
after amplification. On the other hand, the approximate u.c.p.\ lifts
in Theorem~\ref{thm-length-AO}(2) can be tensored with $\id_B$,
yielding both min-continuity and approximate u.c.p.\ lifts
for the stabilized multiplication map. 

\begin{cor}\label{cor-stable-AO}
	We keep the setting of Theorem~\ref{thm-length-AO}, and let $B\subset\mathbb B(\mathcal H_B)$ be a unital $C^*$-algebra.
	Let $\mathcal K_{\ell,B}(G)$ be the hereditary $C^*$-subalgebra of
	\[
	\mathbb B\bigl(\mathcal F_q(H)\otimes L^2(G)\otimes\mathcal H_B\bigr)
	\]
	generated by the projections $\{P_{\le N}\otimes1\otimes1\}_N$, and denote
	the corresponding quotient map on $\mathcal M(\mathcal K_{\ell,B}(G))$ by $q_{G,B}$.
	Define the algebraic maps
	\[
	\begin{aligned}
		\nu_B(z\otimes b)
		&:=q_{G,B}\bigl(\nu_0(z)\otimes b\bigr),
		&&z\in A_q\ota JA_qJ,\quad b\in B,\\
		\nu_B^G(z\otimes b)
		&:=q_{G,B}\bigl(\nu_0^G(z)\otimes b\bigr),
		&&z\in A^{\mathrm u}_q(G)\ota J_GA^{\mathrm u}_q(G)J_G,\quad b\in B,
	\end{aligned}
	\]
	Then assertions~(1) and~(2) of Theorem~\ref{thm-length-AO}
	hold for $\nu_B$ on $D\otimes_{\min}B$, with
	$\mathcal K_\ell(G)$ and $q_G$ replaced by
	$\mathcal K_{\ell,B}(G)$ and $q_{G,B}$, respectively.
	Assertion~(3) also holds for $\nu_B^G$ on $A^{\mathrm u}_q(G)\otm J_GA^{\mathrm u}_q(G)J_G\otm B$.
\end{cor}
\begin{proof}
	Since $\{P_{\le N}\otimes1\}_{N\geq 0}$ and $\{P_{\le N}\otimes1\otimes 1_B\}_{N\geq 0}$ are approximate identities
	for $\mathcal K_\ell(G)$ and $\mathcal K_{\ell,B}(G)$, respectively, we have the natural inclusions
	\[
	\mathcal K_{\ell}(G)\otimes_{\min}B\subset \mathcal K_{\ell,B}(G),
	\qquad
	\mathcal M(\mathcal K_{\ell}(G))\otimes_{\min}B\subset\mathcal M(\mathcal K_{\ell,B}(G)).
	\]
	Consequently, for $T\in\mathcal M(\mathcal K_\ell(G))$ and $b\in B$,
	\[
	\|q_{G,B}(T\otimes b)\|
	\le
	\inf_{k\in\mathcal K_\ell(G)}\|(T-k)\otimes b\|
	=
	\|q_G(T)\|\,\|b\|.
	\]
	
	Let $(\Theta_i)_i$ be the net from
	Theorem~\ref{thm-length-AO}(2).
	By the preceding inclusion, the maps
	\[
	\Theta_{i,B}:=\Theta_i\otimes\id_B:
	D\otimes_{\min}B\longrightarrow\mathcal M(\mathcal K_{\ell}(G))\otimes_{\min}B\subset\mathcal M(\mathcal K_{\ell,B}(G))
	\]
	are u.c.p.
	For each
	\[
	x=\sum_{j=1}^r z_j\otimes b_j
	\in (A_q\ota JA_qJ)\ota B,
	\]
	the preceding quotient norm estimate gives
	\[
	\begin{aligned}
		\|q_{G,B}(\Theta_{i,B}(x))-\nu_B(x)\|
		&\le
		\sum_{j=1}^r
		\|q_G(\Theta_i(z_j))-\nu(z_j)\|\,\|b_j\|\\
		&\longrightarrow0,\qquad \text{as }i\to \infty.
	\end{aligned}
	\]
	It follows that
	$\|\nu_B(x)\|\le\|x\|_{\min}$.
	Thus $\nu_B$ extends uniquely to a unital $*$-homomorphism on
	$D\otimes_{\min}B$.
	By contractivity and density, the same convergence holds for every
	$x\in D\otimes_{\min}B$.
	This proves assertions~(1) and~(2) for $\nu_B$.
	
	Assume that $G$ is amenable.
	We apply the proof of Theorem~\ref{thm-length-AO}(3) with
	\[
	D,\ \beta,\ \mathcal K_\ell(G),\ q_G,\ \nu
	\]
	replaced, respectively, by
	\[
	D_B:=D\otimes_{\min}B,\qquad
	\beta_B:=\beta\otimes\id_B,\qquad
	\mathcal K_{\ell,B}(G),\qquad
	q_{G,B},\qquad
	\nu_B.
	\]
	The implementing unitaries are $W_{g,h}\otimes1$ for
	$(g,h)\in G\times G$, and the corresponding dense $*$-subalgebra is
	\[
	\mathscr C_B
	:=
	C_c(G,A_q)\ota C_c(G,JA_qJ)\ota B
	\subset
	L^1(G\times G,D_B,\beta_B).
	\]
	By the same argument as in the proof of
	Theorem~\ref{thm-length-AO}(3), we obtain
	\[
	\|\nu_B^G(F)\|
	\le
	\|F\|_{L^1(G\times G,D_B)}
	\qquad
	(F\in\mathscr C_B).
	\]
	Since $G\times G$ is amenable and
	\[
	D_B\rtimes_{\beta_B,\mathrm{red}}(G\times G)
	\cong
	A_q(G)\otimes_{\min}J_GA_q(G)J_G\otimes_{\min}B,
	\]
	we obtain the required min-continuous extension of $\nu_B^G$.
	
	If, in addition, $G$ is nondiscrete, the same unitization argument as
	in the proof of Theorem~\ref{thm-length-AO}(3) gives the corresponding
	extension after adjoining units to the two crossed-product factors.
\end{proof}

\section{A characterization of uniform convergence of malleable deformations}
\label{A characterization of uniform convergence of malleable deformations}

In this section, we study the $s$-malleable deformation introduced in
Subsection~\ref{Malleable deformations}.
Throughout this section, fix a strongly continuous orthogonal representation
$U:\mathbb R\to\mathcal O(H_{\mathbb R})$ and a strongly continuous orthogonal representation
$\pi:G\to\mathcal O(H_{\mathbb R})$ of a second countable locally compact group $G$, commuting with $U$.
We consider the associated inclusion
\[
M_q(G)
=
\Gamma_q(H_{\mathbb R},U_t)''\rtimes_{\sigma^\pi}G
\subset
\Gamma_q(\widetilde H_{\mathbb R},U_t)''\rtimes_{\sigma^\pi}G
=
\widetilde M_q(G),
\]
and write $(\alpha,\beta)$ for its $s$-malleable deformation.

For the rotation $R_\theta$ appearing in the definition of $\alpha$, set
\[
V_\theta:=\mathcal F_q(R_\theta)\otimes 1
\in \mathcal U(\mathcal F_q(\widetilde H)\otimes L^2(G))
\qquad (\theta\in\mathbb R),
\]
and for the orthogonal projection $p_H:H\oplus H\to H\oplus 0$, set
\[
e_H:=\mathcal F_q(p_H)\otimes 1
\in \mathbb B(\mathcal F_q(\widetilde H)\otimes L^2(G)).
\]
Then $V_\theta$ is the standard implementation of $\alpha_\theta$, and $e_H$ is the Jones projection for the inclusion
$M_q(G)\subset\widetilde M_q(G)$.
This notation agrees with that used in Lemma~\ref{lem-ao-low-external-legs}(3).

We denote by $P_n$, $P_{\leq N}$, and $P_{>N}$ the length projections on
$\mathcal F_q(H)$. Unlike the notation above, we do not suppress the tensor factor
$1_{L^2(G)}$ when these projections act on
$\mathcal F_q(H)\otimes L^2(G)$.

\subsection{Uniform convergence and length projections}

We first relate uniform convergence of the deformation to length projections.
We use the natural inclusion
\[
L^2(M_q(G))
=\mathcal F_q(H)\otimes L^2(G)
\subset
\mathcal F_q(H\oplus H)\otimes L^2(G)
=L^2(\widetilde M_q(G)),
\]
which intertwines the left actions of $M_q(G)$.

\begin{lem}\label{lem-deformation-uniform-length0}
	For every $\eta\in \mathcal F_q(H)\otimes L^2(G)$ and $\theta\in\mathbb R$, we have
	\begin{enumerate}
		\item $\displaystyle
		e_HV_\theta\eta
		=\sum_{n\ge0}(\cos\theta)^n(P_n\otimes1)\eta$;
		
		\item $\displaystyle
		\|V_\theta\eta-\eta\|^2
		=
		2\sum_{n\ge0}\bigl(1-(\cos\theta)^n\bigr)
		\|(P_n\otimes1)\eta\|^2$.
	\end{enumerate}
\end{lem}

\begin{proof}
	Since $p_HR_\theta$ restricts to multiplication by $\cos\theta$
	on $H\oplus0$, functoriality of second quantization gives (1).
	Taking the inner product with $\eta$ yields (2), since
	$e_H\eta=\eta$ and $V_\theta$ is unitary.
\end{proof}

\begin{lem}\label{lem-deformation-uniform-length1}
	Let $S\subset(M_q(G))_1$ and $\xi\in L^2(M_q(G))$.
	Define $\psi(x):=\langle\xi,x\xi\rangle$ for
	$x\in\widetilde M_q(G)$, and write
	$\|x\|_\psi:=\psi(x^*x)^{1/2}$.
	The following conditions are equivalent:
	\begin{enumerate}
		\item $\displaystyle
		\sup_{a\in S}\|\alpha_\theta(a)-a\|_\psi\longrightarrow0
		\quad\text{as }\theta\to0$;
		
		\item $\displaystyle
		\sup_{a\in S}\|(P_{> N}\otimes1)a\xi\|\longrightarrow0
		\quad\text{as }N\to\infty$.
	\end{enumerate}
\end{lem}

\begin{proof}
	Since $V_\theta a\xi
	=\alpha_\theta(a)V_\theta\xi$, we have
	\[
	\left|
	\|V_\theta a\xi-a\xi\|
	-\|\alpha_\theta(a)-a\|_\psi
	\right|
	\le \|V_\theta\xi-\xi\|.
	\]
	The right-hand side tends to zero independently of $a\in S$.
	Thus (1) is equivalent to
	\[
	\sup_{a\in S}\|V_\theta a\xi-a\xi\|
	\longrightarrow0
	\quad\text{as }\theta\to0.
	\]
	For each fixed $0<|\theta|<\pi/2$, we have
	$1-(\cos\theta)^n\ge1/2$ for all sufficiently large $n$.
	Lemma~\ref{lem-deformation-uniform-length0}(2) therefore shows
	that (1) implies (2).
	Conversely, the same formula gives, for $|\theta|<\pi/2$,
	\[
	\|V_\theta a\xi-a\xi\|^2
	\le
	2\bigl(1-(\cos\theta)^{N+1}\bigr)\|\xi\|^2
	+2\|(P_{> N}\otimes1_{L^2(G)})a\xi\|^2.
	\]
	Taking $N$ large and then $\theta$ small proves that (2) implies (1).
\end{proof}

We record the following standard convexity lemma; see, for instance,
\cite[proof of Theorem~4.3, Step~(1)]{HR10} and
\cite[proof of Theorem~4.2, Step~2]{Ma16a}.

\begin{lem}\label{lem-uniform-convergence-intertwine}
	Let $(M,\psi)$ be a von Neumann algebra with a faithful normal state,
	let $\alpha\in\Aut(M)$, and let $p\in M_\psi$ be a projection.
	Let $A\subset pM_\psi p$ be a von Neumann subalgebra. If
	\[
	\sup_{u\in\mathcal U(A)}
	\|\alpha(u)-u\|_\psi
	<
	\|p\|_\psi,
	\]
	then there exists a nonzero element
	$z\in\alpha(p)Mp$ such that
	\[
	\alpha(a)z=za
	\qquad(a\in A).
	\]
\end{lem}

\subsection{Proof of Theorem~\ref{thmB}: case (a)}

Throughout this subsection, assume that $G=\{e\}$ and write
$M_q=M_q(\{e\})$ and $\widetilde M_q=\widetilde M_q(\{e\})$.

The following is the key lemma of this section. The proof is inspired by
that of \cite[Lemma~8.5]{DP23}. Instead of the conditional expectations used there, we use
the finite-rank contractions appearing in the proof of the Haagerup property
for $q$-Araki--Woods algebras.

\begin{lem}\label{lem-intertwiner-inequality}
	Let $p\in M_q$ be a projection and $A\subset pM_qp$ a finite
	von Neumann subalgebra with expectation $E_A$.
	Let $\psi$ be a faithful normal state on $M_q$ such that $p\in (M_q)_\psi$ and $A\subset p(M_q)_\psi p$. Let
	$\xi_\psi\in L^2(M_q)\subset L^2(\widetilde M_q)$ be its implementing
	vector in the natural positive cone.
	
	Suppose that, for some $\theta\in(0,\pi/2)$, there exists
	$0\ne z\in\alpha_\theta(p)\widetilde M_qp$ such that
	\[
	\alpha_\theta(a)z=za
	\qquad(a\in A).
	\]
	Then there exist a finite-dimensional subspace $K\subset H$ and
	a constant $c>0$ such that
	\[
	\inf_{u\in\mathcal U(A)}
	\|P_{\mathcal F_q(K)}u\xi_\psi\|\ge c.
	\]
	Here $P_{\mathcal F_q(K)}$ denotes the orthogonal projection onto
	$\mathcal F_q(K)$.
\end{lem}

\begin{proof}
	Normalize $z$ so that $\|z\|=1$, and put $m:=\|z\xi_\psi\|>0$.
	Fix $0<\varepsilon<m/2$, to be chosen sufficiently small below.

	By \cite[Proposition~3.17]{HR10}, there exists a net of
	finite-rank contractions $T_i$ on $H$ such that $T_i\to 1$
	strongly and 
	$T_i(j(H_{\mathbb R}))\subset j(H_{\mathbb R})$.
	Although \cite{HR10} assumes
	separability, the proof of this proposition works verbatim
	without this assumption. Hence, by \cite[Theorem~3.4]{Wa17},
	$\Gamma_q(1\oplus T_i)$ defines a vacuum-state-preserving
	normal u.c.p.\ map on $\widetilde M_q$. Set
	\[
	\Psi_i
	:=
	\alpha_\theta\circ\Gamma_q(1\oplus T_i)\circ\alpha_{-\theta},
	\qquad
	F_i
	:=
	\mathcal F_q\bigl(R_\theta(1\oplus T_i)R_{-\theta}\bigr).
	\]
	Then $F_i x\Omega=\Psi_i(x)\Omega$ for $x\in\widetilde M_q$,
	and $F_i\to1$ strongly.
	
	Since $\Gamma_q(1\oplus T_i)$ fixes $M_q$ pointwise, $\Psi_i$ fixes
	$\alpha_\theta(M_q)$ pointwise. Thus $\alpha_\theta(M_q)$ is contained
	in the multiplicative domain of $\Psi_i$, and for $a\in M_q$ and
	$x\in\widetilde M_q$ we have
	\[
	\begin{aligned}
		F_i\alpha_\theta(a)x\Omega
		&=\Psi_i(\alpha_\theta(a)x)\Omega\\
		&=\alpha_\theta(a)\Psi_i(x)\Omega
		=\alpha_\theta(a)F_i x\Omega.
	\end{aligned}
	\]
	By density of $\widetilde M_q\Omega$, it follows that
	\[
	F_i\alpha_\theta(a)=\alpha_\theta(a)F_i
	\qquad(a\in M_q).
	\]
	The multiplicative-domain argument above is precisely what allows these
	u.c.p.\ maps to play the role of the conditional expectations in
	\cite[Lemma~8.5]{DP23}.
	
	Choose $i$ and a finite Wick polynomial
	$y\in\widetilde M_q$ with $\|y\|\le1$ such that
	\[
	\|(F_i-1)z\xi_\psi\|<\varepsilon,
	\qquad
	\|(z-y)\xi_\psi\|<\varepsilon.
	\]
	Since $a\xi_\psi=\xi_\psi a$ for $a\in A$ and $z=zp$,
	right multiplication by $u\in\mathcal U(A)$ gives
	\[
	\|zu\xi_\psi\|=m,
	\qquad
	\|(z-y)u\xi_\psi\|
	\le\|(z-y)\xi_\psi\|<\varepsilon.
	\]
	Moreover, the intertwining relation and the commutation relation above imply
	\[
	\|(F_i-1)zu\xi_\psi\|
	=
	\|\alpha_\theta(u)(F_i-1)z\xi_\psi\|
	<\varepsilon.
	\]
	Consequently,
	\[
	\|F_i yu\xi_\psi\|\ge m-2\varepsilon,
	\qquad
	\|yu\xi_\psi\|\le m+\varepsilon
	\qquad(u\in\mathcal U(A)).
	\]
	Since $y$ is a finite Wick polynomial, we may choose a finite-dimensional
	subspace $K\subset H$ containing
	$\operatorname{Ran}T_i+\operatorname{Ran}T_i^*$ 	
	such that all vectors appearing in the Wick words of $y$ belong to
	$K\oplus K$.
	Then $T_iK\subset K$ and $T_i|_{K^\perp}=0$.
	The generators occurring in $y$ involve only creation and annihilation
	operators with vectors in $K\oplus K$, so
	\[
	yu\xi_\psi\in\mathcal F_q(H\oplus K)
	\qquad(u\in\mathcal U(A)).
	\]
	Put $P:=P_{\mathcal F_q(K\oplus K)}$.
	This projection commutes with those creation and annihilation operators,
	and hence with $y$.
	Since $u\xi_\psi\in\mathcal F_q(H\oplus0)$ and $\|y\|\le1$, we obtain
	\[
	\|Pyu\xi_\psi\|
	=
	\|yP_{\mathcal F_q(K\oplus0)}u\xi_\psi\|
	\le
	\|P_{\mathcal F_q(K)}u\xi_\psi\|.
	\]
	Combining the preceding estimates with
	Lemma~\ref{lem-deformation-inequality} below, we obtain
	\[
	\begin{aligned}
		(m-2\varepsilon)^2
		&\le \|F_i yu\xi_\psi\|^2\\
		&\le
		\cos^2\theta\,\|yu\xi_\psi\|^2
		+(1-\cos^2\theta)\|Pyu\xi_\psi\|^2\\
		&\le
		\cos^2\theta\,(m+\varepsilon)^2
		+(1-\cos^2\theta)
		\|P_{\mathcal F_q(K)}u\xi_\psi\|^2.
	\end{aligned}
	\]
	Choosing $\varepsilon$ sufficiently small that
	$(m-2\varepsilon)^2>\cos^2\theta\,(m+\varepsilon)^2$
	gives the required uniform lower bound.
\end{proof}

We prove the second quantization estimate used above.

\begin{lem}\label{lem-deformation-inequality}
	Let $T\colon H\to H$ be a contraction and $\theta\in(0,\pi/2)$. Put
	\[
	F:=\mathcal F_q\bigl(R_\theta(1\oplus T)R_{-\theta}\bigr)
	\in\B(\mathcal F_q(H\oplus H)).
	\]
	Suppose that a closed subspace $K\subset H$ satisfies
	$TK\subset K$ and $T|_{K^\perp}=0$.
	Then, for every $\eta\in\mathcal F_q(H\oplus K)$, we have
	\[
	\|F\eta\|^2
	\leq
	\cos^2\theta\,\|\eta\|^2
	+
	(1-\cos^2\theta)
	\|P_{\mathcal F_q(K\oplus K)}\eta\|^2.
	\]
\end{lem}

\begin{proof}
	Write
	\[
	H\oplus K
	=
	(K\oplus K)\oplus(K^\perp\oplus0)
	=:
	\widetilde K\oplus K_1^\perp.
	\]
	Since $TK\subset K$, the operator
	$R_\theta(1\oplus T)R_{-\theta}$ leaves $\widetilde K$ invariant.
	Moreover, for $h\in K^\perp$,
	\[
	(1\oplus T)R_{-\theta}(h\oplus0)
	=
	\cos\theta\,h\oplus0,
	\]
	because $T|_{K^\perp}=0$. Hence, on $H\oplus K$, we may factor
	\[
	R_\theta(1\oplus T)R_{-\theta}
	=
	C D,\qquad \text{with }D=1_{\widetilde K}\oplus
	\cos\theta\,1_{K_1^\perp},
	\]
	where $C\colon H\oplus K\to H\oplus H$ agrees with
	$R_\theta(1\oplus T)R_{-\theta}$ on $\widetilde K$ and with
	$R_\theta$ on $K_1^\perp$.
	These restrictions have orthogonal ranges, so $C$ is a contraction.
	It then follows that
	\[
	\|F\eta\|
	\leq
	\|\mathcal F_q(D)\eta\|.
	\]
	The operator $\mathcal F_q(D)$ is the identity on
	$\mathcal F_q(\widetilde K)$ and has norm at most $\cos\theta$
	on its orthogonal complement in $\mathcal F_q(H\oplus K)$.
	Indeed, in the orthogonal decomposition according to the number of
	$K_1^\perp$-legs, it acts by $(\cos\theta)^k$ on the subspace
	with $k$ such legs. Thus,
	\[
	\begin{aligned}
		\|F\eta\|^2
		&\leq
		\|P_{\mathcal F_q(\widetilde K)}\eta\|^2
		+
		\cos^2\theta\,
		\|(1-P_{\mathcal F_q(\widetilde K)})\eta\|^2 \\
		&=
		\cos^2\theta\,\|\eta\|^2
		+
		(1-\cos^2\theta)
		\|P_{\mathcal F_q(\widetilde K)}\eta\|^2.
	\end{aligned}
	\]
	Since $\widetilde K=K\oplus K$, this is the desired estimate.
\end{proof}

\begin{proof}[{\bf Proof of Theorem~\ref{thmB} in case~(a)}]
	In this case, condition~(2) is equivalent to complete atomicity of $A$.
	
	Assume first that (1) holds. If $A$ is not completely atomic, then,
	passing to a nonzero diffuse direct summand, we may assume that
	$A$ is diffuse.
	Since $A$ is finite with expectation, choose a faithful normal
	state $\psi$ on $M_q$ such that $p\in (M_q)_\psi$ and
	$A\subset p(M_q)_\psi p$.
	We extend $\psi$ to $\widetilde M_q$ via the canonical conditional expectation.
	By uniform convergence and
	Lemma~\ref{lem-uniform-convergence-intertwine}, there exist
	$\theta\in(0,\pi/2)$ and a nonzero element
	$z\in\alpha_\theta(p)\widetilde M_qp$ such that
	\[
	\alpha_\theta(a)z=za
	\qquad(a\in A).
	\]
	Lemma~\ref{lem-intertwiner-inequality} provides a finite-dimensional subspace $K\subset H$
	and a constant $c>0$ such that
	\[
	\inf_{u\in\mathcal U(A)}
	\|P_{\mathcal F_q(K)}u\xi_\psi\|\ge c.
	\]
	On the other hand, uniform convergence and
	Lemma~\ref{lem-deformation-uniform-length1} give $N\in\mathbb N$ such that
	\[
	\sup_{u\in\mathcal U(A)}
	\|P_{> N}u\xi_\psi\|<\frac c2.
	\]
	Thus the finite-rank projection
	$Q:=P_{\mathcal F_q(K)}P_{\le N}$ satisfies
	\[
	\inf_{u\in\mathcal U(A)}\|Qu\xi_\psi\|\ge\frac c2.
	\]
	Since $A$ is diffuse, choose a sequence
	$(u_n)\subset\mathcal U(A)$ converging weakly to zero.
	Then $u_n\xi_\psi\to0$ weakly, so
	$\|Qu_n\xi_\psi\|\to0$ as $n\to \infty$ because $Q$ has finite rank.
	This contradiction proves (2).
	
	Conversely, assume that (2) holds.
	Since $A$ is finite and completely atomic, its unit ball is compact
	in the strong$^*$ topology. In particular, $(A)_1\Omega$ is
	norm compact.
	Let $\varphi$ be the vacuum state on $\widetilde M_q$.
	The strong convergence $V_\theta\to1$ is uniform on
	$(A)_1\Omega$, and $V_\theta\Omega=\Omega$. Hence
	\[
	\sup_{a\in(A)_1}\|\alpha_\theta(a)-a\|_{\varphi}
	=
	\sup_{a\in(A)_1}\|(V_\theta-1)a\Omega\|
	\longrightarrow0
	\qquad(\theta\to0).
	\]
	Since $(A)_1$ is self-adjoint, this implies (1).
\end{proof}

\subsection{Proof of Theorem~\ref{thmB}: case (b)}

We use the following doubling argument, which goes back to Popa
\cite{Po03}; see also \cite{HR10,Ma16a} for later variants.
The symmetry $\beta$ allows us to extend intertwiners from small
to arbitrary deformation parameters.

\begin{lem}\label{lem-uniform}
	Let $M\subset\widetilde M$ be an inclusion of $\sigma$-finite
	von Neumann algebras with expectation, equipped with an s-malleable
	deformation $(\alpha_\theta,\beta)$.
	Let $p\in M$ be a nonzero projection and
	$A\subset pMp$ a finite von Neumann subalgebra with expectation.
	Assume that $\alpha$ converges uniformly on the unit ball of $A$ and that
	\[
	A'\cap p\widetilde Mp\subset pMp
	\quad\text{or}\quad
	(A'\cap pMp)'\cap p\widetilde Mp\subset pMp.
	\]
	Then, for every $\theta\in\mathbb R$, there exists a nonzero element
	$w\in\alpha_\theta(p)\widetilde Mp$ such that
	\[
	\alpha_\theta(a)w=wa
	\qquad(a\in A).
	\]
\end{lem}

\begin{proof}
	Choosing a faithful normal state on $\widetilde M$ whose centralizer
	contains $A$, Lemma~\ref{lem-uniform-convergence-intertwine} gives
	the conclusion for all sufficiently small $|\theta|$.
	
	Put $P:=A'\cap pMp$, and let
	$w\in\alpha_s(p)\widetilde Mp$ be a nonzero partial isometry satisfying
	$\alpha_s(a)w=wa$ for all $a\in A$.
	Under the first assumption, $w^*w\in P$, so
	$w\beta(w^*)\ne0$, as in
	\cite[proof of Theorem~4.3, Step~(2)]{HR10}.
	Under the second assumption, the support-projection argument in
	\cite[proof of Theorem~4.2, Step~3]{Ma16a} gives
	$d\in P$ such that $wd\beta(w^*)\ne0$.
	Indeed, the support projection used there belongs to
	$P'\cap p\widetilde Mp$, which is contained in $pMp$ by assumption.
	
	Taking $d=p$ in the first case, we obtain in either case a nonzero
	$\alpha_{2s}(A)$-$A$-intertwiner
	\[
	\alpha_s\bigl(wd\beta(w^*)\bigr)
	\in\alpha_{2s}(p)\widetilde Mp.
	\]
	For any prescribed $\theta\ne0$, start with $s=2^{-n}\theta$
	for sufficiently large $n$ and iterate this construction,
	taking polar parts at each step.
	The case $\theta=0$ is immediate.
\end{proof}

\begin{proof}[{\bf Proof of Theorem~\ref{thmB} in case~(b)}]
	Put $M:=M_q(G)$ and $\widetilde M:=\widetilde M_q(G)$.
	
	Assume first that (1) holds, and fix a nonzero projection
	$r\in A'\cap pMp$.
	The relative-commutant assumption and uniform convergence pass to
	$Ar\subset rMr$.
	Lemma~\ref{lem-uniform}, applied with $\theta=\pi/2$, gives
	\[
	Ar\preceq_{\widetilde M}\alpha_{\pi/2}(Ar)
	\subset\alpha_{\pi/2}(M).
	\]
	Since $\alpha_{\pi/2}(M)$ is the crossed product associated with
	the second copy of $H_{\mathbb R}$,
	Proposition~\ref{prop-two-components-intertwining} yields
	$Ar\preceq_M LG$. This proves (2).

	Conversely, assume that (2) holds.
	Choose a faithful normal state $\psi$ on $\widetilde M$ whose
	centralizer contains $A$.
	For every nonzero projection $r\in A'\cap pMp$, condition~(2)
	provides an integer $n\ge1$, a projection $f\in\M_n(LG)$,
	a normal $*$-homomorphism
	$\gamma\colon A\to f\M_n(LG)f$,
	and a nonzero partial isometry
	$v\in\M_{1,n}(rM)f$ such that
	\[
	av=v\gamma(a)
	\qquad(a\in A).
	\]
	Set $e:=vv^*\in A'\cap pMp$, so that $0\ne e\le r$, and put
	$d_\theta:=\alpha_\theta(v)v^*$.
	Since $\alpha_\theta$ fixes $LG$ pointwise, it fixes $\gamma(A)$, and hence
	\[
	\alpha_\theta(a)d_\theta=d_\theta a
	\qquad(a\in A).
	\]
	Then, for $a\in(A)_1$,
	\begin{align*}
		\|\alpha_\theta(ae)-ae\|_\psi
		&\le
		\|\alpha_\theta(a)e-ea\|_\psi
		+\|\alpha_\theta(e)-e\|_\psi\\
		&\le
		\|\alpha_\theta(a)d_\theta-d_\theta a\|_\psi
		+2\|d_\theta-e\|_\psi
		+\|\alpha_\theta(e)-e\|_\psi.
	\end{align*}
	Since $\alpha_\theta(a)d_\theta=d_\theta a$, while
	$d_\theta\to e$ and $\alpha_\theta(e)\to e$ strongly, the right-hand
	side converges to $0$ uniformly on $(A)_1$ as $\theta\to 0$.
	Since $(Ae)_1$ is self-adjoint, $\alpha$ converges uniformly on the unit ball of $Ae$.
	
	A maximal orthogonal family of such projections $e$ has sum $p$,
	since the same argument applies inside every nonzero projection
	$r\in A'\cap pMp$.
	Approximating $p$ by finite sums of these projections gives uniform
	convergence on the unit ball of $A$, proving (1).
\end{proof}

\section{General consequences of non-uniform convergence}
\label{General consequences of non-uniform convergence}

In this section, we introduce a hereditary $C^*$-algebra determined solely by the convergence of a malleable deformation and discuss its general properties.
In the case of $q$-Araki--Woods algebras, this algebra coincides with the length-compact algebra $\mathcal K_\ell(G)$.

\subsection{A hereditary $C^*$-algebra associated with a deformation}

Let $M\subset\widetilde M$ be an inclusion of $\sigma$-finite von Neumann
algebras with expectation $E_M$, and let
$\alpha\colon\mathbb R\curvearrowright\widetilde M$ be a malleable deformation in the sense of
Subsection~\ref{Malleable deformations}.
For $t\in\mathbb R$, let $V_t\in \mathcal{U}(L^2(\widetilde M))$ be the standard unitary implementing $\alpha_t$.
Let $e_M$ be the Jones projection associated with $E_M$, and consider the operator on $L^2(M)$ given by
\[
V_t^M:=e_MV_t|_{L^2(M)}=e_MV_te_M\in \B(L^2(M)).
\]

\begin{df}\label{def-Kalpha}\upshape
	We define the hereditary $C^*$-algebra associated with the deformation $\alpha$ by
	\[
	\mathcal K_\alpha
	:=
	\left\{
	T\in\B(L^2(M))
	\mathrel{}\middle|\mathrel{}
	\|T(1-V_t^M)\|\to0,\quad \|(1-V_t^M)T\|\to0
	\quad(t\to0)
	\right\}.
	\]
\end{df}

It is immediate that $\mathcal K_\alpha$ is a hereditary $C^*$-subalgebra of $\B(L^2(M))$.
For $q$-Araki--Woods algebras, we have the following.

\begin{lem}\label{lem-Kalpha-Kell}
	For the crossed product $M_q(G) = \Gamma_q(H_\R,U_t)''\rtimes_{\sigma^\pi} G$ of a $q$-Araki--Woods algebra as in Subsection\ref{$q$-Araki--Woods factors} and its malleable deformation $\alpha$, the algebra $\mathcal K_\alpha$ coincides with $\mathcal K_\ell(G)$.
\end{lem}

\begin{proof}
	By Lemma~\ref{lem-deformation-uniform-length0}, if $V_t$ is the second quantization of the rotation matrix, then
	\[
	V_t^M=\sum_{n\geq0}(\cos t)^n(P_n\otimes 1),\qquad
	1-V_t^M=\sum_{n\geq0}\left(1-(\cos t)^n\right)(P_n\otimes 1).
	\]
	The conclusion follows directly from these formulas.
\end{proof}

The following lemma and proposition will be used in the study of a weak mixing-type property and solidity.

\begin{lem}\label{lem-deformation-weakly-mixing}
	Let $p\in M$ be a projection and $A\subset pMp$ a finite von Neumann subalgebra with expectation.
	Assume that $\alpha$ does not converge uniformly on $(A)_1$.
	Then there exist $t_n \to 0$ and a sequence of unitaries $u_n \in A$ such that
	\[
	b_n:=u_n-E_M(\alpha_{t_n}(u_n))\in(M)_2
	\]
	satisfies the following properties.
	Here $\omega$ is a free ultrafilter on $\N$.
	\begin{enumerate}[(i)]
		\item The element $B:=(b_n)_\omega$ defines a nonzero element of $pM^\omega p$.
		
		\item For all $x,y\in M$,
		\[
		xBy = (xu_ny-E_M(\alpha_{t_n}(xu_ny)) )_\omega .
		\]
		
		\item For every $T\in\mathcal K_\alpha$,
		$Tb_n\to0$ $(n\to \infty)$
		strongly in $\B(L^2(M))$.
	\end{enumerate}
\end{lem}

\begin{proof}
	Choose a faithful normal state $\psi\in \widetilde{M}_\ast$ satisfying $\psi\circ E_M = \psi$.
	By assumption, there exist $\delta>0$, $t_n\to0$, and $u_n\in\mathcal U(A)$ such that
	\[
	\|\alpha_{t_n}(u_n)-u_n\|_\psi\geq\delta,\quad \text{for all }n.
	\]
	Define $b_n$ as in the statement.

	Fix any $\xi\in L^2(M)$ and define
	\[
	\beta_{n,\xi}
	:=
	u_n\xi-e_MV_{t_n}u_n\xi
	=(1-V^M_{t_n})(u_n\xi).
	\]
	Then $\|\beta_{n,\xi}-b_n\xi\|
	\leq
	\|V_{t_n}^*\xi-\xi\|$.
	Therefore, for every $T\in\mathcal K_\alpha$,
	\[
	\begin{split}
		\|Tb_n\xi\|
		&\leq
		\|T\beta_{n,\xi}\|
		+
		\|T\|\|\beta_{n,\xi}-b_n\xi\|\\
		&\leq
		\|T(1-V^M_{t_n})\|\|\xi\|
		+
		\|T\|\|V_{t_n}^*\xi-\xi\|
		\longrightarrow0.
	\end{split}
	\]
	This proves (iii).
	
	For simplicity, set
	$\xi:=\xi_\psi$ and
	$\rho_n:=\|V_{t_n}\xi-\xi\|$.
	Define $\beta_{n,\xi}$ as in the proof of (iii), so that 
	$\|\beta_{n,\xi}-b_n\xi\|\leq\rho_n$.
	Moreover,
	\begin{align*}
		\|V_{t_n}u_n\xi-u_n\xi\|
		&\geq
		\|\alpha_{t_n}(u_n)\xi-u_n\xi\|-\rho_n
		\geq
		\delta-\rho_n;\\
		\|V_{t_n}u_n\xi-u_n\xi\|^2
		&= 2\Re \langle (1-V^M_{t_n})u_n\xi,u_n\xi\rangle\\
		&\leq
		2\|(1-V^M_{t_n})u_n\xi\|\|u_n\xi\|
		\leq
		2\|\beta_{n,\xi}\|
	\end{align*}
	and hence
	$\|\beta_{n,\xi}\|
	\geq
	\frac12(\delta-\rho_n)^2$.
	Thus, we have
	\[
	\|b_n\xi\|
	\geq
	\frac12(\delta-\rho_n)^2-\rho_n.
	\]
	Since $\rho_n\to0$, we obtain
	$\|b_n\xi\|\geq\frac{\delta^2}{4}$ for all sufficiently large $n$.
	
	Next, since $A$ is finite with expectation,
	$U:=(u_n)_\omega\in A^\omega\subset
	p(M^\omega)p$ is well defined.
	Since $\psi\circ \alpha_{t_n}\to \psi$ in norm, it follows that there exists an automorphism
	\[
	\Theta\colon
	\widetilde M^\omega\ni(x_n)_\omega
	\longmapsto
	(\alpha_{t_n}(x_n))_\omega
	\in\widetilde M^\omega.
	\]
	Hence
	$\Theta(U)=
	(\alpha_{t_n}(u_n))_\omega
	\in\widetilde M^\omega$, and
	\[
	E_{M^\omega}(\Theta(U))
	=
	(E_M(\alpha_{t_n}(u_n)))_\omega
	\in M^\omega.
	\]
	Therefore,
	$U-E_{M^\omega}(\Theta(U))
	\in M^\omega$, and $(b_n)_n$ is a representing sequence for this element.
	Hence $(b_n)_n$ defines an element of $M^\omega$, which is nonzero by the inequality obtained in the preceding paragraph.
	This proves (i).
	Since $\Theta$ restricts to the identity on $\widetilde M$, (ii) also follows.
\end{proof}

\begin{prop}\label{prop-solid}
	Let $p\in M$ be a projection and $A\subset pMp$ a finite von Neumann subalgebra with expectation, and put $P:=A'\cap pMp$.
	If $\alpha$ does not converge uniformly on $(A)_1$, then there exist a nonzero projection $r\in P'\cap pMp$ and a u.c.p.\ map $\Psi \colon \B(L^2(M))\to \B(rL^2(M))$ satisfying the following conditions:
	\begin{enumerate}[(i)]
		\item $\Psi|_M$ defines a normal map on $M$ and satisfies $\Psi(x)=xr$ $(x\in P)$;
		\item $\Psi(x') = x'r$ $(x'\in M')$;
		\item $\Psi(\mathcal K_\alpha)=0$.
	\end{enumerate}
\end{prop}
\begin{proof}
	By assumption, we may take $B=(b_n)_\omega\in pM^\omega p$ as in Lemma~\ref{lem-deformation-weakly-mixing}. By (ii) of the lemma, we have $B\in P'\cap pM^\omega p$.
	Following the argument of \cite[Theorem 4.1]{HI15b}, we modify $B$ so that $E_\omega(B^*B)$ becomes a projection.
	First,
	\[
	y:=E_\omega(B^*B) = \lim_{n\to \omega} b_n^*b_n \in P'\cap pMp
	\]
	is nonzero, so let $0\neq r\in P'\cap pMp$ be the spectral projection corresponding to $[\frac12\|y\|,\|y\|]$.
	Putting $c:=(yr)^{-1/2}\in P'\cap pMp$, we have
	\[
	E_\omega( (Bc)^*Bc )
	=cE_\omega(B^*B)c
	=cyc
	=r\neq 0.
	\]
	Now set $\widetilde{b}_n := b_nc$ and $\widetilde{B}:=Bc=(\widetilde{b}_n)_\omega \in P'\cap pM^\omega p$, and define a c.c.p.\ map
	\[
	\Psi\colon \B(L^2(M))\to \B(L^2(M));\quad T\mapsto 
	\sigma\text{-weak}\lim_{n\to\omega} \widetilde{b}_n^*T\widetilde{b}_n.
	\]
	Since $\Psi(1)=r$, by replacing the range with $r\B(L^2(M))r = \B(rL^2(M))$, we obtain a u.c.p.\ map.
	Moreover, $\Psi$ is multiplication by $r$ on $M'$, while on $M$ it is given by
	\[
	\Psi(x) = E_\omega(\widetilde{B}^* x \widetilde{B}),\quad x\in M,
	\]
	and hence is normal.
	Since $\widetilde{B}\in P'$, $\Psi$ is multiplication by $r$ on $P$.
	Finally, by Lemma~\ref{lem-deformation-weakly-mixing}(iii), for every $T\in\mathcal K_\alpha$, we have $T\widetilde{b}_n\to0$
	strongly in $\B(L^2(M))$.
	It follows that
	$\Psi(\mathcal K_\alpha)=0$.
\end{proof}

\subsection{S-malleable deformations and stable normalizers}

We first recall a useful subclass of stable normalizers.
Let $A\subset N$ be an inclusion of \emph{finite} von Neumann algebras, and let
$E_{\mathcal Z}\colon A\to\mathcal Z(A)$ denote the center-valued trace.
For $x\in s\mathcal N_N(A)$, let $z_x^r,z_x^l\in\mathcal Z(A)$ be the
support projections of $E_{\mathcal Z}(x^*x)$ and
$E_{\mathcal Z}(xx^*)$, respectively, and put
\[
s\mathcal N^0_N(A)
:=
\left\{
x\in s\mathcal N_N(A)
\mathrel{}\middle|\mathrel{}
\begin{array}{l}
	\text{there exists }\delta>0\text{ such that}\\
	E_{\mathcal Z}(x^*x)\ge\delta z_x^r,\quad
	E_{\mathcal Z}(xx^*)\ge\delta z_x^l
\end{array}
\right\}.
\]
As observed just before \cite[Proposition~3.6]{BHV15},
$s\mathcal N^0_N(A)$ and $s\mathcal N_N(A)$ generate the same von Neumann
algebra.

\medskip

Now let $M\subset^{E_M}\widetilde M$ be an inclusion of $\sigma$-finite
von Neumann algebras equipped with an s-malleable deformation
$(\alpha,\beta)$.
We further assume that $\widetilde{M}$ admits a semifinite trace $\Tr$ such that $\Tr\circ E_M=\Tr$ and that $\alpha,\beta$ preserve $\Tr$.

Then $\beta\circ E_M \circ \beta^{-1}$ is also a $\Tr$-preserving conditional expectation, and hence, by uniqueness, it coincides with $E_M$.
Thus $E_M$ commutes with $\beta$.
It follows that, for $a\in M$,
\[
E_M(\alpha_t(a))
=
\beta\circ E_M(\alpha_t(a))
=
E_M(\alpha_{-t}(\beta(a)))
=
E_M(\alpha_{-t}(a)),
\]
and hence $e_M V_t e_M = e_M V_{-t} e_M$.
Since $V_{-t}=V_t^*$, this shows that $V_t^M$ is self-adjoint for all $t\in \R$.

\medskip
The following proposition is a variant of \cite[Proposition 3.7]{BHV15}.
We use a formulation suited to applying the AO-type
property established in
Section~\ref{An AO-type property for length projections}.
Below, $M^{\op}$ denotes the opposite von Neumann algebra,
and $x^{\op}$ denotes the element corresponding to $x\in M$.
We write $J$ and $J^{\op}$ for the modular conjugations of
$M$ and $M^{\op}$, respectively.
The right action of $M\ovt M^{\op}$ on $L^2(M)\otimes L^2(M^{\op})$ is given by
\[
\xi T=(J\otimes J^{\op})T^*(J\otimes J^{\op})\xi
\qquad
(T\in M\ovt M^\op,\ \xi\in L^2(M)\otimes L^2(M^{\op})).
\]
\begin{prop}\label{prop-deformation-stable-normalizer}
	Assume that $M$ has the complete metric approximation property
	(CMAP).
	Let $p\in M$ be a projection with $\Tr(p)<\infty$, and let $A\subset pMp$ be an amenable von Neumann subalgebra.
	If the deformation $\alpha$ does not converge uniformly on $(A)_1$, then there exists a state $\Omega$ on $\B(L^2(M)\otimes L^2(M^{\op}))$ satisfying the following conditions:
	\begin{enumerate}[(i)]
		\item
		$\Omega(\mathcal K_\alpha\otimes1)=0$ and $\Omega(p\otimes 1)=1$;
		
		\item
		$\Omega|_{M\otimes\C}$ and
		$\Omega|_{JMJ\otimes\C}$ are normal;
		
		\item
		for every partial isometry
		$v\in s\mathcal N^0_{pMp}(A)$, there exist
		\[
		T(v,k)\in
		\bigl(pMp\otimes_{\alg}(pMp)^{\op}\bigr)_1
		\qquad (k\geq1)
		\]
		such that, if $L_v,R_{v,k}\in \B(L^2(M)\otimes L^2(M^{\op}))$ are defined by
		\[
		L_v\xi:=(v\otimes1)\xi
		\quad\text{and}\quad
		R_{v,k}\xi:=(1\otimes v^{\op})\xi T(v,k),
		\qquad
		\text{for}\quad \xi\in L^2(M)\otimes L^2(M^{\op}),
		\]
		then
		\[
		\lim_{k\to\infty}
		\Omega\bigl(
		(L_v-R_{v,k})^*(L_v-R_{v,k})
		\bigr)
		=0
		=
		\lim_{k\to\infty}
		\Omega\bigl(
		(L_v-R_{v,k})(L_v-R_{v,k})^*
		\bigr).
		\]
	\end{enumerate}
	Here, the product $\xi T(v,k)$ in the definition of $R_{v,k}$
	is taken with respect to the right action defined above.
\end{prop}

\begin{proof}
	Rescaling $\Tr$, we may assume that $\Tr(p)=1$. We first apply \cite[Proposition 3.6]{BHV15} to $A\subset pMp$.
	Then there exists a net of unit vectors
	$\eta_j\in L^2(A\overline\otimes A^{\op})^+$
	satisfying the following conditions:
	\begin{enumerate}
		\item[$\rm (a)$] $\|(a\otimes1)\eta_j
		-(1\otimes a^{\op})\eta_j\|
		\to0$,
		\quad $a\in A$;
		\item[$\rm (b)$] $\langle\eta_j,(x\otimes1) \eta_j\rangle
		\to \Tr(x)$,
		\quad $x\in pMp$;
		\item[$\rm (c)$] for every partial isometry $v\in s\mathcal N^0_{pMp}(A)$, there exist
		$
		T(v,k)
		\in
		\bigl(
		pMp\otimes_{\alg}(pMp)^{\op}
		\bigr)_1
		$
		such that, if $L_v,R_{v,k}\in \B(L^2(M)\otimes L^2(M^{\op}))$ are defined as in the statement of the proposition, then
		\[
		\lim_{k\to\infty}\limsup_j
		\|L_v\eta_j -R_{v,k}\eta_j\|=0
		=\lim_{k\to\infty}\limsup_j
		\|L_v^*\eta_j - R_{v,k}^*\eta_j\|.
		\]
		
	\end{enumerate}
	We note that, in (b), for every $x\in \widetilde M$, we also have
	\[
	\langle\eta_j,(x\otimes1)\eta_j\rangle
	=\langle\eta_j,(E_A(pxp)\otimes1)\eta_j\rangle
	\to \Tr(E_A(pxp)) = \Tr(xp),
	\]
	where $E_A:p\widetilde Mp\to A$ is the trace-preserving
	conditional expectation. 
	For each $t\in \R$, we put $D_t:=1-V_t^M$.
	
	\begin{claim}
		There exist $c>0$ and $t_m\to0$ such that
		\[
		\limsup_j\|	(D_{t_m}\otimes1)\eta_j\|>c,
		\qquad \text{for all}\quad  m\geq1.
		\]
	\end{claim}
	\begin{proof}[{\bf Proof of Claim}]
		Since $\alpha$ does not converge uniformly on $(A)_1$, the first part of the proof of \cite[Proposition 3.7]{BHV15} shows that there exist $\varepsilon>0$ and $t_m\to0$ such that
		\[
		\limsup_j\|	(e_M^\perp V_{t_m}\otimes1)\eta_j\|>\varepsilon
		\qquad \text{for all}\quad  m\geq1,
		\]
		where $e_M^\perp = 1-e_M$.
		It is therefore enough to show that, for every $\eta\in L^2(M)\otimes L^2(M^{\op})$,
		\[
		\|
		(e_M^\perp V_t\otimes1)\eta
		\|^2
		\leq
		2\|(D_t\otimes1)\eta\|\,\|\eta\|.
		\]
		Since $V_t^M$ is a self-adjoint contraction,
		\[
		1-(V_t^M)^2
		=
		(1-V_t^M)(1+V_t^M)
		\leq
		2(1-V_t^M)
		=
		2D_t.
		\]
		Hence, for every $\eta\in L^2(M)\otimes L^2(M^{\op})$,
		\begin{align*}
			\|
			(e_M^\perp V_t\otimes1)\eta
			\|^2
			&=
			\|(V_t\otimes1)\eta\|^2
			-
			\|
			(e_M V_t\otimes1)\eta
			\|^2\\
			&=
			\langle
			((1-(V_t^M)^2)\otimes1)\eta,\eta
			\rangle\\
			&\leq
			2\langle(D_t\otimes1)\eta,\eta\rangle
			\leq
			2\|(D_t\otimes1)\eta\|\,\|\eta\|.
		\end{align*}
		This proves the claim.
	\end{proof}

	For every partial isometry $v\in s\mathcal N^0_{pMp}(A)$ and $k\geq1$, put
	\[
	a_{v,k}
	:=
	\limsup_j
	\|L_v\eta_j-R_{v,k}\eta_j\|,
	\qquad
	b_{v,k}
	:=
	\limsup_j
	\|L_v^*\eta_j-R_{v,k}^*\eta_j\|.
	\]
	By (c), both quantities converge to $0$ as $k\to \infty$.
	
	Below, we use the same notation $L_v$ and $R_{v,k}$ for the operators
	defined by the same formulas on
	$L^2(\widetilde M)\otimes L^2(M^{\op})$.
	
	\begin{claim}
		For every finite set
		$\mathcal F
		\subset
		s\mathcal N^0_{pMp}(A)$ consisting of partial isometries,
		$n\in \mathbb N$,
		every finite set $\mathcal G\subset M$,
		and every $\delta>0$, there exists $t_0>0$ such that, putting
		$\eta_j^t:=(V_t\otimes1)\eta_j$,
		for every $|t|<t_0$,
		\[
		\limsup_j
		\|L_v\eta_j^t-R_{v,k}\eta_j^t\|
		\leq
		a_{v,k}+\delta,
		\qquad
		v\in\mathcal F,\ 1\leq k\leq n,
		\]
		\[
		\limsup_j
		\|L_v^*\eta_j^t-R_{v,k}^*\eta_j^t\|
		\leq
		b_{v,k}+\delta,
		\qquad
		v\in\mathcal F,\ 1\leq k\leq n,
		\]
		\[
		\limsup_j
		\|(x\otimes1)\eta_j^t\|
		\leq
		\|xp\|_{\Tr}+\delta,
		\qquad
		x\in\mathcal G.
		\]
	\end{claim}
	\begin{proof}[{\bf Proof of Claim}]
		As mentioned just before equation (3.8) in the proof of \cite[Proposition 3.7]{BHV15}, for every $x\in M$ and $T\in M\otimes_{\alg}M^{\op}$,
		\[
		\lim_{t\to0}\limsup_j
		\left\|
		(x\otimes1)(V_t\otimes1)\eta_j
		-
		(V_t\otimes1)(x\otimes1)\eta_j
		\right\|
		=0;
		\]
		\[
		\lim_{t\to0}\limsup_j
		\left\|
		(V_t\otimes1)(\eta_jT)
		-
		[(V_t\otimes1)\eta_j]T
		\right\|
		=0.
		\]
		Hence there exists $t_0>0$ such that, whenever $|t|<t_0$, for every $y\in \mathcal G\cup \mathcal F\cup \mathcal F^*$, $v\in \mathcal F$, and $k=1,\ldots,n$,
		\[
		\limsup_j
		\left\|
		(y\otimes1)(V_t\otimes1)\eta_j
		-
		(V_t\otimes1)(y\otimes1)\eta_j
		\right\|
		<\frac\delta2;
		\]
		\[
		\limsup_j
		\left\|
		(V_t\otimes1)(\eta_jT(v,k))
		-
		[(V_t\otimes1)\eta_j]T(v,k)
		\right\|
		<\frac\delta2.
		\]
		We may also assume that the second inequality holds after replacing $T(v,k)$ by $T(v,k)^*$.
		
		Applying the first condition to $x\in \mathcal G$, for $|t|<t_0$ we obtain
		\begin{align*}
		\limsup_j\|(x\otimes1)\eta_j^t\|
		&\leq
		\limsup_j\|(V_t\otimes 1)(x\otimes1)\eta_j\|+\frac{\delta}{2}\\
		&=
		\limsup_j\|(x\otimes1)\eta_j\|+\frac{\delta}{2}
		=
		\|xp\|_{\Tr}+\frac{\delta}{2}.
		\end{align*}
		which gives the assertion for $x\in \mathcal G$.
		
		We next consider $v\in \mathcal F$.
		For every $v\in\mathcal F$ and $k=1,\ldots,n$, we estimate $\|L_v\eta_j^t-R_{v,k}\eta_j^t\|$ from above by
		\[
		\left\|
		L_v(V_t\otimes1)\eta_j
		-
		(V_t\otimes1)L_v \eta_j
		\right\|
		+
		\left\|
		(V_t\otimes1)
		\bigl(
		L_v\eta_j-R_{v,k}\eta_j
		\bigr)
		\right\|
		+
		\left\|
		(V_t\otimes1)R_{v,k}\eta_j
		-
		R_{v,k}(V_t\otimes1)\eta_j
		\right\|.
		\]
		After taking $\limsup_j$, the second term gives $a_{v,k}$.
		The first and third terms are bounded by $\delta/2$ by the conditions above.
		This proves the assertion for $v\in \mathcal F$.
		The estimate for the adjoints follows similarly.
	\end{proof}
	
	We now regard the quadruples $i=(\mathcal F,n,\mathcal G,\delta)$ appearing in the preceding claim as a directed set and construct a net.
	
	\begin{claim}
		There exists a net $(\zeta_i)_i$, indexed by quadruples $i=(\mathcal F,n,\mathcal G,\delta)$, such that $c<\|\zeta_i\|\leq 2$ and for any $v\in\mathcal F$,  $1\leq k\leq n$, $x\in\mathcal G$, 
		\[
		\|(L_v-R_{v,k})\zeta_i\|
		<
		2a_{v,k}+4\delta,
		\quad 
		\|(L_v-R_{v,k})^*\zeta_i\|
		<
		2b_{v,k}+4\delta,
		\]
		\[
		\|(x\otimes1)\zeta_i\|
		<
		2\|xp\|_{\Tr}+4\delta.
		\]
		Moreover, $(\zeta_i)_i$ can be chosen as a subnet of $((D_{t_m}\otimes1)\eta_j)_{(m,j)}$.
	\end{claim}
	\begin{proof}[{\bf Proof of Claim}]
		For a quadruple $i=(\mathcal F,n,\mathcal G,\delta)$, apply the preceding claim and choose $t_0>0$.
		Fix $t_m$ such that $|t_m|<t_0$.
		We can then choose an arbitrarily large $j$ such that
		\begin{itemize}
			\item by the first claim,
			$\|(D_{t_m}\otimes1)\eta_j\|>c$;
			
			\item by the preceding claim, for every $v\in\mathcal F$, $1\leq k\leq n$, $x\in \mathcal G$, and $t\in\{0,t_m\}$,
			\[
			\|L_v\eta_j^{t}-R_{v,k}\eta_j^{t}\|
			<
			a_{v,k}+2\delta,
			\quad
			\|L_v^*\eta_j^{t}-R_{v,k}^*\eta_j^{t}\|
			<
			b_{v,k}+2\delta,
			\quad 
			\|(x\otimes1)\eta_j^{t}\|
			<
			\|xp\|_{2,\Tr}+2\delta.
			\]
		\end{itemize}
		For such $m,j$, put
		$\zeta_i
		:=
		(D_{t_m}\otimes1)\eta_j$.
		Since $e_M\otimes1$ commutes with $L_v$ and $R_{v,k}$, combining the inequalities above for $t=0,t_m$ gives
		\begin{align*}
			\|(L_v-R_{v,k})\zeta_i\|
			\leq
			\|(L_v-R_{v,k})\eta_j\|
			+
			\|(L_v-R_{v,k})\eta_j^{t_m}\|
			<
			2a_{v,k}+4\delta,
			\qquad
			v\in\mathcal F,\ 1\leq k\leq n.
		\end{align*}
		This gives the first assertion.
		The second and third assertions follow in the same way.
		This yields the desired net.
		It is clear from the construction that it can be chosen as a subnet of $((D_{t_m}\otimes1)\eta_j)_{(m,j)}$.
	\end{proof}
	
	The net $(\zeta_i)_i$ constructed in the last claim satisfies, for every $x\in M$, every partial isometry $v\in s\mathcal N^0_{pMp}(A)$, and every $k\in \N$,
	\[
	\limsup_i
	\|(L_v-R_{v,k})\zeta_i\|
	\leq
	2a_{v,k},\qquad 
	\limsup_i
	\|(L_v-R_{v,k})^*\zeta_i\|
	\leq
	2b_{v,k},
	\]
	\[
	\limsup_i
	\|(x\otimes1)\zeta_i\|
	\leq
	2\|xp\|_{\Tr}.
	\]
	Consider the vector functionals associated with $\zeta_i$ on $\B(L^2(M)\otimes L^2(M^{\op}))$. Passing to a subnet, we may assume that these vector functionals
	converge weak$^*$ to a positive functional $\Omega$.
	Since $c< \|\zeta_i \|\leq 2$, $\Omega$ is a nonzero positive functional.
	We show that its normalization satisfies all the conclusions of the proposition.
	
	Let $K\in\mathcal K_\alpha$.
	By the definition of $\mathcal K_\alpha$,
	$\|KD_{t_m}\|\to 0$.
	Since $\zeta_i $ is a subnet of $(D_{t_m}\otimes 1)\eta_j$, it follows that $\| (K\otimes 1)\zeta_i\|\to 0$, and hence $\Omega(\mathcal K_\alpha \otimes 1)=0$.
	Next, for each $x\in M$,
	\[
	\Omega(x^*x\otimes 1) = \lim_i \|(x\otimes 1)\zeta_i\|^2 \leq  4 \|xp\|_{\Tr}^2.
	\]
	In particular, $\Omega(p^\perp\otimes 1)=0$, and $\Omega$ is normal on $M\otimes 1$.
	Moreover, since $\zeta_i$ is invariant under $J\otimes J^{\op}$, $\Omega$ is also normal on $JMJ\otimes 1$.
	Finally, for every partial isometry $v\in s\mathcal N^0_{pMp}(A)$ and $k\geq1$,
	\[
	\Omega\bigl(
	(L_v-R_{v,k})^*(L_v-R_{v,k})
	\bigr)
	=
	\lim_i
	\|(L_v-R_{v,k})\zeta_i\|^2
	\leq
	4 a_{v,k}^2,
	\]
	which converges to $0$ as $k\to \infty$.
	Similarly, for $\Omega\bigl(
	(L_v-R_{v,k})(L_v-R_{v,k})^*
	\bigr)$, the same argument gives a bound in terms of $b_{v,k}$, which again converges to $0$.
	This proves all the assertions of the proposition.
\end{proof}

\section{A weak mixing-type property}

In this section, we prove Theorem~\ref{thmF}
and Corollary~\ref{corG}.

\subsection{Stable weak mixing and stable normalizers}

We introduce a weak mixing-type condition adapted to bounded sequences
arising from deformations, whose terms need not be unitary or belong
to the subalgebra under consideration.

\begin{df}\upshape\label{df-stable-weakly-mixing}
	Let $N\subset^{E_N}M$ be an inclusion of von Neumann algebras with expectation.
	Let $0\neq p\in N$ be a projection and let $A\subset pNp$ be a von Neumann subalgebra with expectation.
	We say that the inclusion $A\subset pNp$ is \emph{stably weakly mixing with respect to $N\subset M$}
	if the following condition holds.
	For any nonzero projections $ e\in A$ and 
	$z\in\mathcal Z((eAe)'\cap eNe)$, there exist a free ultrafilter $\omega$ and a bounded sequence
	$(b_n)_n\subset zNz$ satisfying the following conditions:
	\begin{enumerate}[(i)]
		\item
		the sequence $(b_n)_n$ represents a nonzero element
		$B:=(b_n)_\omega$ of $(zNz)^\omega$;
		
		\item $E_\omega(B^*B)\in((eAez)'\cap zNz)'$, and moreover,
		for every partial isometry $v\in M$,
		\[
		v^*v=e\quad \text{and}\quad veAezv^*\subset N
		\quad \Rightarrow\quad
		vBv^*\in N^\omega;
		\]
		
		\item
		for every $C\in N^\omega$ and $x\in M\ominus N:=\ker E_N$, we have $E_\omega(C^*xB)=0$.
	\end{enumerate}
\end{df}

Note that condition (ii) is automatic if the $b_n$ can be chosen from $eAez$.
The following lemma justifies the terminology ``weakly mixing''.

\begin{lem}\label{lem-stable-weakly-mixing}
	Let $N\subset^{E_N}M$ be an inclusion of von Neumann algebras with expectation, and
	let $\omega$ be a free ultrafilter.
	Let $(b_n)_n\subset N$ be a bounded sequence such that $B:=(b_n)_\omega\in N^\omega$ is well defined.
	Then the following conditions are equivalent.
	\begin{enumerate}
		\item
		For every $C\in N^\omega$ and $x\in M\ominus N$, we have
		$E_\omega(C^*xB)=0$.
		
		\item
		For every $x,y\in M\ominus N$, we have $E_{N^\omega} (xBy)=0$.
	\end{enumerate}
	Here $E_{N^\omega}:M^\omega\to N^\omega$ is the conditional expectation given by $E_{N^\omega}((x_n)_\omega)=(E_N(x_n))_\omega$ for $(x_n)_\omega\in M^\omega$.	
\end{lem}

\begin{proof}
	(1)$\Rightarrow$(2).
	Fix arbitrary $x,y\in M\ominus N$, and set $c_n:=E_N(xb_ny)\in N$ and
	$C:=E_{N^\omega}(xB y) =(c_n)_\omega$.
	Since $C\in N^\omega$ and $E_\omega\circ E_{N^\omega} = E_N\circ E_\omega$, we have
	\[
	E_\omega(C^*C)
	=
	E_\omega\bigl(E_{N^\omega}(C^*xB y)\bigr)
	=
	E_N\bigl(E_\omega(C^*xB)y\bigr)
	=
	0.
	\]
	This means $C=0$.
	
	\medskip
	
	(2)$\Rightarrow$(1).	
	First note that for every $x\in M\ominus N$ and $y\in M$, we also have
	$E_{N^\omega} (xBy) =0$.
	Take arbitrary $C=(c_n)_\omega\in N^\omega$ and $x\in M\ominus N$, and
	set $z:=E_\omega(C^*xB)\in M$. For every $y\in M$, we have
	\[
	E_N(zy)
	=
	E_N\bigl(E_\omega(C^*xB)y\bigr)
	=
	E_\omega\bigl(E_{N^\omega}(C^*xB y)\bigr)
	=
	E_\omega\bigl(C^*E_{N^\omega}(xB y)\bigr)=0.
	\]
	In particular, taking $y=z^*$ gives $z=0$.
\end{proof}
\begin{rem}\upshape
	In Lemma~\ref{lem-stable-weakly-mixing}, if the bounded sequence $(b_n)_n$ satisfies (2) for every free ultrafilter $\omega$, then
	for every $x,y\in M\ominus N$, we have
	\[
	E_N(xb_ny)\longrightarrow0
	\quad\text{$\ast$-strongly}\quad (n\to\infty).
	\]
	When the $b_n$ are unitaries, this agrees with the usual weak mixing condition.
\end{rem}

We will use this notion in the following form.

\begin{prop}\label{prop-stably-weakly-mixing-stable-normalizer}
	Let $N\subset^{E_N}M$ be an inclusion of von Neumann algebras with expectation.
	Let $0\neq p\in N$ be a projection and let $A\subset pNp$ be a von Neumann subalgebra with expectation. Assume that $A$ is stably weakly mixing with respect to $N\subset M$.
	
	If a partial isometry $v\in M$ satisfies $v^*v\in A$ and $vAv^*\subset N$, then $v \in N$.
	In particular, $s\mathcal N_{pMp}(A)''\subset pNp$.
\end{prop}
\begin{proof}
	We may assume that $v\neq0$.
	Set $e:=v^*v\in A$, and observe that $eAe\subset eNe$ is also
	stably weakly mixing with respect to $N\subset M$.
	Moreover, since $v^*v\in eAe$ and $v(eAe)v^*=vAv^*\subset N$,
	by replacing $p,A$ with $e,eAe$, we may assume that $v^*v=p$.
	Set
	\[
	Q:=A'\cap pNp,\qquad d:=v-E_N(v).
	\]
	We will show that $d=0$.
	Suppose, towards a contradiction, that $E_N(d^*d)\neq0$.
	
	\begin{claim}
		We have $E_N(d^*d)\in Q$, and if
		$z:=z_Q(E_N(d^*d))\in\mathcal Z(Q)$ denotes its central support projection in $Q$, then
		$v(Az)v^*\subset N$.
	\end{claim}
	\begin{proof}[{\bf Proof of Claim}]
		Let $\theta:=\Ad(v)$. Then for every $a\in A$,
		$da=\theta(a)d$, and this implies $d^*d\in A'\cap pMp$.
		Applying $E_N$ gives
		\[
		E_N(d^*d)\in A'\cap pNp=Q.
		\]
		Thus, $z:=z_Q(E_N(d^*d))\in\mathcal Z(Q)$ is well defined.
		
		Next, putting $z^\perp:=p-z$, we have
		\[
		E_N\bigl(z^\perp d^*dz^\perp\bigr)
		=z^\perp E_N(d^*d)z^\perp=0.
		\]
		This means that $dz^\perp =0$ and 
		$vz^\perp=E_N(v)z^\perp\in N$. It follows that
		\[
		vzv^*
		=
		vpv^*-(vz^\perp )(vz^\perp )^*
		\in N,
		\]
		and therefore $v(Az)v^*=(vAv^*)(vzv^*)\subset N$.	
	\end{proof}
	
	Apply Definition~\ref{df-stable-weakly-mixing} with $e=p$ and the projection $z$ from the claim, and choose
	$\omega$ and $0\neq B\in(zNz)^\omega$ satisfying
	conditions (i)--(iii).
	By the claim and condition (ii), we have $C:=vBv^*\in N^\omega$.
	Since $d=v-E_N(v)\in M\ominus N$, condition (iii) gives $E_\omega(C^*dB)=0$. 
	Set $y:=E_\omega(B^*B)\neq0$. We then compute that 
	\[
	vy
	=E_\omega(vB^*B)
	=E_\omega(C^*vB)
	=E_\omega(C^*E_N(v)B)\in N,
	\]
	and therefore $vy=E_N(vy)=E_N(v)y$ and $dy=0$.
	We obtain
	\[
	E_N(d^*d)y=E_N(d^*dy)=0.
	\]
	By condition (ii), $y$ commutes with $Qz$, hence with $Q$ because $y=zy$.
	Since $z$ is the central support of $E_N(d^*d)$ in $Q$,
	the equality $E_N(d^*d)y=0$ implies $zy=0$.
	This contradicts $zy=y\neq0$, proving that $d=0$ and hence $v\in N$.
	
	The final assertion follows from a standard argument using polar decomposition, since the partial isometries in $s\mathcal N_{pMp}(A)$ generate $s\mathcal N_{pMp}(A)''$.
\end{proof}

\subsection{Proofs of Theorem~\ref{thmF} and Corollary~\ref{corG}}

Theorem~\ref{thmF} follows immediately from the following theorem and Proposition~\ref{prop-stably-weakly-mixing-stable-normalizer}.

\begin{thm}\label{thm-relative-normalizer}
	Let $(H_\R,U_t)$ be a strongly continuous orthogonal representation of $\R$, and let
	$K_\R\subset H_\R$ be a $U_t$-invariant closed real subspace.
	Put $K:=K_\R\otimes_\R\C$ and 
	consider the natural inclusion
	\[
	M_K:=\Gamma_q(K_\R,U_t)''
	\subset
	\Gamma_q(H_\R,U_t)''=:M
	\]
	with expectation $E_K$.
	Let $0\neq p\in M_K$ be any projection and let
	$A\subset pM_Kp$ be any diffuse von Neumann subalgebra
	with expectation.
	Then $A$ is stably weakly mixing with respect to $M_K\subset M$.
\end{thm}

\begin{proof}
	We use the canonical malleable deformation $\alpha$ of the doubled
	inclusion $M\subset\widetilde M$ as in Subsection~\ref{Malleable deformations}, and write
	$E_M:\widetilde M\to M$ for the vacuum-preserving conditional
	expectation.

	Fix arbitrary nonzero projections $e\in A$ and $z\in \mathcal Z((eAe)'\cap eM_Ke )$.
	Since $A_z:=eAez$ is diffuse, there exists a diffuse abelian von Neumann subalgebra with expectation $A_0\subset A_z$. 
	By Theorem~\ref{thmB},
	the deformation $\alpha$ does not converge uniformly on $(A_0)_1$.
	Choose $t_n\to0$ and $u_n\in\mathcal U(A_0)$ witnessing this failure.
	Fix a free ultrafilter $\omega$.
	By Lemma~\ref{lem-deformation-weakly-mixing}, the sequence
	\[
	b_n:=u_n-E_M(\alpha_{t_n}(u_n))\in(M_K)_2
	\]
	satisfies the following conditions:
	\begin{enumerate}
		\item[$\rm(a)$] $B:=(b_n)_\omega$ defines a nonzero element of $zM^\omega z$;
		
		\item[$\rm(b)$] 
		$xBy = (xu_ny-E_M(\alpha_{t_n}(xu_ny)) )_\omega $, \quad for any $x,y\in M$;
		
		\item[$\rm(c)$]  
		$Tb_n\to0$ strongly in $\B(L^2(M))$ as $n\to\infty$,\quad for any $T\in\mathcal K_\alpha$.
	\end{enumerate}
	Noting that $b_n$ and $zb_nz$ define the same element in the ultraproduct, it suffices to show that $B=(b_n)_\omega$ satisfies conditions (i)--(iii) of stable weak mixing. First, (i) follows from (a) and $b_n\in M_K$. Next, (ii) follows immediately from (b).
	
	To verify (iii), fix $C=(c_n)_\omega\in M_K^\omega$
	and $x\in M\ominus M_K$.
	We will show that $E_\omega(C^*xB)=0$.	
	Since the map $x\mapsto E_\omega(C^*xB)$ is normal, it suffices to prove this when $x$ is a Wick polynomial.
	So assume that $x$ is a Wick polynomial satisfying $E_K(x)=0$. Then $x$ can be written as a finite linear combination of elements of the form $W(\xi)$, where $\xi$ is a simple tensor in $\mathcal F_{q}(H)\ominus\mathcal F_{q}(K)$.

	The following claim is proved using Lemma~\ref{lem-relative-wick-vanishing} below.
	
	\begin{claim}
		For any simple tensors $\eta,\eta'\in\mathcal F_{q}(H)$
		whose tensor factors are analytic vectors, we have
		\[
		\left\langle
		xW_r(\eta)b_n\Omega,\,
		W_r(\eta')c_n\Omega
		\right\rangle
		\longrightarrow0\qquad(n\to\infty).
		\]
	\end{claim}
	\begin{proof}[{\bf Proof of Claim}]
		By linearity, we may assume that $x=W(\xi)$ is a Wick word.
		Since $b_n\Omega,c_n\Omega\in\mathcal F_q(K)$, using the notation of Lemma~\ref{lem-relative-wick-vanishing}, we have
		\[		
		\left\langle
		xW_r(\eta)b_n\Omega,\,
		W_r(\eta')c_n\Omega
		\right\rangle
		=
		\left\langle
		T_{\xi,\eta,\eta'}b_n\Omega,\,
		c_n\Omega
		\right\rangle.
		\]
		By Lemmas~\ref{lem-relative-wick-vanishing} and
		\ref{lem-Kalpha-Kell}, we have
		$T_{\xi,\eta,\eta'}\in\mathcal K_\ell=\mathcal K_\alpha$.
		Condition~(c) gives $\|T_{\xi,\eta,\eta'}b_n\Omega\|\to0$ as $n\to\infty$.
		The boundedness of $(c_n)_n$ and the Cauchy--Schwarz inequality
		give the required convergence.
	\end{proof}
	
	In the setting of the claim, we have 
	\[
	\langle E_\omega(C^*xB)\eta,\eta'\rangle
	=
	\lim_{n\to \omega} \langle c_n^*xb_n\eta,\eta'\rangle
	=
	\lim_{n\to \omega}\left\langle
	xW_r(\eta)b_n\Omega,\,
	W_r(\eta')c_n\Omega
	\right\rangle
	=0.
	\]
	Since such $\eta$ and $\eta'$ span a dense subspace of $\mathcal F_q(H)$, it follows that $E_\omega(C^*xB)=0$. Thus (iii) also holds, and the proof is complete.
\end{proof}

\begin{lem}\label{lem-relative-wick-vanishing}
	Let $K_\R\subset H_\R$ be as in
	Theorem~\ref{thm-relative-normalizer}, and let $e_K$ be the
	orthogonal projection onto $\mathcal F_q(K)$.
	Let $\xi,\eta,\eta'$ be simple tensors in $\mathcal F_q(H)$
	such that all tensor factors of $\xi$ belong to $jH$, while
	those of $\eta,\eta'$ belong to $IjH$.
	Assume that $\xi\perp\mathcal F_q(K)$.
	Then
	\[
	T_{\xi,\eta,\eta'}
	:=
	e_KW_r(\eta')^*W(\xi)W_r(\eta)e_K
	\in\mathcal K_\ell.
	\]
\end{lem}

\begin{proof}
	We use the commutation argument from
	\cite[proof of Proposition~6]{SW16}, reformulated modulo
	$\mathcal K_\ell$.
	Since $e_K$ and all left and right creation and annihilation operators
	belong to $\mathcal M(\mathcal K_\ell)$, we have
	$T_{\xi,\eta,\eta'}\in \mathcal M(\mathcal K_\ell)$.
	It is therefore enough to show that $T_{\xi,\eta,\eta'}$ vanishes in the quotient algebra
	$\mathcal M(\mathcal K_\ell)/\mathcal K_\ell$.
	The commutators between left and right creation and annihilation
	operators are either zero or scalar multiples of
	\[
	D_q:=\sum_{n\geq0}q^nP_n\in\mathcal K_\ell.
	\]
	In particular, left and right monomials commute modulo
	$\mathcal K_\ell$.
	
	Decompose the legs of $\xi$ into their $K$- and $K^\perp$-components.
	The $U_t$-invariance of $K_\R$ ensures that the one-particle
	involution in the Wick formula preserves this decomposition.
	Using $\xi\perp\mathcal F_q(K)$ and the Wick formula
	\cite[Lemma~4(1)]{SW16}, we can write $T_{\xi,\eta,\eta'}$
	as a finite linear combination of terms
	\[
	e_KRLR'e_K,
	\]
	where $R,R'$ are right monomials and $L$ is a normally ordered
	left monomial whose vectors lie in $K$ or $K^\perp$,
	with at least one in $K^\perp$.
	If a creation vector of $L$ lies in $K^\perp$, then $e_KL=0$;
	otherwise, an annihilation vector lies in $K^\perp$ and $Le_K=0$.
	
	In the first case, commuting $R$ past $L$ gives
	\[
	e_KRLR'e_K\equiv e_KLRR'e_K=0
	\pmod{\mathcal K_\ell}.
	\]
	In the second case, commuting $R'$ past $L$ gives
	\[
	e_KRLR'e_K\equiv e_KRR'Le_K=0
	\pmod{\mathcal K_\ell}.
	\]
	Thus every term belongs to $\mathcal K_\ell$, proving the assertion.
\end{proof}

\begin{proof}[{\bf Proof of Corollary \ref{corG}}]
	Put $M:=\Gamma_q(H_\R,U_t)''$.
	Assume that $V\neq1$ and that $\sigma^V=\Ad(u)$
	for some $u\in\mathcal U(M)$.
	As in \cite[Lemma~3.2]{HT19}, define an orthogonal representation
	\[
	\rho:\mathbb R\times\mathbb Z\to\mathcal O(H_{\mathbb R}),
	\qquad
	\rho_{(t,n)}=U_tV^n=V^n U_t.
	\]
	
	Assume first that $\rho$ is reducible. Write
	$H_{\mathbb R}=K_{\mathbb R}\oplus K_{\mathbb R}^{\perp}$
	with both summands nonzero and $\rho$-invariant, and put
	\[
	N_1=\Gamma_q(K_{\mathbb R},U|_{K_{\mathbb R}})'',
	\qquad
	N_2=\Gamma_q(K_{\mathbb R}^{\perp},U|_{K_{\mathbb R}^{\perp}})''.
	\]
	Since $V$ leaves both summands invariant, we have
	$uN_i u^*=N_i$ for $i=1,2$. Applying
	Theorem~\ref{thmF} to each of the two decompositions gives
	$u\in N_1\cap N_2=\mathbb C1$. Hence $\sigma^V=\id$, and therefore
	$V=1$.
	
	Assume next that $\rho$ is irreducible. Since $\mathbb R\times\mathbb Z$
	is abelian, \cite[Proposition~3.1]{HT19} implies that
	$\dim H_{\mathbb R}\leq 2$. Hence $\dim H_{\mathbb R}=2$.
	As in \cite[Lemma~3.2]{HT19}, there exist $0<\lambda\leq1$ and
	$\mu>0$ and an orthonormal basis of $H_{\mathbb R}$ such that
	\[
	\rho_{(t,n)}
	=
	\begin{pmatrix}
		\cos(t\log\lambda+n\log\mu)&-\sin(t\log\lambda+n\log\mu)\\
		\sin(t\log\lambda+n\log\mu)&\cos(t\log\lambda+n\log\mu)
	\end{pmatrix},
	\qquad
	(t,n)\in \R\times \Z.
	\]
	
	If $\lambda=1$, then $(U_t)_t$ is trivial and
	$M=\Gamma_q(H_{\mathbb R})''$. For every nonzero
	$\xi\in H_{\mathbb R}$, irreducibility of $\rho$
	implies $\mathbb R\xi\neq\mathbb R V\xi$. Then observe that
	\[
	\sigma^V\bigl(\{W(\xi)\}''\bigr)
	=\{W(V\xi)\}'' \qquad\text{in }M.
	\]
	By \cite[Proposition~4.12]{CSW17}, these two generator masas are not
	unitarily conjugate, contradicting $\sigma^V=\Ad(u)$.

	Suppose finally that $0<\lambda<1$. Then by \cite[Theorem~2.4]{KSW23}, $M$ is a
	type $\mathrm{III}_\lambda$ factor. Putting
	$s=\frac{\log\mu}{\log\lambda}$, we have $V=U_s$ and $\sigma^V=\sigma^\varphi_{-s}$, where $\varphi$
	is the vacuum state. Since $\sigma^V$ is inner and $M$ is of type III$_\lambda$, this means that $s\in \frac{2\pi}{\log\lambda}\mathbb Z$.	
	Hence $V=U_s=1$, a contradiction.
	This proves that $\sigma^V$ is outer whenever $V\neq1$.
	
	For the last assertion, let $g\in G\setminus\{e\}$. Since $\pi$ is
	faithful, $\pi_g\neq1$, and hence $\sigma^{\pi_g}$ is outer by the
	first part. Thus $\sigma^\pi$ is an outer action on the factor $M$,
	and consequently $M\rtimes_{\sigma^\pi}G$ is a factor.
\end{proof}

\section{Strong solidity}
\subsection{Proof of Theorem \ref{thmC}}

We first prove solidity for $q$-Araki--Woods factors.

\begin{thm}\label{thm-solid}
	Let $(H_\R,U_t)$ be a strongly continuous orthogonal representation of $\R$.
	Then $p\Gamma_q(H_\R,U_t)'' p$ is solid for every projection $p\in \Gamma_q(H_\R,U_t)''$. This measn that for every diffuse von Neumann subalgebra $A\subset p\Gamma_q(H_\R,U_t)''p$ with expectation, the relative commutant $A'\cap p\Gamma_q(H_\R,U_t)''p$ is amenable.
\end{thm}

\begin{proof}
	We follow the argument in \cite[Theorem~7.1]{HI15b},
	using the length AO map.
	We only indicate the main steps.
	Put $M:=\Gamma_q(H_\R,U_t)''$. 
	For $J:=J_q$ and $A_q\subset M$, Theorem~\ref{thm-length-AO} gives a unital
	$*$-homomorphism
	\[
	\nu\colon A_q\otimes_{\min}JA_qJ
	\longrightarrow
	\mathcal M(\mathcal K_\ell)/\mathcal K_\ell,
	\qquad
	\nu(a\otimes JbJ)=aJbJ+\mathcal K_\ell.
	\]
	Let $A\subset pMp$ be a diffuse von Neumann subalgebra with expectation.
	We show that $P:=A'\cap pMp$ is amenable.
	By passing to a subalgebra with expectation,
	we may assume that $A$ is a diffuse abelian von Neumann subalgebra.
	Suppose that $P$ is nonamenable.
	Cutting down by a nonzero central projection of $P$,
	we may assume that $P$ has no amenable direct summand.
	
	Let $\alpha$ be the canonical malleable deformation for $M$ as in Subsection~\ref{Malleable deformations}.
	Since $A$ is diffuse,
	Theorem~\ref{thmB}
	shows that $\alpha_t$ does not converge uniformly on $(A)_1$.
	Moreover, $\mathcal K_\alpha=\mathcal K_\ell$ by
	Lemma~\ref{lem-Kalpha-Kell}.
	By Proposition~\ref{prop-solid}, there exist a nonzero projection
	$r\in P'\cap pMp$ and a ucp map
	\[
	\Psi\colon \B(L^2(M))\longrightarrow \B(rL^2(M))
	\]
	whose restriction to $M$ is normal and such that
	\[
	\Psi(x)=xr\quad(x\in P),\qquad
	\Psi(x')=x'r\quad(x'\in M'),\qquad
	\Psi(\mathcal K_\ell)=0.
	\]
	Since $M'$ lies in the multiplicative domain of $\Psi$,
	we have $\Psi(M)\subset rMr$.
	Moreover, the restriction of $\Psi$ to
	$\mathcal M(\mathcal K_\ell)$ factors through the quotient.
	Composing the induced map with $\nu$, we obtain a ucp map
	\[
	\Phi\colon A_q\otimes_{\min}JA_qJ
	\longrightarrow \B(rL^2(M)),
	\qquad
	\Phi(a\otimes JbJ)=\Psi(aJbJ)=\Psi(a)JbJr.
	\]
	Since $A_q$ is exact \cite{KN10} and hence locally reflexive,
	\cite[Lemma~3.6]{HI15b} shows that the normal ucp map
	$\Psi|_M\colon M\to rMr$ admits a ucp extension
	\[
	E\colon \B(L^2(M))\longrightarrow rMr.
	\]
	In particular, $E(x)=xr$ for all $x\in P$.
	As in the final part of the proof of \cite[Theorem~4.1]{HI15b},
	it follows that $P$ has a nonzero amenable direct summand, a contradiction.
\end{proof}

We next establish the key lemma for the proof of Theorem~\ref{thmC}.
Let $G$ be a locally compact amenable group and let
$\pi\colon G\to\mathcal O(H_\R)$ be a strongly continuous orthogonal representation
commuting with $(U_t)_t$.
Put
\[
M:=\Gamma_q(H_\R,U_t)''\rtimes_{\sigma^\pi}G
\subset
\Gamma_q(\widetilde H_\R,U_t)''\rtimes_{\sigma^\pi}G
=:\widetilde M,
\]
and denote by $E_M\colon\widetilde M\to M$ the canonical conditional expectation.
We consider the s-malleable deformation $(\alpha,\beta)$ of the inclusion
$M\subset\widetilde M$ as in Subsection \ref{Malleable deformations}.

Assume that $\widetilde M$ admits a semifinite trace $\Tr$ such that
\[
\Tr\circ E_M=\Tr,
\qquad
\Tr\circ\alpha_t=\Tr\quad(t\in\R),
\qquad
\Tr\circ\beta=\Tr.
\]

\begin{lem}\label{lem-stably-strongly-solid}
	Let $p\in M$ be a projection with $\Tr(p)<\infty$, and let
	$A\subset pMp$ be an amenable von Neumann subalgebra.
	If $\alpha_t$ does not converge uniformly on $(A)_1$ as $t\to0$,
	then $s\mathcal N_{pMp}(A)''$ has a nonzero amenable direct summand.
\end{lem}
\begin{proof}
	We adapt the argument of \cite[Subsection~3.4]{PV12}, as used in
	\cite[Theorem~4.2.1]{Is12}, to stable normalizers, using the stable AO property.
	Since $\Gamma_q(H_\R,U_t)''$ has the CMAP \cite{Av11,ABW17}, see also \cite{Wa19}, and $G$ is amenable,
	$M$ has the CMAP.
	We may therefore take a state
	\[
	\Omega\in\B(L^2(M)\otimes L^2(M^{\op}))^*
	\]
	as in Proposition~\ref{prop-deformation-stable-normalizer}, satisfying conditions (i)--(iii).
	Put
	\[
	A_q(G):=A_q\rtimes_{\sigma^\pi,\mathrm{red}}G,
	\qquad
	B:=\B(L^2(M^{\op})).
	\]
	Let $\mathcal K_{\ell,B}(G)$ be the hereditary $C^*$-algebra
	generated by the projections
	$P_{\leq N}\otimes1_{L^2(G)}\otimes1_B$,
	$N\in\mathbb N_0$, as in Corollary~\ref{cor-stable-AO}.

	Since $P_{\leq N}\otimes1_{L^2(G)}\in \mathcal K_\alpha$ by Lemma \ref{lem-Kalpha-Kell} and
	$\Omega(\mathcal K_\alpha\otimes1_B)=0$ by (i),
	the Cauchy--Schwarz inequality gives $\Omega(\mathcal K_{\ell,B}(G))=0$.
	Thus $\Omega$ induces a state on
	$\mathcal M(\mathcal K_{\ell,B}(G))/\mathcal K_{\ell,B}(G)$.
	Composing it with the stable AO map in Corollary~\ref{cor-stable-AO}, we obtain a state
	\[
	\Phi\colon A_q(G) \otimes_{\min}JA_q(G)J\otimes_{\min}B\to\C,
	\qquad
	\Phi(a\otimes JbJ\otimes c)=\Omega(aJbJ\otimes c).
	\]
	If $A_q(G)$ is non-unital, $\Phi$ extends to a state
	after unitizing the first two tensor factors.
	Since $A_q$ is exact and $G$ is amenable, $A_q(G)$ is exact,
	hence locally reflexive.
	Moreover, the restrictions of $\Omega$ to $M\otimes\C$
	and $JMJ\otimes\C$ are normal by (ii).
	Applying the local reflexivity argument in the proof of
	\cite[Lemma~3.6]{HI15b} successively to the first two tensor factors,
	we obtain a state
	\[
	\Phi\colon M\otimes_{\min}JMJ\otimes_{\min}B\to\C
	\]
	satisfying the same formula for all $a,b\in M$ and $c\in B$.
	Extend $\Phi$ further, by the Hahn--Banach theorem, to a state on
	$\B(L^2(M))\otimes_{\min}JMJ\otimes_{\min}B$ and define
	\[
	E(X):=\Phi(X\otimes1\otimes1),
	\qquad X\in\B(L^2(M)).
	\]
	By (i) and (ii), $E|_M$ is normal and $E(p)=1$.
	
	Fix a partial isometry $v\in s\mathcal N^0_{pMp}(A)$ and take $L_v$ and $R_{v,k}$ for $k\geq1$ as in Proposition~\ref{prop-deformation-stable-normalizer}(iii).
	Then each $R_{v,k}$ 
	belong to $JMJ\otimes_{\alg}B$.
	Put
	\[
	D(v,k):=v\otimes1\otimes1-1\otimes R_{v,k}.
	\]
	By (iii) and the formula for $\Phi$,
	\[
	\begin{aligned}
		\Phi(D(v,k)^*D(v,k))
		&=\Omega((L_v-R_{v,k})^*(L_v-R_{v,k}))\to0,\\
		\Phi(D(v,k)D(v,k)^*)
		&=\Omega((L_v-R_{v,k})(L_v-R_{v,k})^*)\to0,
	\end{aligned}
	\]
	as $k\to\infty$. 
	For $X\in\B(L^2(M))$, the operator $1\otimes R_{v,k}$ commutes
	with $X\otimes1\otimes1$.
	Consequently, the Cauchy--Schwarz inequality gives
	\[
	\begin{aligned}
		|E(vX)-E(Xv)|
		&=\bigl|\Phi([D(v,k),X\otimes1\otimes1])\bigr|\\
		&\leq \|X\|
		\left(
		\Phi(D(v,k)D(v,k)^*)^{1/2}
		+\Phi(D(v,k)^*D(v,k))^{1/2}
		\right)
		\longrightarrow0,
	\end{aligned}
	\]
	as $k\to\infty$.
	Thus $E(vX)=E(Xv)$.
	
	These partial isometries generate $P:=s\mathcal N_{pMp}(A)''$.
	Since $E|_M$ is normal, the same approximation argument as in the last part of the proof of
	\cite[Proposition~3.7]{BHV15} shows that $E$ is $P$-central.	
	Its restriction to $p\B(L^2(M))p$ is therefore a $P$-central state
	whose restriction to $pMp$ is normal.
	It follows that $P$ has a nonzero amenable direct summand.
\end{proof}

\begin{proof}[{\bf Proof of Theorem \ref{thmC}}]
	Put $M:=\Gamma_q(H_\R,U_t)''$. The tracial case follows from
	Theorem~\ref{thmB}
	and Lemma~\ref{lem-stably-strongly-solid} with $G=\{e\}$,
	after cutting down by a central projection of the stable normalizer
	if necessary.
	We may therefore assume that $M$ is a type III factor.

	We first follow the reduction in \cite[Section~4]{BHV15}.
	Since no new argument is involved at this stage, we only indicate the main steps
	and use the standard properties of continuous cores without further comment.
	Since every nonzero corner of $M$ is isomorphic to $M$, and strong
	solidity coincides with stable strong solidity in this case
	\cite[Corollary~5.2]{BHV15}, it suffices to prove that $M$ is strongly solid.
	
	Suppose otherwise.
	By Theorem~\ref{thm-solid} and the reductions in the proof of the
	main theorem of \cite{BHV15}, we may choose a diffuse amenable
	von Neumann subalgebra $A\subset M$ with expectation $E_A$ such that
	\[
	A'\cap M=\mathcal Z(A),
	\qquad
	P:=\mathcal N_M(A)''
	\text{ has no amenable direct summand}.
	\]
	Choose a diffuse abelian von Neumann subalgebra $A_0\subset A$
	with expectation, and a faithful normal state $\psi$ on $M$ such that
	$\psi\circ E_A=\psi$ and $A_0\subset M_\psi$.
	
	Let $M\subset^{E_M}\widetilde M$ be the canonical deformation inclusion,
	and let $\varphi$ be the vacuum state.
	We use the same notation for the extensions of $\psi$ and $\varphi$
	to $\widetilde M$ through $E_M$.
	Put
	\[
	\mathcal M:=C_\varphi(M)
	\subset C_\varphi(\widetilde M)=:\widetilde{\mathcal M}.
	\]
	Using the canonical trace-preserving isomorphisms between the cores
	associated with $\psi$ and $\varphi$, we view
	$\mathcal A:=C_\psi(A)$ as a subalgebra of $\mathcal M$,
	with $L_\psi\R\subset\mathcal A$.
	By \cite[Lemma~4.1]{BHV15} and the argument in the proof of its
	main theorem, we have an inclusion with expectation
	\[
	C_\psi(P)\subset
	\mathcal P:=\mathcal N_{\mathcal M}(\mathcal A)'',
	\]
	and $\mathcal P$ has no amenable direct summand.
	
	The malleable deformation $(\alpha,\beta)$ extends to
	$\widetilde{\mathcal M}$ and preserves its canonical trace $\Tr$.
	The canonical expectation onto $\mathcal M$ also preserves $\Tr$.
	Thus we are in the setting of
	Lemma~\ref{lem-stably-strongly-solid}, with $G=\R$ and
	the modular action of $\varphi$.

	Fix a nonzero projection $r\in L_\psi\R$ with $\Tr(r)<\infty$.
	By \cite[Lemma~3.4]{BHV15},
	\[
	r\mathcal P r
	\subset
	s\mathcal N_{r\mathcal M r}(r\mathcal A r)''.
	\]
	Consequently, the latter algebra has no amenable direct summand.
	We will obtain a contradiction by showing that $\alpha_t$ does not
	converge uniformly on $(r\mathcal A r)_1$.
	
	The formula for the canonical trace gives
	\[
	\|xr\|_{\Tr}^2=\Tr(r)\psi(x^*x),
	\qquad x\in\widetilde M.
	\]
	Since $A_0\subset M_\psi$, the projection $r$ commutes with $A_0$,
	and hence $A_0r\subset r\mathcal A r$.
	For every $u\in\mathcal U(A_0)$, we have
	\[
	\begin{aligned}
		\|\alpha_t(ur)-ur\|_{\Tr}
		&\geq
		\|(\alpha_t(u)-u)r\|_{\Tr}
		-\|\alpha_t(r)-r\|_{\Tr}\\
		&=
		\Tr(r)^{1/2}\|\alpha_t(u)-u\|_{\psi}
		-\|\alpha_t(r)-r\|_{\Tr}.
	\end{aligned}
	\]
	By Theorem~\ref{thmB},
	$\alpha$ does not converge uniformly on $\mathcal U(A_0)$
	in $\|\cdot\|_{\psi}$.
	Since $\|\alpha_t(r)-r\|_{\Tr}\to0$ as $t\to 0$, the preceding estimate shows
	that $\alpha_t$ does not converge uniformly on $(r\mathcal A r)_1$
	in $\|\cdot\|_{\Tr}$.
	
	As $r\mathcal A r$ is amenable,
	Lemma~\ref{lem-stably-strongly-solid} now implies that
	$s\mathcal N_{r\mathcal M r}(r\mathcal A r)''$
	has a nonzero amenable direct summand, a contradiction.
\end{proof}

\section{Connes' $\tau$-invariant and non-isomorphism}

In this section, we compute Connes' $\tau$-invariant
and prove the non-isomorphism result of Theorem~\ref{thmE}.

\subsection{Proof of Proposition~\ref{propD}}

A uniformly bounded net $(x_j)_j$ in a von Neumann algebra $N$
is called \emph{centralizing} if 
$\|x_j\psi-\psi x_j\|_{N_*}\to0$ as $j\to \infty$ for every $\psi\in N_*$.
Following \cite{Co74,Ma16b}, a factor $N$ is called \emph{full}
if every centralizing net $(x_j)_j$ admits a net of scalars
$(\lambda_j)_j$ such that $x_j-\lambda_j1\to0$ strongly as $j\to \infty$.

We equip $\mathrm{Aut}(N)$ with the $u$-topology and
$\mathrm{Out}(N)=\mathrm{Aut}(N)/\mathrm{Int}(N)$ with the
quotient topology, where $\mathrm{Int}(N)$ denotes the group
of inner automorphisms.
If $N$ is full, then $\mathrm{Int}(N)$ is closed in
$\mathrm{Aut}(N)$, so $\mathrm{Out}(N)$ is Hausdorff
\cite[Theorem~3.1]{Co74}\cite[Section~2.4]{Ma16b}.

We will use the following elementary observation.

\begin{lem}\label{lem-compact-tau}
	Let $M$ be a factor such that $\mathrm{Int}(M)$ is closed
	in $\mathrm{Aut}(M)$, and let
	$K\subset\mathrm{Aut}(M)$ be a compact subgroup.
	Put
	\[
	N=M^K
	=\{x\in M\mid \alpha(x)=x\text{ for every }\alpha\in K\}.
	\]
	If $N'\cap M=\mathbb C$, then the restriction of the quotient map
	$\pi|_K\colon K\longrightarrow\mathrm{Out}(M)$
	is a homeomorphism onto its image.
	In particular, for every net $(\alpha_j)_j$ in $K$,
	\[
	[\alpha_j]\longrightarrow[\mathrm{id}]
	\quad\text{in }\mathrm{Out}(M)
	\quad\Longleftrightarrow\quad
	\alpha_j\longrightarrow\mathrm{id}
	\quad\text{in }\mathrm{Aut}(M).
	\]
\end{lem}

\begin{proof}
	Let $\alpha\in K\cap\mathrm{Int}(M)$ and write
	$\alpha=\operatorname{Ad}(u)$ for some $u\in\mathcal U(M)$.
	Since $\alpha$ fixes $N$ pointwise, we have
	$u\in N'\cap M=\mathbb C$, and hence $\alpha=\mathrm{id}$.
	Thus $\pi|_K$ is injective.
	By assumption, $\mathrm{Out}(M)$ is Hausdorff.
	Consequently, $\pi|_K$ is a continuous injection from a
	compact space into a Hausdorff space, and hence it is a
	homeomorphism onto its image.
	The last assertion follows immediately.
\end{proof}

We next establish a spectral gap estimate.
Let $(H_\R,U_t)$ be a strongly continuous orthogonal representation.
Put $H=H_\R\otimes_\R\C$, and write $U_t=A^{it}$ for $t\in \R$.
Set
\[
M:=\Gamma_q(H_\R,U_t)'',
\qquad
\mathcal F_q:=\mathcal F_q(H),
\qquad
\mathcal F_q^+:=\mathcal F_q(H)\ominus\C\Omega.
\]
Following the notation of \cite{HI20}, we define the concatenation maps
\[
\begin{aligned}
	\Phi\colon H\otimes\mathcal F_q
	&\longrightarrow\mathcal F_q^+,
	& \Phi(\xi\otimes\eta)&=\ell_q(\xi)\eta,\\
	\Psi\colon\mathcal F_q\otimes H
	&\longrightarrow\mathcal F_q^+,
	& \Psi(\eta\otimes\xi)&=r_q(\xi)\eta.
\end{aligned}
\]
Here the tensor products in the domains are ordinary Hilbert
space tensor products.
These maps are bounded isomorphisms; see, e.g.,
\cite[Proposition~6.1]{HI20}.
Fix constants $C_1,C_2>0$, depending only on $q$, such that
\[
\|\Phi\|=\|\Psi\|\leq C_1,
\qquad
\|\Phi^{-1}\|=\|\Psi^{-1}\|\leq C_2.
\]

The following spectral gap estimate is a refinement of the
arguments in \cite{Sn03,HI20}.
For $R\geq1$, we denote by $E_R$ the spectral projection of $A$
corresponding to $[R^{-1},R]$.

\begin{lem}\label{lem-spectral-gap}
	For each $C\geq1$ and $d\in\N$ satisfying
	\[
	d-C_1C_2>2CC_1C_2\sqrt d,
	\]
	there exists a constant $\kappa>0$, depending only on
	$q,d,C$, with the following property.
	
	For any two orthonormal families
	$(e_i)_{i=1}^d$ and $(f_i)_{i=1}^d$ in
	$E_{C^2}H_\R$, we have
	\[
	\|\eta-\langle\Omega,\eta\rangle\Omega\|^2
	\leq
	\kappa\sum_{i=1}^d
	\|(W(e_i)-W_r(f_i))\eta\|^2
	\qquad(\eta\in\mathcal F_q).
	\]
\end{lem}
\begin{proof}
	Since $IAI=A^{-1}$, the spectral projection $E_{C^2}$
	commutes with $I$. In particular, $Ie_i=e_i$ and $If_i=f_i$.
	Denoting by $P_e$ and $P_f$ the orthogonal projections
	onto the complex linear spans of $(e_i)_{i=1}^d$ and
	$(f_i)_{i=1}^d$, respectively, we have
	\[
	\|A^{1/2}P_e\|\leq C,
	\qquad
	\|A^{-1/2}P_f\|\leq C.
	\]
	
	We first follow the argument of \cite{HI20}, recalling only
	the estimates that we need below.
	Define $\mathcal M\colon\mathcal F_q\to H\otimes\mathcal F_q$
	by
	\[
	\mathcal M\eta
	=
	\sum_{i=1}^d
	e_i\otimes(W(e_i)-W_r(f_i))\eta\qquad(\eta\in \mathcal F_q)
	\]
	and note that
	$\|\mathcal M\eta\|^2=
	\sum_{i=1}^d
	\|(W(e_i)-W_r(f_i))\eta\|^2$.
	The Wick formulas yield
	$\mathcal M=\mathcal M_{\mathrm{cr}}+\mathcal M_{\mathrm{an}}$,
	where
	\[
	\begin{aligned}
		\mathcal M_{\mathrm{cr}}\eta
		&=
		\sum_{i=1}^d
		e_i\otimes(\ell_q(e_i)-r_q(f_i))\eta,\\
		\mathcal M_{\mathrm{an}}\eta
		&=
		\sum_{i=1}^d
		e_i\otimes
		\bigl(\ell_q(A^{1/2}e_i)^*
		-r_q(A^{-1/2}f_i)^*\bigr)\eta.
	\end{aligned}
	\]
	
	Let $V\colon P_eH\to P_fH$ be the unitary satisfying
	$Ve_i=f_i$, and let
	$\Sigma\colon H\otimes\mathcal F_q\to\mathcal F_q\otimes H$
	be the tensor flip.
	As in the proof of \cite[Lemma~6.4]{HI20}, we have
	\[
	\mathcal M_{\mathrm{an}}^*
	=
	\Phi(A^{1/2}P_e\otimes1_{\mathcal F_q})
	-
	\Psi(1_{\mathcal F_q}\otimes A^{-1/2}VP_e)\Sigma.
	\]
	Consequently,
	\[
	\|\mathcal M_{\mathrm{an}}\|
	\leq
	C_1\bigl(\|A^{1/2}P_e\|
	+\|A^{-1/2}VP_e\|\bigr)
	\leq 2CC_1.
	\]
	
	For the creation part, we adapt the proof of
	\cite[Lemma~6.5]{HI20}.
	Put
	$v=\sum_{i=1}^d e_i\otimes e_i\in H\otimes H$
	and define
	\[
	F=(\langle v,\cdot\rangle\otimes1_{\mathcal F_q})
	(1_H\otimes\Phi^{-1})
	\colon H\otimes\mathcal F_q^+\to\mathcal F_q.
	\]
	Since $\|v\|=\sqrt d$, we have $\|F\|\leq C_2\sqrt d$.
	Also define
	\[
	S_V
	=
	\Psi(1_{\mathcal F_q}\otimes VP_e)\Sigma\Phi^{-1}
	\colon\mathcal F_q^+\to\mathcal F_q^+.
	\]
	Then $\|S_V\|\leq C_1C_2$, and on elementary tensors,
	\[
	S_V(\xi_1\otimes\cdots\otimes\xi_n)
	=
	\xi_2\otimes\cdots\otimes\xi_n\otimes VP_e\xi_1
	\qquad(n\geq1).
	\]
	A direct computation on elementary tensors gives
	\[
	F\left(\sum_{i=1}^d e_i\otimes\ell_q(e_i)\eta\right)
	=d\eta,
	\qquad
	F\left(\sum_{i=1}^d e_i\otimes r_q(f_i)\eta\right)
	=S_V\eta
	\qquad(\eta\in\mathcal F_q^+).
	\]
	Thus
	$F\mathcal M_{\mathrm{cr}}|_{\mathcal F_q^+}
	=
	d\,1_{\mathcal F_q^+}-S_V$,
	and hence
	\[
	\|\mathcal M_{\mathrm{cr}}\eta\|
	\geq L_d\|\eta\|
	\qquad(\eta\in\mathcal F_q^+),
	\qquad
	L_d:=\frac{d-C_1C_2}{C_2\sqrt d}>0.
	\]
	
	So far, we have adapted the estimates of \cite{HI20}
	to two orthonormal families.
	We now use an additional observation to deal with the
	vacuum component.
	
	Let $\eta\in\mathcal F_q$ and put
	$\eta^+=\eta-\langle\Omega,\eta\rangle\Omega$.
	Denote by $P_{>1}$ the orthogonal projection onto
	$\mathcal F_q\ominus(\mathbb C\Omega\oplus H)$, and put
	$Q=1_H\otimes P_{>1}$.
	Since $\mathcal M_{\mathrm{cr}}$ raises the tensor degree
	by one and $\mathcal M_{\mathrm{an}}$ vanishes on the vacuum,
	we have
	\[
	Q\mathcal M\Omega=0,
	\qquad
	Q\mathcal M_{\mathrm{cr}}\eta^+
	=
	\mathcal M_{\mathrm{cr}}\eta^+.
	\]
	Therefore,
	$Q\mathcal M\eta
	=
	\mathcal M_{\mathrm{cr}}\eta^+
	+
	Q\mathcal M_{\mathrm{an}}\eta^+$, and 
	it follows that
	\[
	\begin{aligned}
		\|\mathcal M\eta\|
		&\geq \|Q\mathcal M\eta\|\\
		&\geq
		\|\mathcal M_{\mathrm{cr}}\eta^+\|
		-
		\|Q\mathcal M_{\mathrm{an}}\eta^+\|\\
		&\geq
		(L_d-2CC_1)\|\eta^+\|.
	\end{aligned}
	\]
	Hence, taking
	\[
	\kappa
	=
	(L_d-2CC_1)^{-2}
	=
	\left(
	\frac{d-C_1C_2}{C_2\sqrt d}-2CC_1
	\right)^{-2}>0,
	\]
	we obtain the desired inequality.
\end{proof}

\begin{proof}[{\bf Proof of Proposition~\ref{propD}}]
	The case $q=0$ follows from \cite[Theorem~2.7]{Va05}.
	Although separability is assumed there, the spectral gap
	argument in \cite[Lemma~4.1 and Proposition~4.2]{Va05} uses
	only finitely many test elements at each approximation step
	and works equally well with arbitrary nets.
	We therefore assume that $q\neq0$.
	
	Denote the vacuum state by $\varphi$.
	Since the canonical implementation of $\sigma_t^\varphi$ on
	$\mathcal F_q$ is $\mathcal F_q(U_{-t})$, we have
	\[
	U_{t_j}\longrightarrow1
	\quad\text{strongly}
	\quad\Longleftrightarrow\quad
	\sigma_{t_j}^\varphi\longrightarrow\mathrm{id}
	\quad\text{in }\mathrm{Aut}(M)
	\]
	for every net $(t_j)_j$ in $\mathbb R$.
	
	Suppose first that $H_{\mathbb R}$ is separable and $U$ is
	almost periodic.
	Then $M$ is full by \cite{KSW23,KW24}, and
	$M_\varphi'\cap M=\mathbb C$
	by \cite[proof of Theorem~2.4]{KSW23}.
	Since $\varphi$ is almost periodic,
	\[
	K=\overline{\{\sigma_t^\varphi\mid t\in\mathbb R\}}^{\,u}
	\subset\mathrm{Aut}(M)
	\]
	is a compact subgroup and $M^K=M_\varphi$.
	Lemma~\ref{lem-compact-tau} therefore shows that
	\[
	[\sigma_{t_j}^\varphi]\longrightarrow[\mathrm{id}]
	\quad\Longleftrightarrow\quad
	\sigma_{t_j}^\varphi\longrightarrow\mathrm{id}
	\]
	for every net $(t_j)_j$ in $\mathbb R$.
	This proves the proposition in this case.
	
	We next consider the remaining cases.
	Then either $H_{\mathbb R}$ is nonseparable or $U$ has a
	nonzero weakly mixing part.
	In either case, we may choose $R>1$ such that
	\[
	E_RH\text{ is infinite dimensional},
	\qquad
	E_R=1_{[R^{-1},R]}(A).
	\]
	By \cite[Lemma~6.3]{HI20}, we may choose an orthonormal
	family $e_1,\ldots,e_d$ in $E_RH_\R$,
	where $d$ is large enough to apply Lemma~\ref{lem-spectral-gap}
	with $C=R^{1/2}$.
	As $U_t$ commutes with $A$ and $I$, the families
	$(e_k)_{k=1}^d$ and $(U_{-t}e_k)_{k=1}^d$ satisfy the
	hypotheses of that lemma uniformly in $t$.
	Hence there is a constant $\kappa>0$ such that
	\[	
	\|\eta-\langle\Omega,\eta\rangle\Omega\|^2
	\leq
	\kappa\sum_{k=1}^d
	\|(W(e_k)-W_r(U_{-t}e_k))\eta\|^2
	\]
	for all $t\in\mathbb R$ and $\eta\in\mathcal F_q$.
	Using $W_r(U_{-t}e_k)\in M'$ and
	$W_r(U_{-t}e_k)\Omega
	=
	\sigma_t^\varphi(W(e_k))\Omega$,
	taking $\eta=x\Omega$ gives
	\[
	\|x-\varphi(x)1\|_\varphi^2
	\leq
	\kappa\sum_{k=1}^d
	\|W(e_k)x
	-x\sigma_t^{\varphi}(W(e_k))\|_\varphi^2
	\qquad(x\in M).
	\]
	Since every centralizing net asymptotically commutes
	$*$-strongly with each fixed element of $M$, applying the
	inequality to such a net and its adjoint shows that $M$ is full.
	
	Assume that $[\sigma_{t_j}^\varphi]\to[\mathrm{id}]$ for a net
	$(t_j)_j$ in $\mathbb R$.
	To prove that $\sigma_{t_j}^\varphi\to\mathrm{id}$, it suffices
	to show that every subnet has a further subnet converging to
	$\mathrm{id}$.
	Starting with an arbitrary subnet, the openness of the quotient
	map $\mathrm{Aut}(M)\to\mathrm{Out}(M)$ allows us to pass to a
	further subnet and relabel so that there exist unitaries
	$u_j\in\mathcal U(M)$ satisfying
	\[
	\alpha_j:=
	\operatorname{Ad}(u_j)\circ\sigma_{t_j}^\varphi
	\longrightarrow\mathrm{id}
	\quad\text{in }\mathrm{Aut}(M),\quad \text{as }j\to \infty.
	\]
	The $\sigma^\varphi$-invariance of $\varphi$ gives
	$\|\varphi\circ\operatorname{Ad}(u_j)-\varphi\|\to0$.
	Together with $\alpha_j(x)\to x$ $*$-strongly, this implies that, for each $x\in M$,
	\[
	\|xu_j-u_j\sigma_{t_j}^\varphi(x)\|_\varphi
	\longrightarrow0,\qquad \text{as }j\to \infty.
	\]
	Applying the spectral gap inequality above, we obtain
	\[
	\|u_j-\varphi(u_j)1\|_\varphi^2
	\leq
	\kappa\sum_{k=1}^d
	\|W(e_k)u_j
	-u_j\sigma_{t_j}^\varphi(W(e_k))\|_\varphi^2
	\longrightarrow0,\qquad \text{as }j\to\infty.
	\]
	Since the $u_j$ are unitaries, this implies 
	$\operatorname{Ad}(u_j)\to\mathrm{id}_M$, and therefore
	$\sigma_{t_j}^\varphi\to\mathrm{id}_M$.
	
	Combining this with the first equivalence proves the formula
	for $\tau(M)$.
\end{proof}

\subsection{Proof of Theorem~\ref{thmE}}

The following result is the $q$-Araki--Woods analogue of the free Araki--Woods case of \cite[Theorem~5.1]{HSV16}; when $q=0$, it coincides with that result.

\begin{lem}\label{lem-state-intertwining}
	Let $M=\Gamma_q(H_\R,U_t)''$ and let $\varphi$ be the vacuum state.
	Let $\psi$ be a faithful normal state on $M$, and let
	$z\in\mathcal Z(M_\psi)$ be the nonzero central projection such that
	$M_\psi(1-z)$ is amenable and $M_\psi z$ has no amenable direct
	summand.
	Then there exists a sequence of nonzero partial isometries $v_n\in M$ such
	that
	\[
	q_n:=v_nv_n^*\in M_\psi,
	\qquad 
	p_n:=v_n^*v_n\in M_\varphi,
	\qquad
	\sum_n q_n=z
	\]
	and
	\[
	q_n[D\psi:D\varphi]_t
	=
	\lambda_n^{it}v_n\sigma_t^\varphi(v_n^*),
	\qquad
	\psi(q_n\, \cdot \, q_n)
	=
	\lambda_n v_n\varphi v_n^*,
	\]
	where $\lambda_n=\psi(q_n)/\varphi(p_n)$. Here $[D\psi:D\varphi]_t$ denotes the Connes cocycle.
\end{lem}
\begin{proof}
	We follow the proof of \cite[Theorem~5.1]{HSV16} and use the same
	notation as there. It suffices to verify the only step in that proof
	where the specific deformation properties of free Araki--Woods
	factors are used.
	
	We use the canonical identification
	\[
	\Pi_{\varphi,\psi}:C_\psi(M)\to C_\varphi(M).
	\]
	It is enough to prove the following claim.
	
	\begin{claim}
		For any nonzero projection $e\in M_\psi z$ and any nonzero
		finite trace projection $r_0\in L_\psi(\R)$, we have
		\[
		\Pi_{\varphi,\psi}(L_\psi(\R))r
		\preceq_{C_\varphi(M)}
		L_\varphi(\R),
		\]
		where $r=\Pi_{\varphi,\psi}(er_0)$.
	\end{claim}
	
	\begin{proof}[{\bf Proof of Claim}]
		Let $M\subset^{E_M}\widetilde M$ and $C_\varphi(M)\subset C_\varphi(\widetilde{M})$
		be the inclusions associated with the malleable deformation as in
		Subsection~\ref{Malleable deformations}. We extend $\varphi$ and
		$\psi$ to $\widetilde M$ via $E_M$, keeping the same notation.
		The canonical isomorphism $\Pi_{\varphi,\psi}$ then extends to
		\[
		\Pi_{\varphi,\psi}:C_\psi(\widetilde M)
		\to C_\varphi(\widetilde M),
		\]
		whose restriction to $C_\psi(M)$ is the above identification.
		In what follows, we use this identification and omit
		$\Pi_{\varphi,\psi}$ from the notation. In particular, we write $r=er_0$ and regard $L_\psi(\R)r$
		as a von Neumann subalgebra of $rC_\varphi(M)r$.
		
		Since $eM_\psi e$ has no amenable direct summand, it is diffuse.
		By Theorem~\ref{thmF}, we have
		\[
		(eM_\psi e)'\cap e\widetilde Me\subset eMe.
		\]
		Using this inclusion, we obtain
		\begin{align*}
			&[(L_\psi\R e)'\cap eC_\varphi(M)e]'
			\cap eC_\varphi(\widetilde M)e \\
			&\qquad=
			[eM_\psi e\ovt L_\psi\R]'
			\cap eC_\varphi(\widetilde M)e \\
			&\qquad=
			[(eM_\psi e)'\cap e\widetilde M_\psi e]
			\ovt L_\psi\R \\
			&\qquad=
			\mathcal Z(M_\psi)e\ovt L_\psi\R.
		\end{align*}
		Cutting down by $r$, we get
		\[
		[(L_\psi\R r)'\cap rC_\varphi(M)r]'
		\cap rC_\varphi(\widetilde M)r
		=
		[\mathcal Z(M_\psi)\ovt L_\psi\R]r
		\subset rC_\varphi(M)r.
		\]
		Thus condition~(b) of Theorem~\ref{thmB} is satisfied.
		
		Suppose, toward a contradiction, that
		$L_\psi(\R)r\not\preceq_{C_\varphi(M)}L_\varphi(\R)$.
		By Theorem~\ref{thmB}, the malleable deformation does not
		converge uniformly on the unit ball of $L_\psi(\R)r$.
		Since this algebra is amenable,
		Lemma~\ref{lem-stably-strongly-solid} implies that
		$s\mathcal N_{rC_\varphi(M)r}(L_\psi(\R)r)''$
		has a nonzero amenable direct summand.
		On the other hand,
		\[
		(eM_\psi e\ovt L_\psi\R)r
		=
		(eM_\psi e)\ovt L_\psi(\R)r_0
		\]
		is contained in this algebra with expectation and has no amenable
		direct summand. This gives a contradiction and proves the claim.
	\end{proof}
	
	The remainder of the proof is exactly the same as that of
	\cite[Theorem~5.1]{HSV16}. Namely,
	\cite[Theorem~3.1]{HSV16} gives the required partial isometry for
	every nonzero $q\in M_\psi z$, and the usual maximality argument
	gives the sequence $(v_n)_n$ satisfying $\sum_nq_n=z$.
\end{proof}

In the remainder of this subsection, we assume that all real Hilbert spaces are separable, since we will use the measure-theoretic description of orthogonal representations.
For a positive nonsingular self-adjoint operator $A$ on a separable
Hilbert space, we denote by $\mathcal C(A)$ the measure class of the
spectral measure of $\log A$. More precisely,
\[
\mathcal C(A)
=
\{E\subset\R \mid E \text{ is Borel and }
1_E(\log A)=0\}.
\]
For a $\sigma$-finite Borel measure $\mu$ on $\R$, we denote by
$\mathcal C(\mu)$ its measure class.
For a sequence $(\mu_n)_n$ of measures, we denote by
$\bigvee_n\mu_n$ any measure such that
$\mathcal C\left(\bigvee_n\mu_n\right)
=
\bigcap_n\mathcal C(\mu_n)$.

Let $U:\R\to\mathcal O(H_\R)$ be a strongly continuous orthogonal
representation. Write 
$U_t=A^{it}$ for $t\in \R$ on $H=H_\R\otimes_\R \C$.
We call $\mathcal C(A)$ the
spectral measure class of $U$. It is symmetric, and we may choose a
finite symmetric Borel measure $\mu$ on $\R$ such that
$\mathcal C(\mu)=\mathcal C(A)$.
We call
$\mathcal C\left(\bigvee_{n\geq1}\mu^{*n}\right)$
the \emph{joint measure class} of $U$. This does not depend on the choice of
the measure $\mu$ representing the spectral measure class of $U$.

We denote by $H_\R^{\mathrm{ap}}$ the almost periodic part of
$(H_\R,U_t)$. Recall that, for a $q$-Araki--Woods algebra
$M=\Gamma_q(H_\R,U_t)''$ with vacuum state $\varphi$, the centralizer
$M_\varphi$ only depends on the almost periodic part of $U$. In
particular, $M_\varphi$ is abelian if
$\dim H_\R^{\mathrm{ap}}\leq1$, whereas it is a nonamenable
$\mathrm{II}_1$ factor if
$\dim H_\R^{\mathrm{ap}}\geq2$; see e.g., \cite[Theorem~2.2 and the proof of Theorem~2.4]{KSW23}.
For a faithful normal state $\rho$ and a nonzero projection $p$
in its centralizer, we denote by $\rho_p$ the normalized
restriction of $\rho$ to $pMp$.

The following theorem is the $q$-Araki--Woods analogue of
\cite[Theorem~5.3]{HSV16}; when $q=0$, it coincides with that result
except for the final assertion concerning a state-preserving isomorphism.

\begin{thm}\label{thm-joint-measure-class}
	For $i=1,2$, let $M_i=\Gamma_q(H_\R^{(i)},U_t^{(i)})''$
	be $q$-Araki--Woods factors, where each $H_\R^{(i)}$ is separable.
	Let $\varphi_i$ be their vacuum
	states. Assume that $\dim H_\R^{(2),\mathrm{ap}}\geq2$.
	If $M_1$ is $\ast$-isomorphic to $M_2$, then there exist nonzero projections
	$p_i\in (M_i)_{\varphi_i}$ for $i=1,2$
	and a state-preserving isomorphism
	\[
	(p_1M_1p_1,(\varphi_1)_{p_1})
	\cong
	(p_2M_2p_2,(\varphi_2)_{p_2}).
	\]
	In particular, the joint measure classes of $U^{(1)}$ and
	$U^{(2)}$ are equal.
\end{thm}

\begin{proof}
	The proof is the same as the first part of
	\cite[Theorem~5.3]{HSV16}, using
	Lemma~\ref{lem-state-intertwining} in place of
	\cite[Theorem~5.1]{HSV16}. We include the short argument.
	
	Let $\pi:M_1\to M_2$ be an isomorphism and put
	$\psi=\varphi_2\circ\pi$.
	Since $\dim H_\R^{(2),\mathrm{ap}}\geq2$, the centralizer
	$(M_1)_\psi\cong(M_2)_{\varphi_2}$ is a nonamenable
	$\mathrm{II}_1$ factor. By
	Lemma~\ref{lem-state-intertwining}, there exist nonzero projections
	$p_1\in(M_1)_{\varphi_1}$ and $e\in(M_1)_\psi$ such that
	\[
	(p_1M_1p_1,(\varphi_1)_{p_1})
	\simeq
	(eM_1e,\psi_e)
	\]
	state-preservingly. Putting $p_2=\pi(e)$ gives the first assertion.
	
	In particular,
	$p_1(M_1)_{\varphi_1}p_1
	\simeq
	p_2(M_2)_{\varphi_2}p_2$
	is a nonamenable $\mathrm{II}_1$ factor. Hence
	$\dim H_\R^{(1),\mathrm{ap}}\geq2$, and
	$(M_1)_{\varphi_1}$ is also a $\mathrm{II}_1$ factor.
	Therefore, \cite[Lemma~2.1]{HSV16}, applied twice, yields
	\[
	\mathcal C(\Delta_{\varphi_1})
	=
	\mathcal C(\Delta_{(\varphi_1)_{p_1}})
	=
	\mathcal C(\Delta_{(\varphi_2)_{p_2}})
	=
	\mathcal C(\Delta_{\varphi_2}).
	\]
	Finally, since 
	$\Delta_{\varphi_i}^{it}=
	\mathcal F_q(U_{-t}^{(i)})$ for $i=1,2$ and $t\in \R$,
	\[
	\mathcal C(\Delta_{\varphi_i})
	=
	\mathcal C\left(
	\bigvee_{n\geq1}\mu_i^{*n}
	\right),
	\]
	where $\mu_i$ is any finite symmetric Borel measure representing
	the spectral measure class of $U^{(i)}$. Thus the joint measure
	classes of $U^{(1)}$ and $U^{(2)}$ are equal.
\end{proof}

\begin{proof}[{\bf Proof of Theorem~\ref{thmE}}]
	We apply the measure-theoretic construction in
	\cite[Example~5.4]{HSV16}, including the addition of a common
	Lebesgue component.
	This gives a continuum of strongly continuous orthogonal
	representations on separable real Hilbert spaces, each with a
	Lebesgue component and an almost periodic part of dimension
	at least two, whose joint measure classes are pairwise distinct.
	
	Since each representation contains a Lebesgue component, the
	corresponding $q$-Araki--Woods factors are all of type
	$\mathrm{III}_1$, and Proposition~\ref{propD} shows that their
	$\tau$-invariant is the usual topology on $\R$.
	Finally, Theorem~\ref{thm-joint-measure-class} implies that
	these factors are pairwise non-isomorphic.
\end{proof}

\appendix

\section{Almost periodic approximation}
\label{Almost periodic approximation}

We record in this appendix a general approximation procedure relating arbitrary
$q$-Araki--Woods algebras to those associated with almost periodic representations,
which may also be useful in other contexts.

\begin{lem}\label{lem-almost-periodic-approximation-of-generator}
	Let $H_\R$ be a real Hilbert space, put $H=H_\R\otimes_\R\C$, and denote
	by $I$ the canonical conjugation on $H$.
	Let $A$ be a positive injective self-adjoint operator on $H$ satisfying
	$IAI=A^{-1}$.
	There exist positive injective self-adjoint operators $A_n$ on $H$ and
	positive invertible contractions $T_n\in\B(H)$ with the following properties.
	\begin{enumerate}
		\item
		Each $A_n$ satisfies $IA_nI=A_n^{-1}$ and is diagonalizable with at most
		countably many eigenvalues. In particular, $U_t^{(n)}:=A_n^{it}$ defines
		an almost periodic strongly continuous orthogonal representation on $H_\R$.
		
		\item
		For every $t\in\R$, we have
		\[
		\|A_n^{it}-A^{it}\|\longrightarrow0
		\qquad\text{as }n\to\infty.
		\]
		
		\item
		Writing
		$j:=\left(\frac{2}{1+A^{-1}}\right)^{1/2}$ and 
		$j_n:=\left(\frac{2}{1+A_n^{-1}}\right)^{1/2}$,
		we have
		\[
		\|j_n-j\|\longrightarrow0
		\qquad\text{as }n\to\infty.
		\]
		
		\item
		We have
		\[
		T_nj=e^{-1/(2n)}j_n,
		\qquad
		e^{-1/n}1\leq T_n\leq1,
		\qquad
		\|T_n-1\|\longrightarrow0
		\quad\text{as }n\to\infty.
		\]
		
		\item
		Both $A_n$ and $T_n$ are Borel functions of $A$.
		In particular, every bounded operator commuting with all $A^{it}$
		also commutes with all $A_n^{it}$ and with $T_n$.
	\end{enumerate}
\end{lem}

\begin{proof}
	Put $B=\log A$, so that $IBI=-B$.
	Choose an odd, countably valued Borel function $h_n:\R\to\R$ such that
	\[
	|h_n(x)-x|\leq\frac1n
	\qquad(x\in\R).
	\]
	For instance, set
	\[
	h_n(x)=\frac{\lfloor nx\rfloor}{n}
	\qquad(x\geq0)
	\]
	and extend $h_n$ to an odd function on $\R$.
	Set
	\[
	B_n=h_n(B),
	\qquad
	A_n=e^{B_n}.
	\]
	These operators need not be bounded.
	Since $h_n$ has countable range, $A_n$ is diagonalizable with at most
	countably many eigenvalues.
	Moreover, functional calculus gives
	$IB_nI=-B_n$ and $IA_nI=A_n^{-1}$.
	For every $t\in\R$,
	\[
	\|A_n^{it}-A^{it}\|
	\leq
	\sup_{x\in\R}|e^{ith_n(x)}-e^{itx}|
	\leq\frac{|t|}{n}\longrightarrow0,\quad \text{as }n\to\infty .
	\]
	Define
	\[
	r(x)=\left(\frac{2}{1+e^{-x}}\right)^{1/2},
	\]
	so that $j=r(B)$ and $j_n=r(B_n)$.
	Since
	$(\log r)'(x)=\frac{1}{2(1+e^x)}$,
	the function $\log r$ is $1/2$-Lipschitz, and therefore
	\[
	e^{-1/(2n)}
	\leq\frac{r(h_n(x))}{r(x)}
	\leq e^{1/(2n)}
	\qquad(x\in\R).
	\]
	As $r(x)\leq\sqrt2$, this implies
	\[
	\|j_n-j\|
	\leq\sqrt2\bigl(e^{1/(2n)}-1\bigr)
	\longrightarrow0,\quad \text{as }n\to \infty.
	\]
	Now put
	\[
	t_n(x)=e^{-1/(2n)}\frac{r(h_n(x))}{r(x)},
	\qquad
	T_n=t_n(B).
	\]
	The preceding bounds give
	\[
	e^{-1/n}1\leq T_n\leq1,
	\qquad
	\|T_n-1\|\leq1-e^{-1/n},
	\]
	while the definition gives $T_nj=e^{-1/(2n)}j_n$.
	Finally, $A_n$ and $T_n$ are Borel functions of $A$, which proves the
	last assertion.
\end{proof}

\begin{prop}\label{prop-almost-periodic-approximation-qaw}
	Let $(H_\R,U_t)$ be a strongly continuous
	orthogonal representation.
	There exist almost periodic strongly continuous orthogonal representations
	$(H_\R,U_t^{(n)})$ and vacuum-state-preserving normal u.c.p.\ maps
	\[
	\theta_n:M:=\Gamma_q(H_\R,U_t)''
	\longrightarrow
	M^{(n)}:=\Gamma_q(H_\R,U_t^{(n)})'',
	\qquad
	\rho_n:M^{(n)}\longrightarrow M
	\]
	with the following properties.
	
	Using the associated embeddings $j$ and $j_n$, realize all these algebras
	on $\mathcal F_q(H)$, where $H=H_\R\otimes_\R\C$.
	Denote their Wick words by $W$ and $W_n$, respectively, and put
	\[
	A_q=C^*\{1,W(j(\xi))\mid\xi\in H_\R\},
	\qquad
	A_q^{(n)}=C^*\{1,W_n(j_n(\xi))\mid\xi\in H_\R\}.
	\]
	With these realizations, the modular conjugations associated with the
	vacuum states are given by the same operator $J_q$ on $\mathcal F_q(H)$.
	\begin{enumerate}
		\item
		We have $\theta_n(A_q)\subset A_q^{(n)}$ and
		$\rho_n(A_q^{(n)})\subset A_q$.
		More precisely, with $c_n:=e^{-1/(2n)}$, we have
		\[
		\begin{aligned}
			\theta_n\bigl(W(j\xi_1\otimes\cdots\otimes j\xi_r)\bigr)
			&=
			c_n^rW_n(j_n\xi_1\otimes\cdots\otimes j_n\xi_r),\\
			\rho_n\bigl(W_n(j_n\xi_1\otimes\cdots\otimes j_n\xi_r)\bigr)
			&=
			c_n^rW(j\xi_1\otimes\cdots\otimes j\xi_r)
		\end{aligned}
		\]
		for all $r\geq1$ and $\xi_1,\ldots,\xi_r\in H$.

		\item
		In the common Fock representation,
		$\theta_n(x)\to x$ $*$-strongly as $n\to\infty$ for every $x\in M$, and
		$\|\theta_n(a)-a\|\to0$ as $n\to\infty$ for every $a\in A_q$.
		Moreover,
		\[
		\rho_n\circ\theta_n=\mathcal P_{1/n},
		\qquad
		\theta_n\circ\rho_n=\mathcal P_{1/n}^{(n)},
		\]
		where, for $t\geq0$, $\mathcal P_t$ and $\mathcal P_t^{(n)}$ are the radial
		u.c.p.\ maps on $M$ and $M^{(n)}$, respectively, multiplying
		Wick words of length $r$ by $e^{-tr}$.
		In particular, $\{\rho_n\circ\theta_n\}_n$ converges to the identity
		pointwise in norm on $A_q$ and pointwise $*$-strongly on $M$.
		
		\item
		Let $G$ be a locally compact group, and let
		$\pi:G\to\mathcal O(H_\R)$ be a strongly continuous orthogonal
		representation commuting with $U$.
		Then $\pi$ commutes with every $U^{(n)}$.
		Writing $\sigma^\pi$ and $\sigma^{\pi,(n)}$ for the induced
		actions on $M$ and $M^{(n)}$, we have
		\[
		\theta_n\circ\sigma_g^\pi
		=
		\sigma_g^{\pi,(n)}\circ\theta_n,
		\qquad
		\rho_n\circ\sigma_g^{\pi,(n)}
		=
		\sigma_g^\pi\circ\rho_n
		\qquad(g\in G).
		\]
	\end{enumerate}
\end{prop}

\begin{proof}
	Write $U_t=A^{it}$ and let $I$ be the canonical conjugation on $H$.
	Choose $A_n,j_n,T_n$ as in
	Lemma~\ref{lem-almost-periodic-approximation-of-generator}, and put
	\[
	U_t^{(n)}=A_n^{it},
	\qquad
	c_n=e^{-1/(2n)},
	\qquad
	\widehat T_n=c_n^2T_n^{-1}.
	\]
	Then
	\[
	T_nj=c_nj_n,
	\qquad
	\widehat T_nj_n=c_nj,
	\qquad
	c_n^2\,1\leq T_n,\widehat T_n\leq1.
	\]
	In particular, both $T_n$ and $\widehat T_n$ are contractions converging
	to $1$ in norm.
	
	Let $S$ and $S_n$ be the closed antilinear involutions on $H$ given by
	\[
	\begin{aligned}
		&\operatorname{dom}(S)=jH,
		&&S(j\xi)=jI\xi,\\
		&\operatorname{dom}(S_n)=j_nH,
		&&S_n(j_n\xi)=j_nI\xi
	\end{aligned}
	\qquad(\xi\in H).
	\]
	The preceding identities show that
	\[
	T_n(jH)=j_nH,
	\qquad
	\widehat T_n(j_nH)=jH,
	\]
	and that
	\[
	S_nT_n=T_nS\quad\text{on }jH,
	\qquad
	S\widehat T_n=\widehat T_nS_n\quad\text{on }j_nH.
	\]
	Thus, by the second quantization theorem
	\cite[Theorem~3.4]{Wa17}
	(see also \cite[Corollary~3.16]{HR10} for $q=0$),
	these contractions induce vacuum-state-preserving normal u.c.p.\ maps
	\[
	\theta_n:M\longrightarrow M^{(n)},
	\qquad
	\rho_n:M^{(n)}\longrightarrow M.
	\]
	Their action on Wick words is obtained by applying $T_n$ and
	$\widehat T_n$, respectively, to each tensor factor.
	Their GNS implementations are
	$\mathcal F_q(T_n)$ and $\mathcal F_q(\widehat T_n)$.
	
	(1) 
	The asserted Wick-word formulas follow from
	$T_nj=c_nj_n$ and $\widehat T_nj_n=c_nj$.
	They also give the inclusions of the corresponding $C^*$-algebras,
	since Wick polynomials are norm dense.
	
	(2) 
	We prove convergence in the common Fock representation.
	For $\xi_1,\ldots,\xi_r\in H_\R$, the Wick formula expresses
	$W_n(j_n\xi_1\otimes\cdots\otimes j_n\xi_r)$ 
	as a finite sum of products of $\ell_q(j_n\xi_i)$ and their adjoints,
	with coefficients independent of $n$.
	Since $\|j_n-j\|\to0$ and $c_n\to1$, norm continuity of the creation
	operators gives $\|\theta_n(a)-a\|\to0$ for every Wick polynomial $a$.
	By contractivity of $\theta_n$, this holds for every $a\in A_q$.
	
	Since $T_n\to1$ in norm, we have $\mathcal F_q(T_n)\to1$ strongly.
	Consequently,
	\[
	\theta_n(x)\Omega
	=
	\mathcal F_q(T_n)x\Omega
	\longrightarrow x\Omega,\quad\text{as }n\to \infty
	\qquad(x\in M).
	\]
	To pass from convergence on $\Omega$ to strong convergence, note that
	$M$ and all $M^{(n)}$ have the same modular conjugation $J_q$
	with respect to $\Omega$, given by
	\[
	J_q\Omega=\Omega,
	\qquad
	J_q(\eta_1\otimes\cdots\otimes\eta_r)
	=
	I\eta_r\otimes\cdots\otimes I\eta_1.
	\]
	For $b\in A_q$, put $b_n=\theta_n(b)$.
	Then $\|b_n-b\|\to0$ and $Jb_nJ\in(M^{(n)})'$.
	For every $x\in M$, it follows that
	\[
	\begin{aligned}
		\|(\theta_n(x)-x)JbJ\Omega\|
		&\leq
		\|\theta_n(x)J(b-b_n)J\Omega\|
		+
		\|Jb_nJ\theta_n(x)\Omega-JbJx\Omega\|\\
		&\leq
		2\|x\|\,\|b_n-b\|
		+
		\|b_n\|\,\|(\theta_n(x)-x)\Omega\|\\
		&\longrightarrow0,\qquad \text{as }n\to \infty.
	\end{aligned}
	\]
	Since $JA_qJ\Omega$ is dense in $\mathcal F_q(H)$ and $\theta_n$ is a contraction,
	we obtain $\theta_n(x)\to x$ strongly.
	Applying the same argument to $x^*$ gives $*$-strong convergence.
	
	Since
	\[
	\widehat T_nT_n=T_n\widehat T_n=c_n^2\,1=e^{-1/n}1,
	\]
	the composition identities
	$\rho_n\circ\theta_n=\mathcal P_{1/n}$ and
	$\theta_n\circ\rho_n=\mathcal P_{1/n}^{(n)}$
	hold on Wick words and hence, by normality, on the whole algebras.
	The asserted convergence of $\rho_n\circ\theta_n=\mathcal P_{1/n}$
	follows by the same argument from its action on Wick words and the strong
	convergence of its GNS implementation
	$\mathcal F_q(e^{-1/n}1)$.
	
	(3) If $\pi$ commutes with $U$, then
	Lemma~\ref{lem-almost-periodic-approximation-of-generator}
	implies that $\pi$ commutes with every $U^{(n)}$, as well as with
	$T_n$ and $\widehat T_n$.
	The equivariance identities follow on Wick words and extend by
	normality.
\end{proof}

\section{Intertwining between orthogonal summands in $q$-Araki--Woods factors}
\label{Intertwining between orthogonal summands in $q$-Araki--Woods factors}

Let $G$ be a second countable locally compact group. For $i=1,2$, let $U^{i}\colon \R\curvearrowright H_\R^{i}$ and $\pi^{i}\colon G\curvearrowright H_\R^{i}$ be strongly continuous orthogonal representations satisfying $[U^i_t,\pi^i_g]=0$ for all $t\in \R$ and $g\in G$. We set
\[
H_\R=H_\R^{1}\oplus H_\R^{2},
\qquad
U_t=U_t^{1}\oplus U_t^{2},
\qquad
\pi_g=\pi_g^{1}\oplus\pi_g^{2}
\]
and denote the corresponding $q$--Araki--Woods algebras and crossed products by
\[
M_q=\Gamma_q(H_\R,U_t)'',
\qquad
M_q(G)=M_q\rtimes_{\sigma^\pi}G,
\]
\[
M_q^i=\Gamma_q(H_\R^i,U^i_t)'',
\qquad
M_q^i(G)=M_q^i\rtimes_{\sigma^{\pi^i}}G,
\qquad i=1,2.
\]
The purpose of this appendix is to prove the following proposition, which may be regarded as a generalization of \cite[Proposition~4.10]{CSW17}. Unlike in the setting considered there, we deal with arbitrary finite subalgebras with expectation associated with possibly infinite-dimensional orthogonal summands. Consequently, the finite-dimensionality argument used in \cite[Proposition~4.10]{CSW17} is no longer directly available, and we need to keep track of the relevant Wick-word contractions more explicitly.

\begin{prop}\label{prop-two-components-intertwining}
	Let $p\in M_q^1(G)$ be a projection and $A\subset pM_q^1(G)p$ a finite von Neumann subalgebra with expectation.
	If $A\preceq_{M_q(G)}M_q^2(G)$,
	then $A\preceq_{M_q^1(G)}LG$.
\end{prop}

Put $H=H_\R\otimes_\R\C$ and
$H^i=H_\R^i\otimes_\R\C$ for $i=1,2$.
Let $H_{\rm an}$ be the space of vectors analytic for $U$, and set
$H_{\rm an}^i:=H^i\cap H_{\rm an}$.
We denote the vacuum state on $M_q$ by $\varphi$
and use the same notation for its restrictions to $M_q^i$.
Let 
\[
\W(M_q)
:=A_{q,\mathrm{alg}}
=\operatorname{span}
\left\{
W(\xi_1\otimes\cdots\otimes\xi_n)
\mathrel{}\middle|\mathrel{}
n\geq 0,\ \xi_1,\ldots,\xi_n\in jH
\right\}
\]
be the space of Wick polynomials. We denote by $\W_{\rm an}(M_q)\subset \W(M_q)$ the subspace obtained by choosing each $\xi_k$ to be in $H_{\rm an}$. Similarly, $\W_{\rm an}(M_q^i)\subset \W(M_q^i)$ is defined for $i=1,2$. 

We denote the canonical conditional expectations onto
$M_q^i$, $M_q^i(G)$, and $LG$ by
$E_{M_q^i}$, $E_{M_q^i(G)}$, and $E_{LG}$, respectively,
using the same notation for their restrictions when appropriate.

To prove Proposition~\ref{prop-two-components-intertwining}, we indeed prove the following lemma.

\begin{lem}\label{lem-finite-wick-sandwich}
	For any $x\in \W(M_q)$ and $y\in \W_{\rm an}(M_q)$, there exist
	\[
	a_\nu\in\W(M_q),
	\quad
	b_\nu\in\W_{\mathrm{an}}(M_q),
	\quad
	a_\nu^{(i)}\in\W(M_q^i),
	\quad
	b_\nu^{(i)} \in\W_{\mathrm{an}}(M_q^i),\quad i=1,2,
	\]
	for $\nu=1,\ldots,N$ such that, for every $v\in M_q^1(G)$,
	\[
	E_{M_q^2(G)}(xvy)
	=
	\sum_{\nu=1}^N
	a_\nu^{(2)}
	E_{LG}\left(
	a_\nu E_{LG}(a_\nu^{(1)} \, v \, b_\nu^{(1)})b_\nu
	\right)b_\nu^{(2)}.
	\]
\end{lem}

\begin{proof}[{\bf Proof of Lemma~\ref{lem-finite-wick-sandwich} $\Rightarrow$ Proposition~\ref{prop-two-components-intertwining}}]
	Suppose that the conclusion does not hold. By Popa's intertwining criterion, there exists a net
	$(u_i)_i\subset\mathcal U(A)$
	such that
	\[
	E_{LG}(b^*u_i a)\longrightarrow0
	\qquad\text{$\sigma$--strongly},\quad \text{as }i\to \infty
	\]
	for all $a,b\in M^1_q(G)$.
	Applying this to the right-hand side of the conclusion of Lemma~\ref{lem-finite-wick-sandwich}, and denoting by $\lambda_g$ $(g\in G)$ the group unitaries, we obtain
	\[
	E_{M_q^2(G)}
	\bigl((y\lambda_t)^*u_i(x\lambda_s)\bigr)
	=
	\lambda_t^* E_{M_q^2(G)}
	\bigl(y^*u_ix\bigr)\lambda_s
	\longrightarrow0
	\qquad\text{$\sigma$--strongly},\quad \text{as }i\to \infty
	\]
	for any
	$y\in\W(M_q)$, $x\in\W_{\mathrm{an}}(M_q)$, 
	$s,t\in G$.
	Since the elements $y\lambda_t$ and $x\lambda_s$ generate dense subspaces, it follows that $A\not\preceq_{M_q(G)}M_q^2(G)$, contradicting the assumption.
\end{proof}

It remains to prove Lemma~\ref{lem-finite-wick-sandwich}.
For this, we need several lemmas. 
For integers $a,b,c\geq 0$, let
\[
R_{a,b,c}
\colon
H_q^{\otimes a}\otimes H_q^{\otimes b}\otimes H_q^{\otimes c}
\longrightarrow
H_q^{\otimes(a+b+c)}
\]
be the operator induced by the identity map on the algebraic tensor product $H^{\otimes_\alg a+b+c}$.
When there are two indices, we similarly write $R_{a,b}$.
Next, write $U_t = A^{\ri t} $ $(t\in \R)$ on $H$, and let $S:=I A^{-1/2}$ be the Tomita operator, where $I$ denotes the canonical involution on $H=H_\R\otimes_{\R}\C$. For each $d\geq 1$, define
\[
S_d(\xi_1\otimes\cdots\otimes\xi_d)
:=
S\xi_d\otimes\cdots\otimes S\xi_1.
\]
Wet $S_0(z\Omega)=\overline z\Omega$
for $z\in\C$.
Set $m^0=\id_{\C}$ and
\[
m^d(\xi\otimes\eta)
:=
\langle S_d\xi,\eta\rangle_q,\quad \xi,\eta\in (jH)^{\otimes_{\alg} d}.
\]
As observed in the proof of \cite[Lemma 3.4]{DT25}, the arguments of
\cite[Theorem 3.3]{EP03} and \cite[Proposition 4.9]{CIW19} yield the following formula.

\begin{lem}\label{lem-Wick-formula-triple}
	Let $n,p,m\geq 0$. Let 
	$\xi\in (jH)^{\otimes_{\alg} n}$, 
	$\rho\in (jH)^{\otimes_{\alg} p}$, 
	$\eta\in (jH)^{\otimes_{\alg} m}$ be simple tensors.
	Then $W(\xi)W(\rho)W(\eta)$ is equal to
	\[
	\sum_{\substack{a,b,c\geq 0\\
			a+b\leq n,\;
			b+c\leq p,\;
			a+c\leq m}}
	q^{a(p-b-c)}
	W\Bigl(
	m_{13}^a m_{12}^b m_{23}^c
	\bigl(
	R_{n-a-b,a,b}^*(\xi)
	\otimes
	R_{b,p-b-c,c}^*(\rho)
	\otimes
	R_{c,a,m-a-c}^*(\eta)
	\bigr)
	\Bigr).
	\]
	Here $m_{13}^a$ denotes the map obtained by applying $m^a$ in the $H_q^{\otimes a}$ directions of $R_{n-a-b,a,b}^*(\xi)$ and $R_{c,a,m-a-c}^*(\eta)$, and $m_{12}^b$ and $m_{23}^c$ are defined similarly.
\end{lem}

\begin{lem}\label{lem-map-product-conjugate}
	The following set is closed under products and complex conjugation:
	\[
	\mathcal C
	:=
	\operatorname{span}
	\left\{
	G\ni g\mapsto
	\varphi\bigl(a^*\sigma_g^\pi(b)\bigr)
	\mathrel{}\middle|\mathrel{}
	a\in\W(M_q),\ 
	b\in\W_{\mathrm{an}}(M_q)
	\right\}.
	\]
	Here products and complex conjugation are understood in the natural sense for functions on $G$.
\end{lem}
\begin{proof}
	For simple tensors $\xi \in (jH)^{\otimes_{\alg} n}$ and $\eta\in (H_{\rm an})^{\otimes_\alg k}$, we have
	\[
	\varphi\bigl(W(\xi)^*\sigma_g^\pi(W(\eta))\bigr)
	=
	\delta_{n,k}\langle \xi, \pi^{\otimes n}_g\eta\rangle_q.
	\]
	It follows that
	\[
	\mathcal C
	=
	\operatorname{span}
	\left\{
	G\ni g\mapsto
	\langle \xi, \pi^{\otimes n}_g\eta\rangle_0
	\mathrel{}\middle|\mathrel{}n\geq 0,\ 
	\xi \in (jH)^{\otimes_{\alg} n}, \ \eta\in (H_{\rm an})^{\otimes_\alg n}
	\right\}.
	\]
	Indeed, the $q$-inner product can be expressed in terms of the $0$-inner product and the invertible positive operator $P_q^{(n)}=\sum_{\sigma\in S_n}q^{\operatorname{inv}(\sigma)} U_\sigma\in\B(H^{\otimes n})$. 
	Here $U_\sigma$ is the permutation operator defined by
	\[
	U_\sigma(\xi_1\otimes\cdots\otimes\xi_n)
	=
	\xi_{\sigma(1)}\otimes\cdots\otimes\xi_{\sigma(n)}.
	\]
	Moreover, $(P_q^{(n)})^{-1}$ belongs to the finite-dimensional $C^*$-algebra generated by the $U_\sigma$, and hence is a finite linear combination of the $U_\sigma$, from which the assertion follows.
	Using this expression in terms of the $0$-inner product, it is immediate that $\mathcal C$ is closed under products.
	
	To show that $\mathcal C$ is closed under complex conjugation, we use the KMS condition. Indeed, for $a\in\W(M_q)$ and $ b\in\W_{\mathrm{an}}(M_q)$,
	\[
	\overline{\varphi\bigl(a^*\sigma_g^\pi(b)\bigr)}
	=
	\varphi\bigl(\sigma_g^\pi(b^*)a\bigr)
	=
	\varphi\left(
	a\,\sigma_{-i}^\varphi
	\bigl(\sigma_g^\pi(b^*)\bigr)
	\right).
	\]
	Since $\sigma_{-i}^\varphi
	\bigl(\sigma_g^\pi(b^*)\bigr)=\sigma_g^\pi
	\bigl(\sigma_{-i}^\varphi(b^*)\bigr)$, and since $a^*\in\W(M_q)$ and
	$\sigma_{-i}^\varphi(b^*)\in\W_{\rm an}(M_q)$,
	the last expression belongs to $\mathcal C$.
\end{proof}

\begin{lem}\label{lem-Wick-word-group}
	For any simple tensors $\xi\in (jH)^{\otimes_{\alg} n}$ and $\eta \in (H_{\rm an})^{\otimes_\alg k}$, there exist finitely many
	\[
	f_\mu\in\mathcal C,
	\quad
	\xi_\mu \in (jH)^{\otimes_{\alg} n_\mu}
	\quad
	(n_\mu\leq n),
	\quad
	\eta_\mu \in (H_{\rm an})^{\otimes_\alg k_\mu}
	\quad
	(k_\mu\leq k)
	\]
	such that
	\[
	W\bigl(\xi\otimes \pi_g^{\otimes k}\eta\bigr)
	=
	\sum_\mu
	f_\mu(g)\,
	W(\xi_\mu)\sigma_g^{\pi}(W(\eta_\mu))\quad \text{for all }g\in G.
	\]
\end{lem}
\begin{proof}
	We may assume that $\xi$ and $\eta$ are elementary tensors. We prove the assertion by induction on $n+k$. The cases $n=0$ or $k=0$ are immediate. Suppose that $n,k\geq 1$.
	Recall that the Wick formula gives (see, e.g., the proof of \cite[Theorem 4.3]{ABW17})
	\[
	W(\zeta_1\otimes\cdots\otimes\zeta_N)
	=
	W(\zeta_1)W(\zeta_2\otimes\cdots\otimes\zeta_N)
	-
	\sum_{r=2}^N
	q^{r-2}
	\langle S\zeta_1,\zeta_r\rangle_q
	W(\zeta_2\otimes\cdots\otimes\widehat{\zeta_r}
	\otimes\cdots\otimes\zeta_N)
	\]
	for any $\zeta_k\in j(H)$. It is then straightforward to show that
	\[
	W(\xi)\sigma_g^\pi(W(\eta))
	=
	W(\xi)W\bigl(\pi_g^{\otimes k}\eta\bigr)
	=
	W\bigl(\xi\otimes\pi_g^{\otimes k}\eta\bigr)
	+
	\sum_{\mu}
	f_\mu(g)
	W\bigl(
	\xi_\mu\otimes
	\pi_g^{\otimes k_\mu}\eta_\mu
	\bigr)
	\]
	for finitely many $f_\mu\in \mathcal C$, 
	$\xi_\mu\in (jH)^{\otimes_{\alg} n_\mu}$, and 
	$\eta_\mu\in(H_{\rm an})^{\otimes_\alg k_\mu}$, with 
	$n_\mu\leq n$, $k_\mu\leq k$, and 
	$n_\mu+k_\mu<n+k$.
	We apply the induction hypothesis to each
	$W\bigl(
	\xi_\mu\otimes
	\pi_g^{\otimes k_\mu}\eta_\mu
	\bigr)$. Since $\mathcal C$ is closed under products, this equality yields the conclusion.
\end{proof}

We now prove Lemma~\ref{lem-finite-wick-sandwich}.

\begin{proof}[{\bf Proof of Lemma~\ref{lem-finite-wick-sandwich}}]
	We may assume that $x,y$ are of the following form: $x=W(\xi)$, $y=W(\eta)$ and
	\[
	\xi=\xi_1\otimes\cdots\otimes\xi_n\in (jH)^{\otimes_{\alg} n},
	\qquad 
	\eta=\eta_1\otimes\cdots\otimes\eta_m\in H_{\rm an}^{\otimes_\alg m},
	\]
	where each $\xi_i,\eta_j$ belongs to either $H^1$ or $H^2$. 
	We first fix $g\in G$ and $p\geq 0$ and determine the form of the linear operator defined by
	\[
	\Phi_{p,g}\colon (jH^1)^{\otimes_\alg p}\ni \rho \mapsto E_{M_q^2}\bigl(xW(\rho) \sigma_g^\pi(y)\bigr)=E_{M_q^2}\bigl(W(\xi)W(\rho)W(\pi^{\otimes m}_g\eta)\bigr).
	\]
	We use the set $\mathcal C$ in Lemma \ref{lem-map-product-conjugate}.
	
	\begin{claim}
		There exist 
		\[
		h_{\nu,p}\in\mathcal C,
		\quad 
		\beta_{\nu,p}^{(1)}\in
		(jH^1)^{\otimes_\alg b_{\nu,p}},
		\quad
		\alpha_{\nu,p}^{(1)}\in
		(H_{\rm an}^1)^{\otimes_\alg c_{\nu,p}},
		\quad
		b_{\nu,p}+c_{\nu,p}=p,
		\]
		\[
		\beta_{\nu,p}^{(2)}\in
		(jH^2)^{\otimes_\alg k_{\nu,p}},
		\quad
		\alpha_{\nu,p}^{(2)}\in
		(H_{\rm an}^2)^{\otimes_\alg \ell_{\nu,p}}\quad
		k_{\nu,p}+\ell_{\nu,p}\leq n+m,
		\]
		for $\nu=1,\ldots,N_p$ 
		such that, for every $g\in G$ and every simple tensor $\rho\in (jH^1)^{\otimes_{\alg} p}$,
		\[
		\Phi_{p,g}(\rho)
		=
		\sum_{\nu=1}^{N_p}
		h_{\nu,p}(g)
		\varphi\left(
		W(\beta_{\nu,p}^{(1)}\otimes
		\pi_g^{\otimes c_{\nu,p}}\alpha_{\nu,p}^{(1)})^*
		W(\rho)
		\right)
		W\left(
		\beta_{\nu,p}^{(2)}\otimes
		\pi_g^{\otimes \ell_{\nu,p}}\alpha_{\nu,p}^{(2)}
		\right).
		\]
		Furthermore, if $p>n+m$, then $\Phi_p=0$.
	\end{claim}
	\begin{proof}[{\bf Proof of Claim}]
		Let (A) denote the formula obtained by applying Lemma~\ref{lem-Wick-formula-triple} to $W(\xi)W(\rho)W(\pi^{\otimes m}_g\eta)$. Applying $E_{M_q^2}$ to (A), since $\rho$ has only $H^1$ components, all terms that do not satisfy
		\[
		p-b-c=0,\quad R_{b,p-b-c,c}^*(\rho)=R_{b,0,c}^*(\rho)=R_{b,c}^*(\rho)
		\]
		vanish. Hence, in what follows, it suffices to consider only the case $b+c=p$. In particular, we have $\Phi_{p,g}=0$ if $p>n+m$.
		
		Fix $a,b,c$. We expand the following as finite sums of elementary tensors:
		\[
		R_{n-a-b,a,b}^*(\xi)
		=
		\sum_\alpha \xi_0\otimes\xi_{a}\otimes\xi_{b},\quad 
		R_{c,a,m-a-c}^*(\eta)
		=
		\sum_\beta\eta_{c}\otimes\eta_{a}\otimes\eta_0.
		\]
		Since $R^*$ is given by a finite sum of weighted shuffles and $\eta$ is made of analytic vectors, $\eta_c,\eta_a,\eta_0$ are elementary tensors of analytic vectors. Moreover,
		\[
		R_{c,a,m-a-c}^*
		\bigl(\pi_g^{\otimes m}\eta\bigr)
		=
		\sum_\beta \pi_g^{\otimes c}\eta_{c}
		\otimes
		\pi_g^{\otimes a}\eta_{a}
		\otimes
		\pi_g^{\otimes(m-a-c)}\eta_0
		\]	
		holds. Using these expressions, we compute the scalar coefficients arising from the application of $m^a_{13}m_{12}^{b}m_{23}^{c}$ in (A).
		First, applying $m^a_{13}$ gives
		\[
		h_{a}(g)
		:=
		\left\langle
		S_a\xi_{a},
		\pi_g^{\otimes a}\eta_{a}
		\right\rangle_q,\quad h_a \in \mathcal C.
		\]
		Next, applying $m_{12}^{b}m_{23}^{c}$ gives
		\[
		m_{12}^b m_{23}^c
		\left(
		\xi_{b}
		\otimes
		R_{b,c}^*(\rho)
		\otimes
		\pi_g^{\otimes c}\eta_{c}
		\right)
		=
		\varphi\left(
		W(S_b\xi_{b}\otimes\pi_g^{\otimes c}S_c^*\eta_{c})^*
		W(\rho)
		\right).
		\]
		This equality can be verified by writing $R_{b,c}^*(\rho)$ as a finite sum of elementary tensors. Since $\eta_c$ is analytic, it belongs to the domain of $S_c^*$ and
		$S_c^*\eta_c$ is again analytic.
		We may restrict ourselves to the case where this value is nonzero. Since $\rho\in (H^1)^{\otimes p}$, it is enough to consider the case
		$S_b\xi_{b}\otimes\pi_g^{\otimes c}S_c^*\eta_{c}\in (H^1)^{\otimes p}$ and $b+c=p$.
		
		Thus, we have obtained
		\[
		\Phi_{p,g}(\rho)
		=
		\sum_{a,b,c}\sum_{\alpha,\beta} h_{a}(g)\varphi\left(
		W(S_b\xi_{b}\otimes\pi_g^{\otimes c}S_c^*\eta_{c})^*
		W(\rho)
		\right)
		E_{M^2_q}\left(W(\xi_0\otimes \pi_g^{\otimes(m-a-c)}\eta_0)\right).
		\]
		As we are applying $E_{M^2_q}$, we may assume that $\xi_0,\eta_0$ are made only of vectors from $H^2$, and this completes the proof.
	\end{proof}
	
	\begin{claim}
		There exist finitely many
		\[
		a_\nu\in\W(M_q),
		\quad
		b_\nu\in\W_{\mathrm{an}}(M_q),
		\quad
		a_\nu^{(i)}\in\W(M_q^i),
		\quad
		b_\nu^{(i)} \in\W_{\mathrm{an}}(M_q^i),\quad i=1,2,
		\]
		such that, for every $g\in G$ and $z\in \W(M_q^1)$,
		\[
		E_{M_q^2}\bigl(W(\xi)z W(\pi^{\otimes m}_g\eta)\bigr)
		=
		\sum_\nu
		\varphi\left(
		a_\nu\sigma_g^\pi(b_\nu)
		\right)
		\varphi\left(
		a_\nu^{(1)}
		z\sigma_g^{\pi}(b_\nu^{(1)})
		\right)
		a_\nu^{(2)}
		\sigma_g^{\pi}(b_\nu^{(2)}).
		\]
	\end{claim}
	\begin{proof}[{\bf Proof of Claim}]
		For each $p\leq n+m$, take
		$h_{\nu,p},
		\beta_{\nu,p}^{(i)},
		\alpha_{\nu,p}^{(i)}$
		as in the previous claim. 
		We define a linear map $\Psi_{p,g}\colon \W(M_q^1)\to M_q^2$ by replacing
		$W(\rho)$ on the right-hand side of the preceding formula in the last claim with
		$z\in\W(M_q^1)$. 
		Put $\Psi_{p,g}=0$ if $p>n+m$.		

		Let $z\in \W(M_q^1)$ and write $z=\sum_{p=0}^k z_p$,  
		where each $z_p$ is a homogeneous Wick polynomial of degree $p$. Since $b_{\nu,p}+c_{\nu,p}=p$, we have $\Psi_{p,g}(z) = \Psi_{p,g}(z_p)$. Hence,
		\[
		\begin{aligned}
			&E_{M_q^2}\bigl(W(\xi)zW(\pi^{\otimes m}_g\eta)\bigr)
			\\
			&\quad=
			\sum_{p=0}^{k} E_{M_q^2}\bigl(W(\xi)z_p W(\pi^{\otimes m}_g\eta)\bigr)
			=\sum_{p=0}^{k} \Psi_{p,g}(z_p)
			=\sum_{p=0}^{n+m} \Psi_{p,g}(z)	\\
			&\quad=
			\sum_{p=0}^{n+m}\sum_{\nu=1}^{N_p}
			h_{\nu,p}(g)
			\varphi\left(
			W(\beta_{\nu,p}^{(1)}\otimes
			\pi_g^{\otimes c_{\nu,p}}\alpha_{\nu,p}^{(1)})^*
			z
			\right)
			W\left(
			\beta_{\nu,p}^{(2)}\otimes
			\pi_g^{\otimes \ell_{\nu,p}}\alpha_{\nu,p}^{(2)}
			\right).
		\end{aligned}
		\]
		Thus, we have determined the form of this linear operator on $\W(M_q^1)$.
		
		Now, for each
		$W\left(
		\beta_{\nu,p}^{(1)}\otimes
		\pi_g^{\otimes c_{\nu,p}}\alpha_{\nu,p}^{(1)}
		\right)$
		on the right-hand side, applying Lemma~\ref{lem-Wick-word-group} on $H^1$ gives
		\[
		W\left(
		\beta_{\nu,p}^{(1)}\otimes
		\pi_g^{\otimes c_{\nu,p}}\alpha_{\nu,p}^{(1)}
		\right)
		=
		\sum_\mu
		f_\mu^{(1)}(g)\,
		W(\xi_\mu^{(1)})\sigma_g^{\pi}(W(\eta_\mu^{(1)})).
		\]
		Since each $\eta_\mu^{(1)}$ is analytic, the KMS condition shows that one can choose $\widetilde{\xi}^{(1)}_\mu,\widetilde{\eta}^{(1)}_\mu$ so that, for every $z\in \W(M_q^1)$ and $g\in G$, 
		\begin{align*}
			\varphi\left(
			W(\beta_{\nu,p}^{(1)}\otimes
			\pi_g^{\otimes c_{\nu,p}}\alpha_{\nu,p}^{(1)})^*
			z
			\right)
			&=
			\sum_\mu
			\overline{f^{(1)}_{\mu}(g)}
			\varphi\left(
			\sigma_g^{\pi}(W(\eta^{(1)}_\mu))^*
			W(\xi^{(1)}_\mu)^*
			z\right)\\
			&=
			\sum_\mu
			\overline{f^{(1)}_{\mu}(g)}
			\varphi\left(
			W( \widetilde{\xi}^{(1)}_\mu)
			z
			\sigma_g^{\pi}(W(\widetilde{\eta}^{(1)}_\mu))
			\right),
		\end{align*}
		where $\widetilde{\eta}^{(1)}_\mu$ is again analytic. Moreover, by Lemma~\ref{lem-map-product-conjugate}, we have
		$\overline{f^{(1)}_{\mu}}\in\mathcal C$.

		Finally, we apply Lemma~\ref{lem-Wick-word-group} on $H^2$ to each
		$W\left(
		\beta_{\nu,p}^{(2)}\otimes
		\pi_g^{\otimes \ell_{\nu,p}}\alpha_{\nu,p}^{(2)}
		\right)$.
		Since $\mathcal C$ is closed under products by Lemma~\ref{lem-map-product-conjugate}, we obtain the conclusion.
	\end{proof}

	It remains to rewrite the formula for crossed products. Take
	$a_\nu,b_\nu,a_\nu^{(i)},b_\nu^{(i)}$
	as in the last claim.
	For any $v=z\lambda_g$, $z\in \W(M_q^1)$, $g\in G$, it is immediate that
	\[
	a_\nu^{(2)}
	E_{LG}\left(
	a_\nu
	E_{LG}\left(
	a_\nu^{(1)}vb_\nu^{(1)}
	\right)b_\nu
	\right)b_\nu^{(2)}
	=
	\varphi\left(
	a_\nu^{(1)}
	z\sigma_g^{\pi}(b_\nu^{(1)})
	\right)
	\varphi\left(
	a_\nu\sigma_g^\pi(b_\nu)
	\right)
	a_\nu^{(2)}
	\sigma_g^{\pi}(b_\nu^{(2)})
	\lambda_g.
	\]
	Summing over $\nu$, the equality in the last claim gives
	\[
	E_{M_q^2}\bigl(W(\xi)zW(\pi^{\otimes m}_g\eta)
	\bigr)\lambda_g 
	= E_{M_q^2(G)}(xvy).
	\]
	This proves the conclusion for $v=z\lambda_g$.
	Since the linear span of elements of this form is $\sigma$--weakly dense in $M_q^1(G)$, the conclusion follows.
\end{proof}

\small{
	
}

\end{document}